\documentclass[twoside,a4paper]{article} 
\usepackage{amssymb}
\usepackage{a4wide} 
\usepackage{latexsym} 
\usepackage{amsmath} 
\usepackage{theorem} 
\usepackage{mathrsfs} 
\usepackage{xspace}  
\usepackage{fancyhdr} 
\usepackage{float} 
\usepackage{listings}
\usepackage{enumerate}
\usepackage{morefloats}
\usepackage{tikz}
\usepackage{accents}
\usepackage{graphicx} 
\usepackage{caption}
\usepackage{rotating}
\usepackage{subcaption}
\definecolor{SussexFlint}{rgb}{.00,.19,.21}
\definecolor{SussexGrey}{rgb}{.51,.58,.49}
\definecolor{SussexOrange}{rgb}{.94,.29,.00}
\definecolor{SussexYellow}{rgb}{1.00,.73,.00}
\definecolor{SussexRed}{rgb}{.94,.01,.49}
\definecolor{SussexPurple}{rgb}{.48,.06,.44}
\definecolor{SussexGreen}{rgb}{.00,.58,.46}
\definecolor{SussexBlue}{rgb}{.00,.58,.65}
\colorlet{a}{SussexOrange}
\colorlet{b}{SussexRed}
\colorlet{c}{SussexYellow}
\colorlet{d}{SussexPurple}
\colorlet{e}{SussexGreen}
\colorlet{f}{SussexBlue}
\colorlet{g}{SussexGrey}
\colorlet{h}{white}
\colorlet{i}{black}
\colorlet{j}{SussexFlint}
\numberwithin{equation}{section}
\newtheorem{theorem}{Theorem}[section]
\newtheorem{lemma}[theorem]{Lemma}
\newtheorem{corollary}[theorem]{Corollary}
\newtheorem{definition}[theorem]{Definition}%
\newtheorem{remark}[theorem]{Remark}
\newcommand{\qed}{\hfill$\square$}

\newcommand{\reals}{\ensuremath{\mathbb{R}}\xspace}
\newcommand{\R}[1]{\ensuremath{\reals^{#1}}\xspace}

\newcommand{\avg}[1]{\ensuremath{\langle\!\langle#1\rangle\!\rangle} }
\newcommand{\jump}[1]{\ensuremath{[\![#1]\!]} }

\newcommand{\Th}{\ensuremath{\mathscr{T}_h}}
\newcommand{\Ta}{\ensuremath{\mathbf{T}}}
\newcommand{\mc}{\ensuremath{\mathcal{H}}}
\newcommand{\Tr}{\ensuremath{\operatorname{Tr}}}

\newcommand{\dgcomp}{\ensuremath{\dg^{\operatorname{comp}}}}

\newcommand{\dgcompo}{\ensuremath{\dgo^{\operatorname{comp}}}}

\newcommand{\nablaTa}{\ensuremath{\nabla_\Ta}}
\newcommand{\Eb}{\ensuremath{\mathscr{E}^b_h}}
\newcommand{\Ei}{\ensuremath{\mathscr{E}^i_h}}
\newcommand{\X}{\ensuremath{\Omega}}
\newcommand{\Y}{\ensuremath{\Upsilon}}
\newcommand{\Eib}{\ensuremath{\mathscr{E}^{i,b}_h}}
\newcommand{\Vb}{\ensuremath{\mathscr{V}^b_h}}
\newcommand{\dgoo}{\ensuremath{V_{h,0}}}
\newcommand{\gradTa}{\ensuremath{\nabla_\Ta}}
\newcommand{\divTa}{\ensuremath{\operatorname{div}_\Ta}}
\newcommand{\betaperp}{\ensuremath{\beta^\perp}}
\newcommand{\dg}{\ensuremath{V_{h,p}}}
\newcommand{\dgo}{\ensuremath{V_{h,p,0}}}
\newcommand{\Aobl}{\ensuremath{{A_h}}}

\newcommand{\Bhalfobl}{\ensuremath{{B_{h,1/2}}}}

\newcommand{\Bhtheobl}{\ensuremath{{B_{h,\theta}}}}

\newcommand{\Jhobl}{\ensuremath{{J_{h}}}}

\newcommand{\Bstarobl}{\ensuremath{{B_{h,*}}}}

\newtheorem{assumption}[theorem]{Assumption}
\title{A DGFEM for Uniformly Elliptic Two Dimensional Oblique boundary-value problems}
\author{Ellya Kawecki\footnote{EK acknowledges support of the Engineering and Physical Sciences Research Council
[EP/L015811/1].}}
\begin{document}
\maketitle
\begin{abstract}
In this paper we present and analyse a discontinuous Galerkin finite element method (DGFEM) for the approximation of solutions to elliptic partial differential equations in nondivergence form, with oblique boundary conditions, on curved domains. In ``E. Kawecki, \emph{A DGFEM for Nondivergence Form Elliptic Equations with Cordes Coefficients on Curved Domains}", the author introduced a DGFEM for the approximation of solutions to elliptic partial differential equations in nondivergence form, with Dirichlet boundary conditions. In this paper, we extend the framework further, allowing for the oblique boundary condition. The method also provides an approximation for the constant occurring in the compatibility condition for the elliptic problems under consideration. 
\end{abstract}
\begin{section}{Introduction}\label{section:1}
The model problem that we consider in this paper is the following oblique boundary-value problem, find $u:\Omega\to\mathbb R$ such that
\begin{equation}\label{1}
\left\{\begin{aligned}
\sum_{i.j=1}^2A_{ij}D^2_{ij}u & = f,\,\,\mbox{in}\,\,\Omega,\\
\beta\cdot\nabla u &\,\,\mbox{is constant},\,\,\mbox{on}\,\,\partial\Omega,
\end{aligned}
\right.
\end{equation}
where $\Omega\subset\mathbb R^2$ is a given, convex $C^2$ domain, $f\in L^2(\Omega)$, and $A\in L^\infty(\Omega;\mathbb R^{2\times 2}_{\operatorname{Sym}})$ satisfies, for some constant $\lambda>0$, 
\begin{equation}\label{A:assumption}
x^TA(\xi)x\ge\lambda|x|^2\,\,\forall x\in\mathbb R^2,\,\,\mbox{for a.e.}\,\,\xi\in\Omega.
\end{equation}
The constant present in the boundary condition of~(\ref{1}) is there to absorb potentially arising compatibility conditions (consider solving the Poisson problem with a homogeneous Neumann boundary condition imposed). Finally, we assume that the vector-valued function $\beta\in C^1(\partial\Omega;\mathbb S^1)$ (hereby called the ``oblique vector").

The oblique boundary-value problem appears in several interesting applications, often dependent upon which dimension, $d$, is considered, and whether or not the oblique boundary-value problem is \emph{strict}. For $d\ge2$, and $\beta\in C^1(\partial\X;\mathbb S^d)$, the boundary-value problem is referred to as \emph{strictly} oblique, if there exits a constant $\delta>0$ such that
\begin{equation}\label{strictobl}
\beta\cdot n_{\partial\X}\ge\delta\quad\mbox{on }\partial\X,
\end{equation}
where $n_{\partial\X}:\partial\X\to\mathbb S^1$ is the unit outward normal to $\partial\X$. If~(\ref{strictobl}) may only hold with $\delta=0$, then the boundary-value problem is called \emph{degenerate} oblique. 

For $d\ge3$,~(\ref{strictobl}) is necessary for the well-posedness of the boundary-value problem, with~\cite{MR1435286} (pg 13--14) providing counter examples to uniqueness for the Poisson problem, in the case that the oblique vector becomes tangential to the boundary, even on a set of zero boundary measure.

That said, the degenerate (or tangential) oblique problem (falling into the class of degenerate elliptic problems), arises naturally in the (geodetic) problem of determining the gravitational fields of celestial bodies~\cite{MR2180076}. This problem was discovered by Poincar\'e~\cite{Poincare-paper} during his work on the theory of tides. In the case that $d=2$, the oblique boundary-value problem arises in systems of conservation laws in~\cite{MR2019399,MR1733695}, where the latter focuses on a mixed elliptic-hyperbolic problem that requires the boundary condition to be strictly oblique. For an overview of the case  $d=2$ one should refer to~\cite{MR871676}, and for the case $d\ge 3$, one should seek~\cite{MR3059278}.  A particular, and broad subclass of the oblique boundary-value problem is the case when $\beta\equiv n_{\partial\X}$, which is in fact the Neumann boundary-value problem.

The author's interest in this type of boundary-value problem, stems from applications to fully nonlinear second order elliptic partial differential equations (PDEs). In particular equations of Monge--Amp\`ere (MA) and Hamilton--Jacobi--Bellman (HJB) type. Upon linearising such equations (for instance by the application of Newton's method), one arrives at an infinite sequence of problems of the form~(\ref{1}), and as such, the linear theory contained in this paper will be applicable when considering these nonlinear problems. The MA problem arises in areas such as optimal transport and differential geometry, and has been an area of interest, both from an analytical and a numerical computation point of view for many years,  see~\cite{MR1624426,MR1426885,MR0180763} and~\cite{MR2916369,MR3162358}; while the HJB problem arises in applications to mean field games, engineering, physics, economics, optimal control, and finance~\cite{MR2179357,lachapelle2011mean}, where~\cite{MR3196952,MR3033005} mark recent developments in the numerical analysis of such problems. 

It is clear that the linearisation of HJB and MA type equations results in a sequence of nondivergence form elliptic equations. What is not immediately clear, is how the oblique boundary condition may also arise. 
In the applications outlined above (geodetic problems, and conservation laws), the oblique \emph{boundary condition} arises, but these problems are not typically cast in nondivergence form.

The nondivergence form oblique boundary-value problem~(\ref{1}) arises in the linearisation of the MA optimal transport problem, which can be posed as follows: given two uniformly convex domains, $\X,\Y\subset\mathbb R^2$, and uniformly positive functions $f_1:\X\to\mathbb R$, $f_2:\Y\to\mathbb R$, find $u:\X\to\mathbb R$ such that
\begin{equation}\label{MA:OT}
\left\{
\begin{aligned}
\operatorname{det}D^2u(x) - \frac{f_1(x)}{f_2(\nabla u(x))}& = 0\quad x\in\X,\\
B_\Upsilon(\nabla u(x)) & = 0\quad x\in\partial\X.
\end{aligned}
\right.
\end{equation}
where $B_\Upsilon:\mathbb R^2\to\mathbb R$ is a convex defining function for $\Y$ (take, for example the signed distance function to $\partial\Y$). For simplicity, we assume that $f_2\equiv c$ is constant. Upon applying Newton's method to~(\ref{MA:OT}), we arrive at the following sequence of inhomogeneous oblique boundary-value problems: given $u_n\in C^2(\overline{\X})$, uniformly convex, find $u_{n+1}$ such that
\begin{equation}\label{newton:MA}
\left\{
\begin{aligned}
\operatorname{Cof}D^2u_n(x):D^2u_{n+1}& = \operatorname{Cof}D^2u_n(x):D^2u_{n}+\frac{f_1(x)}{c}-\operatorname{det}(D^2u_n)\quad x\in\X,\\
DB_\Upsilon(\nabla u_n(x))\!\cdot\!\nabla u_{n+1} & = DB_\Upsilon(\nabla u_n(x))\!\cdot\!\nabla u_{n} -B_\Upsilon(\nabla u_n(x))\quad x\in\partial\X,
\end{aligned}
\right.
\end{equation}
where $\operatorname{Cof}$ denotes the cofactor matrix.
Notice that there is no free constant in the boundary condition for~(\ref{newton:MA}). This is due to the fact that the functions $f_1$, $f_2$ are assumed to satisfy the compatibility condition
$\int_\X f_1 = \int_\Y f_2,$
which is necessary for the existence of a unique (up to a constant) convex solution to the MA optimal transport problem. One can see that the boundary condition in~(\ref{newton:MA}) is of the (inhomogeneous) oblique type, with $\beta := DB_\Upsilon(\nabla u_n(x))/|DB_\Upsilon(\nabla u_n(x))|$. However, when considering the discretisation of~(\ref{newton:MA}), the gradient of an arbitrary finite element function may take arbitrary values in $\mathbb R^2$, but, even when $\partial\Y$ is smooth, the derivative of the defining function, $B_\Upsilon$, may not be well defined everywhere. This can occur when one takes $B_\Upsilon$ to be the signed distance function to $\partial\Y$ (for example, the signed distance function to the unit disk is not smooth at the origin); one may appeal to a different representation of this function in terms of the supporting hyperplanes of $\Y$. That is,
\begin{equation}\label{BCrep}
B(q) := \sup_{\|n\|=1}\{q\cdot n-H^*(n)\},\mbox{ where }H^*(n):=\sup_{q\in\partial\Y}\{q\cdot n\}.
\end{equation}
See~\cite{benamou2012viscosity} for further justification of this representation. In~\cite{benamou2012viscosity}, the authors use this representation of the boundary condition (with further modifications, so that the supremum in~(\ref{BCrep}) is taken over vectors that make an acute angle with the unit normal to $\partial\X$), and iteratively solve (via a finite difference method) the MA problem with \emph{Neumann} boundary conditions. 
Utilising~(\ref{BCrep}), we arrive at the following MA problem:
\begin{equation}\label{MA:OT:HJ}
\left\{
\begin{aligned}
\operatorname{det}D^2u(x) - \frac{f_1(x)}{c}& = 0\quad x\in\X,\\
\sup_{\|n\|=1}\{\nabla u\cdot n-H^*(n)\} & = 0\quad x\in\partial\X.
\end{aligned}
\right.
\end{equation}
This brings us to the intersection of MA and HJB problems. Taking into account the representation~(\ref{BCrep}), one can see that the application of Newton's method in~(\ref{newton:MA}) is not justified, in general, since we may not be able to make sense of the derivative of $B_\Upsilon$. However, a semismooth Newton's method may be justified, which relies only on the subdifferential of the operator under consideration. This type of linearisation scheme is implemented in~\cite{MR3196952} in order to design a convergent scheme for the DG approximation of HJB equations with \emph{Dirichlet} boundary conditions. To move into the context of HJB equations,  one may utilise a longstanding result proven by Krylov~\cite{MR901759}, that allows one to characterise the MA problem~(\ref{MA:OT:HJ}), as the following HJB problem
\begin{equation}\label{MA:OT:HJB}
\left\{
\begin{aligned}
\sup_{W\in X}\left\{-W:D^2u+2(\operatorname{det}W)^{1/2} \frac{f_1^{1/2}(x)}{c^{1/2}}\right\}& = 0\quad\mbox{in }\X,\\
\sup_{\|n\|=1}\{\nabla u\cdot n-H^*(n)\} & = 0\quad\mbox{on }\partial\X,
\end{aligned}
\right.
\end{equation}
where $X:=\{W\in\mathbb R^{2\times 2}:W=W^T,W\ge 0,\operatorname{Trace}(W)=1\}$. Applying a semismooth Newton's method to the operator in the domain and on the boundary, we arrive at a sequence of nondivergence form oblique boundary-value problems.

Exporting the ideas of this paper to the context of~(\ref{MA:OT:HJB}) will be the focus of future work. This will involve recasting the HJB problem in terms of one with uniformly elliptic coefficients (the set $X$ contains degenerate matrices), and, furthermore, a proof of convergence of the semismooth Newton's method with a Hamilton--Jacobi type boundary condition. A key motivation for Krylov's HJB formulation of the MA equation, is that the numerical solution may in fact be unique, which is not always the case when considering numerical methods for the MA equation (see~\cite{MR3049920}, for example).

Another motivation for considering HJB type equation with oblique boundary conditions arises in the area of mean field games, in the modelling of pedestrian crowds (the elliptic counterpart of this problem is to determine the long term behaviour of the crowd). In this case the HJB equation is coupled with a Fokker--Planck equation (see~\cite{lachapelle2011mean}), and the boundary condition models the interaction the crowd will have with the boundary of the space that they inhabit. In the case of oblique boundary conditions, one can think of this interaction as the requirement that the direction in which one is incentivised to exit the domain is determined by the oblique vector.

Now that we have discussed why one is motivated to approximate solutions problems of the form~(\ref{1}), it is pertinent to discuss why standard approaches may not apply in this case.
The problem~(\ref{1}) poses several difficulties, both analytically and numerically. The problem is in nondivergence form (and due to our assumptions, it \emph{cannot} be written in divergence form, in general, and so the standard weak formulation cannot be used here), and as such, well-posedness is not guaranteed (see~\cite{MR1814364,MR2260015,MR1684729}), and the use of standard conforming finite element methods is also not applicable (as they also rely on the weak formulation of the PDE).

 As we will see, on $C^2$ domains, one can show well-posedness by assuming that the coefficients satisfy the ``Cordes condition" (defined in Section~\ref{section:3}), and a condition between the oblique vector, and the curvature of the boundary. In~\cite{MR2260015}, well-posedness is proven using the method of continuity; in this paper we use a different method, analogous to that seen in~\cite{MR3077903} Theorem 3. This technique relies upon a variant of the Miranda--Talenti estimate (see Lemma~\ref{FullMT}). The motivation for the use of a different technique is the extension of this technique to problems of HJB type, for which the method present in~\cite{MR2260015} does not apply.

The Miranda--Talenti estimate is known to hold for $d\ge 2$ when considering functions in a suitable Sobolev space, whose trace or normal derivative vanishes on $\partial\Omega$. However, for $d\ge 3$, this estimate remains an open problem if one assumes that the oblique derivative ($\beta\cdot\nabla u$) is constant on the boundary~\cite{MR1776331}, restricting us to a two-dimensional framework.

In terms of PDE analysis, our goal is to prove existence of a unique strong solution ($H^2$-regular, so that the PDE is satisfied a.e.) $u\in H$ satisfying
 \begin{equation}\label{agammadef}
 A_\gamma(u,v):=\int_\X\gamma A:D^2u\,\Delta v =\int_\X\gamma f\Delta v=:\ell(v)\quad\forall v\in H,
 \end{equation}
 where $H$ is a particular subset of $H^2(\X)$ (see Section~\ref{section:2}), and $\gamma$ is a uniformly positive renormalisation factor defined in Section~\ref{section:3}. This motivates the construction of the finite element method given in Section~\ref{section:4}, with the goal of finding a suitable numerical analogue of the operator $A_\gamma$.
 
 This is another example where the techniques present in this paper divert from the standard concepts of conforming finite element methods. A direct discretisation of the operator $A_\gamma$ in a finite element setting (by replacing the derivatives present in $A_\gamma$ with piecewise derivatives) would fail, since general finite element functions do not satisfy the Miranda--Talenti estimate, and thus the corresponding operator would, in general, fail to be coercive. The alternative would be to construct a finite element space whose members are both $H^2$-regular, and satisfy the oblique boundary condition. In particular, it is not immediately clear how one would do the latter. We instead design a discontinuous Galerkin scheme that numerically enforces the Miranda--Talenti estimates (via internal and boundary jump penalisation terms), and is also consistent (i.e., that the true solution to the PDE is also a solution to the finite element method, assuming the solution has sufficient broken-Sobolev regularity, allowing one to substitute it into the formulation), resulting in optimal a priori error estimates.
 
 It is important to consider the expected regularity of solutions to~(\ref{1}). In general, since the coefficient matrix, $A$ is only $L^\infty(\X)$, and the right-hand side function $f\in L^2(\X)$, standard elliptic regularity theory implies that a solution would, in general, be at most $H^2$-regular. In Section~\ref{section:5} we provide an a priori error estimate for solutions of such regularity. In the case that solutions are smoother, we are able to prove a priori error estimates in a $H^2$-style norm, that are optimal with respect to the polynomial degree of the finite element space. Even though the linearised problems may admit solutions of lower regularity, the limiting problems (by way of Newton's method) may indeed be of higher regularity. For example, in the case that $\partial\X$, $\partial\Y$ are $C^{2,1}$ regular, the MA optimal transport problem~(\ref{MA:OT}) admits a $C^{3,\alpha}(\X)\cap C^{2,\alpha}(\overline{\X})$ convex solution (see~\cite{MR1454261}).


Interestingly, our method also gives an approximation of the constant that arises in the compatibility condition for problems of the form~(\ref{1}). In the case of conormal boundary-value problems, i.e.,
\begin{equation*}
\left\{\begin{aligned}
-\nabla\cdot(A\nabla u) & = f,\,\,\mbox{in}\,\,\Omega,\\
(A\nabla u)\cdot n_{\partial\X} &\,\,\mbox{is constant on}\,\,\partial\Omega,
\end{aligned}
\right.
\end{equation*}
it is quite straightforward (by an application of the Divergence Theorem) to determine the constant in terms of the function $f$. In the problems we consider, however, there does not appear to be any explicit relationship between the constant and the right hand side, in general.

Our approach extends the framework of~\cite{MR3077903} and~\cite{Kawecki:77:article:On-curved} (the first of which applies to the Dirichlet boundary condition on polytopal domains, and the second applies to the Dirichlet boundary condition on curved, uniformly convex domains) to the oblique case. \begingroup\color{black} To the author's knowledge, there is a sparse amount of work on finite element methods for oblique boundary-value problems present in the existing literature, and, as such, this paper provides the first DGFEM for the approximation of solutions to~(\ref{1}). As such, a motivation of this paper is to widen the scope of the current numerical framework for oblique boundary-value problems. 

The papers \cite{MR1331351,favskova2010finite,MR789678} provide examples of finite element approximations of oblique boundary-value problems, where~\cite{favskova2010finite} and~\cite{MR789678} apply to a particular geodetic and free boundary problem respectively. Close to the timing of this paper, the author of~\cite{gallistl2018numerical} introduced a mixed finite element method for the approximation of solutions to~(\ref{1}), and proved a priori and a posteriori error estimates. The hypotheses of~\cite{gallistl2018numerical} with respect to the data for problem~(\ref{1}) are identical to that of this paper (i.e., the assumptions upon $A$, $\partial\X$, $f$, and $\beta$), however, the consideration of curved finite elements provides a difference in our computationally sufficient assumptions. This stems from the fact that the finite element functions that we consider are polynomials on the flat reference simplex, whereas the finite element functions in~\cite{gallistl2018numerical} are polynomials on the curved physical element. In particular, computationally,~\cite{gallistl2018numerical} assumes that $\partial\X\in C^2$, whereas we require that $\partial\X$ is also piecewise $C^3$. However, in~\cite{gallistl2018numerical} the function approximating $\nabla u$ must be piecewise affine, whereas we are able to employ piecewise polynomials of degree $p\ge2$ (as we are approximating $u$, one may consider piecewise quadratics in our case to be the analogue of piecewise affine in the case of~\cite{gallistl2018numerical}).\endgroup
 
This paper is organised as follows: In Section~\ref{section:2} we begin by introducing the notation needed, as well as the function spaces used in the analysis of the PDE. In Section~\ref{section:3}, we prove several important estimates, such as the Miranda--Talenti estimate, and then proceed to prove an existence and uniqueness result for problems of the form~(\ref{1}). In Section~\ref{section:4} we prove an important consistency result, and proceed to prove a stability result for our numerical scheme; this stability result is then used as a main tool in the proof of existence and uniqueness of a numerical solution. In Section~\ref{section:5} we prove an error estimate that is optimal in terms of the mesh size, for piecewise sufficiently smooth solutions, as well as an estimate for solutions of conformal regularity. \begingroup\color{black} In Section~\ref{Experiments:sec:6} we perform several numerical experiments, validating the estimates of Section~\ref{section:5}. \endgroup
Section~\ref{Conclusion} is the final section, where we give concluding remarks on what has been accomplished in this paper.
\end{section}
\begin{section}{Set-up}\label{section:2}
Before we prove existence and uniqueness of solutions to~(\ref{1}), we must first determine a suitable function space in which to do so. In general, we will consider strong solutions of~(\ref{1}). That is, where all weak derivatives of the solution, up to second order, are square integrable, so that~(\ref{1}) may hold a.e. in $\Omega$.
\subsection{Function spaces}
We consider the standard Sobolev space
$$H^2(\Omega):=\{v\in L^2(\Omega):D^\alpha v\in L^2(\Omega),\,\,\forall\,1\le|\alpha|\le 2\}.$$
We define the following subset (determined by the oblique vector $\beta$) of $H^2(\Omega)$:
$$H^2_\beta(\Omega):=\bigcup_{c\in\mathbb R}\{v\in H^2(\Omega):\beta\cdot\nabla v\mbox{ is constant on }\partial\Omega\},$$
as well as a further subspace
$$H^2_{\beta,0}(\Omega):=\left\{v\in H^2_\beta(\Omega):\int_\Omega v = 0\right\}.$$
%
Furthermore, we will endow all three spaces with the following norm:
\begin{equation*}
\begin{aligned}
\|u\|_{H^2(\Omega)}^2 &:=\sum_{|\alpha|\le 2}\|D^\alpha u\|^2_{L^2(\X)}.
\end{aligned}
\end{equation*}
We will seek to prove existence and uniqueness of solutions to~(\ref{1}) in $H^2_{\beta,0}(\Omega)$.
\subsection{Tangential differential operators}
As mentioned previously, we denote by $n_{\partial\X}=(n_1,n_2)^T:\partial\Omega\to\mathbb S^1$, the unit outward normal to $\partial\Omega$. We define the unit tangent vector $\Ta:\partial\Omega\to\mathbb S^1$ as follows,
$$\Ta := (-[n_{\partial\X}]_2,[n_{\partial\X}]_1)^T.$$
For a given vector-valued function $v$ we denote:
$$v^\perp:=(-v_2,v_1)^T.$$
\begin{definition}[Tangential gradient]
We define the tangential gradient, $\nabla_\Ta:H^s(\partial\Omega)\to H^{s-1}_\Ta(\partial\Omega)$, as follows:
$$\nabla_\Ta v = \nabla v-\frac{\partial v}{\partial n_{\partial\X}}n_{\partial\X},$$
where $\frac{\partial v}{\partial n_{\partial\X}}=\nabla v\cdot n_{\partial\X}$, and $H^s_\Ta(\partial\Omega):=\{v\in H^s(\partial\Omega)\times H^s(\partial\Omega):v\cdot n_{\partial\X}= 0\}$.
\end{definition}
\begin{definition}[Directional derivative]\label{lemma1}
For a sufficiently smooth function $v$, and a vector $s$, we define the directional derivative of $u$ with respect to $s$ as follows:
$$\partial_su :=\nabla u\cdot s.$$
\end{definition}
\begin{definition}[Tangential Laplacian]
We define the tangential Laplacian, $\Delta_\Ta:H^s(\partial\Omega)\to H^{s-2}(\partial\Omega)$, as follows:
\begin{equation*}
\begin{aligned}
\Delta_\Ta &:= \operatorname{div}_\Ta\nabla_\Ta\\
& = \nabla_\Ta\cdot\nabla_\Ta.
\end{aligned}
\end{equation*}
\end{definition}
\begin{definition}[Mean curvature]
We define the mean curvature, $\mathcal{H}_{\partial\X}:\partial\Omega\to\mathbb R$, of $\partial\Omega$ as follows:
$$\mathcal{H}_{\partial\X}:=\nabla_\Ta\cdot n_{\partial\X}.$$
\end{definition}
\end{section}
\begin{section}{PDE analysis}\label{section:3}
We begin this section by defining the ``Cordes condition", which is crucial to the technique we will use to analyse problems in the form of~(\ref{1}); remarkably, the condition is a consequence of uniform ellipticity in two dimensions.
\begin{remark}[Minimal domain and oblique vector regularity]
We will always assume that the domain, $\X$, has a $C^2$ boundary, and that the oblique vector, $\beta\in C^1(\partial\X;\mathbb S^1)$.
\end{remark}
\begin{definition}[Cordes condition]
Let $L:=\sum_{i,j=1}^2A_{ij}D^2_{ij}$. The operator $L$ satisfies the Cordes condition if there exists $\varepsilon\in(0,1)$ such that
\begin{equation}\label{2}
\frac{|A|^2}{(\operatorname{Tr}(A))^2}\le\frac{1}{1+\varepsilon}\,\,\mbox{ a.e. in }\X.
\end{equation}
\end{definition}
\begin{lemma}
Assume that $A$ is uniformly elliptic. Then, $A$ satisfies the Cordes condition~(\ref{2}).
\end{lemma}
\emph{Proof:} Denote the lower and upper ellipticity constants of $A$ by $\lambda$ and $\Lambda$, respectively. Then,
$$\frac{|A|^2}{(\operatorname{Tr}(A))^2} = \frac{\lambda^2+\Lambda^2}{(\lambda+\Lambda)^2}=\frac{1}{1+2\lambda\Lambda/(\lambda+\Lambda)^2}.$$
Choosing $\varepsilon\in(0,2\lambda\Lambda/(\lambda+\Lambda)^2)$ yields the claim.$\quad\quad\square$
\begin{lemma}\label{lemma4}
Assume that $L$ satisfies the Cordes condition~(\ref{2}). Then, defining
$$\gamma:=\frac{\Tr A}{|A|^2},$$
we have that
$$|\gamma Lv-\Delta v|\le\sqrt{1-\varepsilon}|D^2v|,$$
for any $v\in H^2(\Omega)$.
\end{lemma}
\emph{Proof:} See~\cite{MR3077903}, Lemma 1.
\begin{definition}[Oblique angle]
We define the ``oblique angle", $\Theta:\partial\Omega\to\mathbb R$, to be the (anticlockwise) oriented angle between the oblique vector, $\beta$, and the unit outward normal, $n$.
\end{definition}
The following lemma \begingroup\color{black}will allow \endgroup us to prove the Miranda--Talenti estimate, which can be found in~\cite{MR2260015}. It also generalises the techniques present in~\cite{MR2260015} to the following bilinear form $\mathscr{B}:H^2(\X)\times H^2(\X)\to\R{}$
$$\mathscr{B}(u,v):=\int_\X D^2_{11}uD^2_{22}v+D^2_{22}uD^2_{11}v-2D^2_{12}uD^2_{12}v,$$
which will be used to prove an important consistency result in Section~\ref{section:4}.
\begin{lemma}\label{obl:lemma:applicabletotriangles}
Assume that $E\subset\R2$ is a bounded, Lipschitz, piecewise $C^{2}$ domain, and that $\beta\in C^1(\Gamma_n;\mathbb S^1)$ for each $C^{2}$ portion $\Gamma_n$ of $\partial E$, $n=1,\ldots,N$, $N\in\mathbb N$. Then, for any $u,v\in H^s(E)$, $s>5/2$, we have that
\begin{equation}\label{obl:identity:ref:forcorollary}
\begin{aligned}
&\int_E D^2_{11}uD^2_{22}v+D^2_{22}uD^2_{11}v-2D^2_{12}uD^2_{12}v\\
&~~~~~~~=\int_{\partial E}\left(\beta_1{\partial_{{\Ta}}}\beta_2-\beta_2{\partial_{{\Ta}}}\beta_1\right)(\betaperp\cdot\nabla u\,\betaperp\cdot\nabla v+\beta\cdot\nabla u\,\beta\cdot\nabla v)\\
&~~~~~~~~~~~~~+\int_{\partial E}\left({\partial_{{\Ta}}}(\betaperp\cdot\nabla u)\beta\cdot\nabla v-{\partial_{{\Ta}}}(\beta\cdot\nabla u)\betaperp\cdot\nabla v\right).
\end{aligned}
\end{equation}
\end{lemma}
\emph{Proof:} Let us momentarily assume that $u,v\in C^3(\overline{E})$, and note that an application of integration by parts gives us
\begin{equation}\label{obl:T}
\begin{aligned}
&\int_E D^2_{11}uD^2_{22}v+D^2_{22}uD^2_{11}v-2D^2_{12}uD^2_{12}v\\
&~~~~~~=\int_E D_1(D_1vD^2_{22}u-D_2vD^2_{12}u)-D_2(D_1vD^2_{21}u-D_2vD^2_{11}u)\\
&~~~~~~= \int_{\partial E}(D_1vD^2_{22}u-D_2vD^2_{12}u)n_1-(D_1vD^2_{21}u-D_2vD^2_{11}u)n_2.
\end{aligned}
\end{equation}
Now, denoting $C^1_u:=\betaperp\cdot\nabla u$, $C^1_v:=\betaperp\cdot\nabla v$, $C^2_u:=\beta\cdot\nabla u$, and $C^2_v:=\beta\cdot\nabla u$, we obtain the following linear systems on $\Gamma_n$, $n=1,\ldots,N$:
\begin{equation}\label{obl:star:u}
\left\{
\begin{aligned}
\beta_1D_1u+\beta_2D_2u & = C^2_u,\\
-\beta_2D_1u+\beta_1D_2u & = C^1_u,
\end{aligned}
\right.
\end{equation}
\begin{equation}\label{obl:star:v}
\left\{
\begin{aligned}
\beta_1D_1v+\beta_2D_2v & = C^2_v,\\
-\beta_2D_1v+\beta_1D_2v & = C^1_v,
\end{aligned}
\right.
\end{equation}
with corresponding unique solutions
\begin{equation}\label{obl:star:star}
\left\{
\begin{aligned}
D_1u = -C^1_u\beta_2+C^2_u\beta_1,\quad D_2u = C^1_u\beta_1+C^2_u\beta_2,\\
D_1v = -C^1_v\beta_2+C^2_v\beta_1,\quad D_2v = C^1_v\beta_1+C^2_v\beta_2.
\end{aligned}
\right.
\end{equation}
Substituting in the values of $D_1v$ and $D_2v$ present in~(\ref{obl:star:star}) into~(\ref{obl:T}), (denoting $n_i = [n_{\partial\X}]^i$, $i=1,2$) we obtain the following:
\begin{equation}\label{obl:TT}
\begin{aligned}
&\int_E D^2_{11}uD^2_{22}v+D^2_{22}uD^2_{11}v-2D^2_{12}uD^2_{12}v\\
&~~~~~~=\int_{\partial E}-C^1_v\beta_2D^2_{22}u\,n_1+C^2_v\beta_1D^2_{22}u\,n_1-C^1_v\beta_1D^2_{12}u\,n_1-C^2_v\beta_2D^2_{12}u\,n_1\\
&~~~~~~~~~~~~+\int_{\partial E}C^1_v\beta_2D^2_{21}u\,n_2-C^2_v\beta_1D^2_{21}u\,n_2+C^1_v\beta_1D^2_{11}n_2+C^2_v\beta_2D^2_{11}u\,n_2\\
&~~~~~~ = \int_{\partial E}C^1_v(-\beta_2D^2_{22}u\,n_1-\beta_1D^2_{12}u\,n_1+\beta_2D^2_{21}u\,n_2+\beta_1D^2_{11}u\,n_1)\\
&~~~~~~~~~~~~+\int_{\partial E}C^2_v(\beta_1D^2_{22}u\,n_1-\beta_2D^2_{12}u\,n_1-\beta_1D^2_{21}n_2+\beta_2D^2_{11}u\,n_2).
\end{aligned}
\end{equation}
Taking the directional derivative of the equations in~(\ref{obl:star:u}) with respect to $\Ta$ yields (denoting $\dot{~}:={\partial_{{\Ta}}}$)
\begin{equation*}
\left\{
\begin{aligned}
\dot{\beta_1}D_1u+\beta_1D^2_{11}u(-n_2)+\beta_1D^2_{12}u(n_1)+\dot{\beta_2}D_2u+\beta_2D^2_{12}u(-n_2)+\beta_2D^2_{22}u(n_1) &\!=\!\dot{C}^2_u,\\
-\dot{\beta_2}D_1u-\beta_2D^2_{11}u(-n_2)-\beta_2D^2_{12}u(n_1)+\dot{\beta_1}D_2u+\beta_1D^2_{12}u(-n_2)+\beta_1D^2_{22}u(n_1) &\!=\!\dot{C}^1_u,
\end{aligned}
\right.
\end{equation*}
thus,
\begin{equation}\label{obl:star:star:star}
\left\{
\begin{aligned}
-(\beta_1D^1_{11}u\,n_2-\beta_1D^2_{12}u\,n_1+\beta_2D^2_{12}u\,n_2-\beta_2D^2_{22}u\,n_1)& = \dot{C}^2_u-(\dot{\beta_1}D_1u+\dot{\beta_2}D_2u),\\
(\beta_1D^2_{22}u\,n_1-\beta_2D^2_{12}u\,n_1-\beta_1D^2_{12}u\,n_2+\beta_2D^2_{11}u\,n_2)& = \dot{C}^1_u+\dot{\beta_2}D_1u-\dot{\beta_1}D_2u.
\end{aligned}
\right.
\end{equation}
Substituting the values for $D_1u$ and $D_2u$ present in~(\ref{obl:star:star}) gives us
\begin{equation}\label{obl:sub:1}
\begin{aligned}
\dot{\beta_1}D_1u+\dot{\beta_2}D_2u & = -C^1_u\beta_2\dot{\beta_1}+C^2_u\beta_1\cdot{\beta_1}+C^1_u\beta_1\dot{\beta_2}+C^2_u\beta_2\dot{\beta_2}\\
& = C^1_u(\beta_1\dot{\beta_2}-\beta_2\dot{\beta_1})+C^2_u(\beta_1\dot{\beta_1}+\beta_2\dot{\beta_2})\\
& = C^1_u(\beta_1\dot{\beta_2}-\beta_2\dot{\beta_1});
\end{aligned}
\end{equation}
note that the latter equality holds, since
 $$\beta_1\dot{\beta_1}+\beta_2\dot{\beta_2}=\frac{1}{2}{\partial_{{\Ta}}}(|\beta|^2) = \frac{1}{2}{\partial_{{\Ta}}}(1) = 0.$$
 Similarly, we obtain
 \begin{equation}\label{obl:sub:2}
 \dot{\beta_2}D_1u-\dot{\beta_1}D_2u = C_u^2(\beta_1\dot{\beta_2}-\beta_2\dot{\beta_1}).
 \end{equation}
Substituting~(\ref{obl:sub:1}) and~(\ref{obl:sub:2}) into~(\ref{obl:star:star:star}), and then substituting the resulting equations into~(\ref{obl:TT}), we obtain
\begin{equation*}
\begin{aligned}
&\int_E D^2_{11}uD^2_{22}v+D^2_{22}uD^2_{11}v-2D^2_{12}uD^2_{12}v\\
&~~~~~~~~~~~+\int_{\partial E}C^1_v(C^1_u(\beta_1\dot{\beta_2}-\beta_2\dot{\beta_1})-\dot{C}^1_u)+C^2_v(C^2_u(\beta_1\dot{\beta_2}-\beta_2\dot{\beta_1})+\dot{C}^2_u)\\
&~~~~~~=\int_{\partial E}(\beta_1\dot{\beta_2}-\beta_2\dot{\beta_1})(C^1_uC^1_v+C^2_uC^2_v)+\dot{C}^1_uC^2_v-\dot{C}^2_uC^1_v\\
&~~~~~~=\int_{\partial E}\left(\beta_1{\partial_{{\Ta}}}\beta_2-\beta_2{\partial_{{\Ta}}}\beta_1\right)(\betaperp\cdot\nabla u\,\betaperp\cdot\nabla v+\beta\cdot\nabla u\,\beta\cdot\nabla v)\\
&~~~~~~~~~~~~~+\int_{\partial E}\left({\partial_{{\Ta}}}(\betaperp\cdot\nabla u)\beta\cdot\nabla v-{\partial_{{\Ta}}}(\beta\cdot\nabla u)\betaperp\cdot\nabla v\right),
\end{aligned}
\end{equation*}
which is exactly~(\ref{obl:identity:ref:forcorollary}). This identity extends to $u,v\in H^s(E)$, $s>5/2$, by density, which concludes the proof.$\quad\quad\square$
\begin{corollary}\label{cor:1thu}
Assume that $\X\subset\R2$ is a bounded $C^2$ domain, and that $\beta\in C^1(\partial\X;\mathbb S^1)$. Then, for all $u,v\in H^s(\X)$, $s>5/2$, we have that
\begin{equation}\label{true:identity}
\begin{aligned}
&\int_\X D^2_{11}uD^2_{22}v+D^2_{22}uD^2_{11}v-2D^2_{12}uD^2_{12}v = \\
&~~~~~~~\int_{\partial\X}\left(\partial_\Ta\Theta+\mc_{\partial\X}\right)(\betaperp\cdot\nabla u\,\betaperp\cdot\nabla v+\beta\cdot\nabla u\,\beta\cdot\nabla v)\\
&~~~~~~~~~~~~~+\int_{\partial\X}\left(\partial_\Ta(\betaperp\cdot\nabla u)\beta\cdot\nabla v-\partial_\Ta(\beta\cdot\nabla u)\betaperp\cdot\nabla v\right).
\end{aligned}
\end{equation}
\end{corollary}
\emph{Proof:} This follows from applying Lemma~\ref{obl:lemma:applicabletotriangles}, and applying the following equality
\begingroup\color{black}
\begin{equation}\label{star:13:16}
\beta_1\partial_\Ta\beta_2-\beta_2\partial_\Ta\beta_1=\partial_\Ta\Theta+\mc_{\partial\X}\quad\mbox{on }\partial\X.
\end{equation}
\endgroup
Note that the above equality follows from identities (1.83) and (1.84) on page 48 in~\cite{MR2260015} (note that to pass from the present notation, to the notation found in~\cite{MR2260015}, one must denote $\ell:=\beta$, $\theta:=\Theta$, and $\chi:=-\mc_{\partial\X}$.)$\quad\quad\square$
\begin{corollary}\label{cor:applicabletotriangles}
Assume that $E\subset\R2$ is a bounded, Lipschitz, piecewise $C^{2}$ domain, and that $\beta\in C^1(\Gamma_n;\mathbb S^1)$ for each $C^{2}$ portion $\Gamma_n$ of $\partial E$, $n=1,\ldots,N$, $N\in\mathbb N$. Then, for any $u,v\in H^s(E)$, $s>5/2$, we have that
\begin{equation*}
\begin{aligned}
\int_ED^2u:D^2v&+\int_{\partial E}\left(\beta_1\partial_\Ta\beta_2-\beta_2\partial_\Ta\beta_1\right)(\betaperp\cdot\nabla u\,\betaperp\cdot\nabla v+\beta\cdot\nabla u\,\beta\cdot\nabla v)\\
&~~~~~~+\int_{\partial E}\left(\partial_\Ta(\betaperp\cdot\nabla u)\beta\cdot\nabla v-\partial_\Ta(\beta\cdot\nabla u)\betaperp\cdot\nabla v\right)\\
&~~~~~~~~~~~~=\int_E\Delta u\,\Delta v.
\end{aligned}
\end{equation*}
\end{corollary}
\emph{Proof:} First note that for $u,v\in H^2(E)$,
$$D^2u:D^2v+D^2_{11}uD^2_{22}v+D^2_{22}uD^2_{11}v-2D^2_{12}uD^2_{12}v = \Delta u\,\Delta v,$$
and apply Lemma~\ref{obl:lemma:applicabletotriangles}.$\quad\quad\square$
\begin{lemma}[Miranda--Talenti estimate]\label{MT}
Let $\X\subset \mathbb R^2$ be a $C^2$ domain, and assume that $\beta\in C^1(\partial\Omega;\mathbb S^1)$. Furthermore, assume that
$$\partial_\Ta\Theta+\mc_{\partial\X}\ge 0\,\,\mbox{on}\,\,\partial\Omega.$$
Then, we have that
\begin{equation}\label{MT1}|u|^2_{H^2(\Omega)}:=\int_\Omega\sum_{i,j=1}^2(D^2_{ij}u)^2\le\int_\Omega(\Delta u)^2=\|\Delta u\|_{L^2(\Omega)}^2,
\end{equation}
for all $u\in H^2_\beta(\Omega)$.
\end{lemma}
\emph{Proof:} First, we note that since $\X\subset\mathbb R^2$ is a $C^2$ domain, it is both Lipschitz continuous and piecewise $C^2$. Furthermore, as $\beta\in C^1(\partial\X;\mathbb S^1)$, it follows that $\beta\in C^1(\Gamma_{\partial\X};\mathbb S^1)$ for each $C^2$ portion $\Gamma_{\partial\X}$ of $\partial\X$; indeed one may take $\Gamma_{\partial\X} = \partial\X$. Now, let us assume that $u\in C^3(\overline{\X})\cap H^2_{\beta,0}(\X)$, so that $u\in H^s(\X)$, with $s>5/2$. Setting $v=u$, it then follows that the hypotheses of Corollaries~\ref{cor:1thu} and~\ref{cor:applicabletotriangles} are satisfied. It then follows that
\begin{equation*}
\begin{aligned}
\int_\X D^2u:D^2u&+\int_{\partial\X}\left(\partial_\Ta\Theta+\mc_{\partial\X}\right)(\betaperp\cdot\nabla u\,\betaperp\cdot\nabla u+\beta\cdot\nabla u\,\beta\cdot\nabla u)\\
&~~~~~~+\int_{\partial\X}\left(\partial_\Ta(\betaperp\cdot\nabla u)\beta\cdot\nabla u-\partial_\Ta(\beta\cdot\nabla u)\betaperp\cdot\nabla u\right)\\
&~~~~~~~~~~~~=\int_\X\Delta u\,\Delta u.
\end{aligned}
\end{equation*}
Since $\beta\cdot\nabla u|_{\partial\X} = C$ for some constant $C$, it follows that $$\partial_\Ta(\beta\cdot\nabla u)|_{\partial\X} = 0.$$ Furthermore, we see that
$$\betaperp\cdot\nabla u\,\betaperp\cdot\nabla u+\beta\cdot\nabla u\,\beta\cdot\nabla u=|\nabla u|^2.$$
Thus, we obtain
$$\int_\X |D^2u|^2+\int_{\partial\X}\left(\partial_\Ta\Theta+\mc_{\partial\X}\right)|\nabla u|^2+C\partial_\Ta(\betaperp\cdot\nabla u) = \int_\X(\Delta u)^2.$$
Since $\partial\X$ is a compact hyper-surface, an application of integration by parts yields
$$\int_{\partial\X}C\partial_\Ta(\betaperp\cdot\nabla u) = 0.$$
Finally, since 
$$\partial_\Ta\Theta+\mc_{\partial\X}\ge 0\,\,\mbox{on}\,\,\partial\Omega,$$
we obtain
$$|u|^2_{H^2(\Omega)}:=\int_\Omega\sum_{i,j=1}^2(D^2_{ij}u)^2=\int_\X|D^2u|^2\le\int_\Omega(\Delta u)^2=\|\Delta u\|_{L^2(\Omega)}^2,$$
as desired.$\quad\quad\square$
\begin{lemma}[Gradient estimate]\label{MTgradient}
Under the assumptions of Lemma~\ref{MT}, with the additional assumption that
$$\partial_\Ta\Theta+\mc_{\partial\X} > 0\,\,\mbox{on}\,\,\partial\Omega,$$
we have that
$$\int_\Omega|\nabla u|^2\le C\int_\Omega(\Delta u)^2,$$
for all $u\in H^2_\beta(\Omega)$, where $C$ is a constant independent of $u$.
\end{lemma}
\emph{Proof:} See~\cite{MR2260015}, Lemma 1.5.8.$\quad\quad\square$
\begin{corollary}\label{FullMT}
Under the assumptions of Lemma~\ref{MTgradient}, we have that
\begin{equation}\label{MT2}
\|u\|_{H^2(\Omega)}^2\le C\|\Delta u\|_{L^2(\Omega)}^2,
\end{equation}
for all $u\in H^2_{\beta,0}(\Omega)$, where the constant $C$ is independent of $u$.
\end{corollary}
\emph{Proof:} Applying Lemmas~\ref{MT} and~\ref{MTgradient}, we obtain{}{}
\begin{equation}\label{star}
\begin{aligned}
\|u\|_{H^2(\Omega)}^2 & = \|u\|_{L^2(\Omega)}^2 +\|\nabla u\|_{L^2(\Omega)}^2 +\|D^2u\|_{L^2(\Omega)}^2 \\
&\le \|u\|_{L^2(\Omega)}^2+(C+1)\|\Delta u\|_{L^2(\Omega)}^2.
\end{aligned}
\end{equation}
Now, by Poincar\'e's inequality, we obtain
\begin{equation*}
\begin{aligned}
\|u\|_{L^2(\Omega)}^2 & = \left\|u-\frac{1}{|\Omega|}\int_\Omega u\right\|_{L^2(\Omega)}^2\\
&\le C\int_\Omega|\nabla u|^2\\
&\le C\|\Delta u\|^2_{L^2(\Omega)},
\end{aligned}
\end{equation*}
combining this with~(\ref{star}), we obtain the desired result.$\quad\quad\square$
\begin{theorem}\label{linear:existence+uniqueness}
Under the assumptions of Lemma~\ref{MTgradient}, there exists a unique $u\in H^2_{\beta,0}(\Omega)$ that is a strong solution of~(\ref{1}).
\end{theorem}
\emph{Proof:} Define $A_\gamma:H^2_{\beta,0}(\Omega)\times H^2_{\beta,0}(\Omega)\to\mathbb R$, and $l_\gamma:H^2_{\beta,0}(\Omega)\to\mathbb R$ as follows:
$$A_\gamma(u,v):=\int_\Omega\gamma A:D^2u\,\Delta v\,\,\quad\forall u,v\in H^2_{\beta,0}(\Omega),$$
and
$$l_\gamma(v):=\int_\Omega\gamma f\Delta v\,\,\quad\forall v\in H^2_{\beta,0}(\Omega).$$
We will first prove that there exists a unique $u\in H^2_{\beta,0}(\Omega)$ such that
$$A_\gamma(u,v) = l_\gamma(v)\,\,\quad\forall v\in H^2_{\beta,0}(\Omega).$$
We see that
\begin{equation*}
\begin{aligned}
A_\gamma(u,v)&\le\|\gamma\|_{L^\infty(\Omega)}\|A\|_{L^\infty(\Omega)}|u|_{H^2(\Omega)}\|\Delta v\|_{L^2(\Omega)}\\
&\le C\|u\|_{H^2(\Omega)}\|v\|_{H^2(\Omega)},
\end{aligned}
\end{equation*}
similarly
$$l_\gamma(v)\le C\|v\|_{H^2(\Omega)}.$$
Futhermore,
\begin{equation*}
\begin{aligned}
A_\gamma(u,u)& = \int_\Omega(\gamma Lu-\Delta u)\Delta u+\int_\Omega(\Delta u)^2\\
&\ge -\sqrt{1-\varepsilon}\int_\Omega|D^2u||\Delta u|+\int_\Omega(\Delta u)^2\\
&\ge -\sqrt{1-\varepsilon}\|D^2u\|_{L^2(\Omega)}\|\Delta u\|_{L^2(\Omega)}+\|\Delta u\|_{L^2(\Omega)}^2\\
&\ge(1-\sqrt{1-\varepsilon})\|\Delta u\|_{L^2(\Omega)}^2\\
&\ge C\|u\|_{H^2(\Omega)}^2.
\end{aligned}
\end{equation*}
Note that the second inequality follows from Lemma~\ref{lemma4}, and the final inequalities follow from Lemma~\ref{MT} and Corollary~\ref{FullMT}.

We have shown that $A_\gamma$ is both coercive and bounded, and $l_\gamma$ is bounded. In view of the fact that $H^2_{\beta,0}(\Omega)$ is a Hilbert space when endowed with the inner product $$\langle u,v\rangle_\Delta:=\int_\Omega\Delta u\,\Delta v\,\,\quad\forall u,v\in H^2_{\beta,0}(\Omega),$$
the Lax Milgram Theorem (\cite{MR2597943} Section 6.2.1) yields existence and uniqueness of $u\in H^2_{\beta,0}(\Omega)$ such that
$$A_\gamma(u,v) = l_\gamma(v)\,\,\quad\forall v\in H^2_{\beta,0}(\Omega),$$
i.e.,
$$\int_\Omega\gamma A:D^2u\,\Delta v  = \int_\Omega\gamma f\,\Delta v\,\,\quad\forall v\in H^2_{\beta,0}(\Omega).$$
Now, it follows from~\cite{MR2260015} Page 56, that $\Delta:H^2_{\beta,0}(\Omega)\to L^2(\Omega)$ is a surjection, and thus we obtain
$$\int_\Omega\gamma (A:D^2u-f) v  = 0\,\,\quad\forall v\in L^2(\Omega).$$
An application of the Fundamental Lemma of the calculus of variations yields
$$\gamma(A:D^2u-f) = 0\,\,\quad\mbox{a.e. in}\,\,\Omega,$$
and since $\gamma$ is uniformly positive, we obtain
$$A:D^2u = f\,\,\quad\mbox{a.e. in}\,\,\Omega.\quad\quad\square$$
\end{section}
\begin{section}{Numerical method}\label{section:4}
\subsection{Notation}
\begin{definition}[Edge and vertex sets]
Given a triangulation $\Th$, we denote by \begingroup\color{black}$\Eb$, the set of boundary edges of $\Th$, by $\mathscr{E}_h^i$ the set of interior edges of $\Th$, by $\mathscr{E}_h^{i,b}:=\mathscr{E}_h^{i}\cup\mathscr{E}_h^{b}$\endgroup, and by $\mathscr{V}^b_h$ the set of boundary vertices of $\Th$.
\end{definition}
\begingroup\color{black}
\begin{definition}[Piecewise $C^{k}$ domain]\label{pwck21:05}
A domain $\Omega\subset\mathbb R^d$ is piecewise $C^{k}$ for $k\in\mathbb N$, if we may express the boundary of $\Omega$, $\partial\Omega$, as a finite union
\begin{equation}\label{piecewiserep}
\partial\Omega = \bigcup_{n=1}^N\overline{\Gamma_n},
\end{equation}
where each $\Gamma_n\subset\mathbb R^d$ is of zero $d$-dimensional Lebesgue measure, and admits a local representation as the {}graph of a uniformly $C^{k}$ function. That is, for each $n$, and at each $x\in\Gamma_n$ there exists an open neighbourhood $V_n$ of $x$ in $\mathbb R^d$ and an orthogonal coordinate system $(y^n_1,\ldots,y^n_d)$, such that $$V_n=\{(y^n_1,\ldots,y^n_d):-a_j^n<y_j^n<a_j^n,1\le j\le d\};$$
as well as a uniformly $C^k$ function $\varphi_n$ defined on $V'_n=\{(y_1^n,\ldots,y^n_{d-1}):-a_j^n<y_j^n<a_j^n,1\le j\le d-1\}$ and such that
\begin{align*}
&|{\varphi_n}({y^{n}}')|\le a^n_d/2\mbox{ for every }{y^n}'=(y^n_1,\ldots,y^n_{d-1})\in V_{n'},\\ 
&\X\cap V = \{y^n=({y^n}',y^n_d)\in V:y^n_d<\varphi_n({y^n}')\},\\
&\Gamma_n\cap V = \{y^n=({y^n}',y^n_d)\in V:y^n_d=\varphi_n({y^n}')\}.
\end{align*}
\end{definition}
\endgroup
\emph{Jump and average operators.}
For each face $F=\overline{K}\cap\overline{K'}$ for some $K,K'\in\Th$ (in the case that $F\in\mathscr{E}^b_h$ take $F=\partial K\cap\partial\X$), with corresponding unit normal vector $n_F$ (which, for convention is chosen so that it is the outward normal to $K$), we define the jump operator, $\jump{\cdot}$ over $F$, by
\begin{equation*}
\jump{\phi} = \left\{
\begin{aligned}
&(\phi|_K)|_F-(\phi|_{K'})|_F\,\,\mbox{if}\,\,F\in\Ei\\
&(\phi|_K)|_F\,\,\mbox{if}\,\,F\in\Eb,\\
\end{aligned}
\right.
\end{equation*} 
and the average operator, $\avg{\cdot}$, by
\begin{equation*}
\avg{\phi} = \left\{
\begin{aligned}
&\frac{1}{2}((\phi|_K)|_F+(\phi|_{K'})|_F)\,\,\mbox{if}\,\,F\in\Ei\\
&(\phi|_K)|_F\,\,\mbox{if}\,\,F\in\Eb.\\
\end{aligned}
\right.
\end{equation*} 
Each vertex $e=\overline{F}\cap\overline{F'}$ for some $F,F'\in\Eb$. We thus define the jump and average over a vertex $e\in\Vb$ analogously.

For two matrices $A,B\in\mathbb R^{m\times n}$, we set $A:B:=\sum_{i,j=1}^{m,n}A_{ij}B_{ij}$. For an element $K$, we define the bilinear form $\langle\cdot,\cdot\rangle_K$ by
\begin{equation*}
\langle u,v\rangle_K :=
\int_Ku:v\,\,\mbox{if}\,\,u,v\in L^2(K;\mathbb R^{m\times n}).
\end{equation*} 
Any ambiguity in this notation will be resolved by arguments of the bilinear form. The bilinear forms $\langle\cdot,\cdot\rangle_{\partial K}$ and $\langle\cdot,\cdot\rangle_{F}$ for $F\in\Eib$, are defined similarly.
\subsection{Exact domain approximation} 
We will continue this section by providing the details of~\cite{MR1014883}, which provides us with a notion of exact domain approximation
\begin{definition}[Curved $d$-simplex]\label{curveddsimp}
An open set $K\subset\mathbb R^d$ is called a curved $d$-simplex if there exists a $C^1$ mapping $F_K$ that maps a straight reference $d$-simplex $\hat K$ onto $K$, and that is of the form
\begin{equation}\label{non:aff:map}
F_K=\tilde{F}_K+\Phi_K,
\end{equation}
where 
\begin{equation}\label{def:affine-map}
\tilde{F}_K:\hat{x}\mapsto\tilde{B}_K\hat x+\tilde{b}_K
\end{equation}
is an invertible map and $\Phi_K\in C^1(\hat{K};\mathbb R^d)$ satisfies
\begin{equation}\label{CK:def}
C_K:=\sup_{\hat{x}\in\hat K}\|D\Phi_K(\hat x)\tilde{B}_K^{-1}\|<1,
\end{equation}
where $\|\cdot\|$ denotes the induced Euclidean norm on $\mathbb R^{d\times d}$.
\end{definition}

\begin{definition}[Class $m$ curved $d$-simplex]
A curved $d$-simplex $K$ is of class $C^m$, $m\ge1$, if the mapping $F_K$ is of class $C^m$ on $\hat K$.
\end{definition}
\begin{definition}[Mesh size] For each element $K\in\Th$, let $h_K:=\operatorname{diam}\tilde{K}\ge C_{\mathcal{F}}\|\tilde{B}_K\|$ (where $\tilde{K}=\tilde{B}_K(\hat K)$). It is assumed that $h = \max_{K\in\Th} h_K$ for each mesh $\Th$.
\begingroup\color{black}Furthermore, for each face $F\in\mathscr{E}^{i,b}_h$, we define
\begin{equation}\label{hF:def}
\tilde{h}_F:=\left\{\begin{array}{l l}
\min(h_K,h_{K'}) & \mbox{if}~F\in\mathscr{E}^i_h,\\
h_K & \mbox{if}~F\in\mathscr{E}^b_h.
\end{array}\right.
\end{equation}
where $K$ and $K'$ are such that $F=\partial K\cap\partial K'$ if $F\in\mathscr{E}^i_h$, or $F\subset\partial K\cap\partial\X$ if $F\in\mathscr{E}^b_h$. Finally, for each $e\in\Vb$, we define
\begin{equation}\label{hedef}
h_e:=\min_{F\in\Ei:\overline{F}\cap e\ne\emptyset}\tilde{h}_F.
\end{equation}
\endgroup
\end{definition}
\begingroup\color{black}
\begin{definition}
The family $(\Th)_h$ of meshes is said to be regular if there exist two constants, $\sigma$ and $c$, independent of $h$, such that, for each $h$, any $K\in\Th$ satisfies
\begin{equation}\label{shape-reg}
h_K/\rho_K\le\sigma,
\end{equation}
where $\rho_K$ is the diameter of the sphere inscribed in $\tilde{K}$.
Furthermore, we have
\begin{equation}\label{regularpart}
\sup_h\sup_{K\in\Th} C_K\le c<1.
\end{equation}
\end{definition}
\begin{definition}\label{regoforderm}
The family $(\Th)_h$ of meshes is said to be regular of order $m$ if it is regular and if, for each $h$, any $K\in\Th$ is of class $C^{m+1}$, with
\begin{equation}\label{regoforderms}
\sup_h\sup_{K\in\Th}\sup_{\hat x\in\hat K}\|D^lF_K(\hat x)\|\|\tilde{B}_K\|^{-l}<\infty,\quad 2\le l\le m+1.
\end{equation}
\end{definition}
\begin{assumption}\label{Meshconds}
The meshes are allowed to be irregular, i.e., there may be hanging nodes. We assume that there is a uniform upper bound on the number of faces composing the boundary of any given element; in other words, there is a $C_\mathcal{F}>0$, independent of $h$, such that
\begin{equation}\label{meshcond1}
\max_{K\in\mathscr{T}_h}\operatorname{card}\{F\in\mathscr{E}^{i,b}_h:F\subset\partial K\}\le C_\mathcal{F}\quad\forall K\in\mathscr{T}_h,~\forall h>0.
\end{equation}
We assume that any two elements sharing a face have commensurate diameters, i.e., there is a $C_{\mathcal{T}}\ge 1$, independent of $h$, such that
\begin{equation}\label{meshcond2}
\max(h_K,h_{K'})\le C_\mathcal{T}\min(h_K,h_{K'}),
\end{equation}
for any $K$ and $K'$ in $\mathscr{T}_h$ that share a face. 

Finally, we assume that each $F\in\Eb$ satisfies
\begin{equation}\label{Fcontained}
F = F\cap\Gamma_n,
\end{equation}
for some $n\in\{1,\ldots,N\}$, with $\Gamma_n$ given as in~(\ref{piecewiserep}). This implies that each boundary face is completely contained in a boundary portion $\Gamma_n$, as well as ensuring that our approximation of the domain $\X$ is exact.
\end{assumption}
\begin{remark}
The assumptions on the mesh given by Assumption~\ref{Meshconds}, 
in particular~(\ref{meshcond2}), show that if $F$ is a face of $K$, then
\begin{equation}\label{meshcondcons}
h_K\le C_\mathcal{T}\tilde{h}_F.
\end{equation}
\end{remark}
\endgroup
\subsection{Finite element spaces}
 For each $K\in\Th$, we denote by $\mathbb P^{p}(K)$ the space of all polynomials on $K$ with total degree less than or equal to $p$. The discontinuous Galerkin finite element space $V_{h,\mathbf p}$ is defined by
\begin{equation}
\dg \begingroup\color{black}:=\endgroup\{v\in L^2(\Omega):~v|_K=\rho\circ F_K^{-1},\,\rho\in\mathbb P^{p}(\hat K),~\forall K\in\Th\},
\end{equation}
we also define the subspace, $\dgo$, of $\dg$ as follows
$$\dgo:=\{v\in\dg:\int_\Omega v = 0\}.$$
Let $\mathbf s=(s_K:K\in\Th)$ denote a vector of nonnegative real numbers and let $r\in[1,\infty]$.\\
The broken Sobolev space $W^{\mathbf s,r}(\Omega;\Th)$ is defined by
\begin{equation}
W^{\mathbf s,r}(\Omega;\Th):=\{v\in L^2(\Omega):~v|_K\in W^{s_K,r}(K)~\forall K\in\Th\}.
\end{equation}
We denote $H^\mathbf{s}(\Omega;\Th):=W^{\mathbf s,2}(\Omega;\Th)$, and set $W^{s,r}(\Omega;\Th):=W^{\mathbf s,r}(\Omega;\Th)$, in the case that $s_K=s,~s\ge 0$, for all $K\in\Th$. For $v\in W^{1,r}(\Omega;\Th)$, let 
$\nabla_hv\in L^r(\Omega;\mathbb R^d)$ denote the discrete (also known as broken) gradient of $v$, i.e., $(\nabla_hv)|_K=\nabla(v|_K)$ for all $K\in\Th$. Higher order discrete derivatives are defined in a similar way. We define a norm on $W^{s,r}(\Omega;\Th)$ by
\begin{equation}
\|v\|^r_{W^{s,r}(\Omega;\Th)}:=\sum_{K\in\Th}\|v\|^r_{W^{s,r}(K)}
\end{equation}
with the usual modification when $r=\infty$.

The following Lemma is a direct application of Lemma 4 from~\cite{MR3077903}.
\begin{lemma}
Let $\Omega$ be a bounded domain, and let $\mathscr{T}_h$ be a mesh on $\Omega$ consisting of \emph{possibly curved} simplices. Then, for each $K\in\mathscr{T}_h$ and each face $F\subset\partial K$ that belongs to $\Ei$, the following identities hold:
\begin{equation}\label{identities}
\begin{aligned}
\tau_F(\nabla v)  &= \gradTa(\tau_Fv)+\left(\tau_F\frac{\partial v}{\partial n_F}\right)n_F\quad\forall v\in H^s(K),\,s>3/2,\\
\tau_F(\Delta v)  &= \operatorname{div}_T\gradTa(\tau_Fv)+\tau_F\frac{\partial}{\partial n_F}(\nabla v\cdot n_F),\quad\forall v\in H^s(K),\,s>5/2.
\end{aligned}
\end{equation}
\end{lemma}
\subsection{The design of the numerical method}\label{obldes}
\begingroup\color{black}
We shall now discuss how the terms in the definition bilinear form that defines the finite element method of this paper, arise. We are motivated by the desire to numerically enforce the Miranda--Talenti (MT) estimates~(\ref{MT1}) and~(\ref{MT2}), whilst producing a  scheme that is both consistent and symmetric (the latter occurs when the operator $A\!:\!D^2$ is isotropic).

We solve for both $u_h\in\dgo:=\dg\cap L^2_0(\X)$, which approximates the strong solution $u\in H^2_{\beta,0}(\X)$ of~(\ref{1}), and $c_h\in\dgoo$, which approximates the compatibility constant of~(\ref{1}), that is, $c_h$ approximates the value of $C=\beta\cdot\nabla u|_{\partial\X}$ (the value of $C$ is a priori \emph{unknown}). As such, our finite element space will be
$$M_h:=\dgo\times\dgoo.$$
 We first note that the bilinear form, $A_h:M_h\times M_h\to\mathbb R$, that we use to define the finite element method, will take the following structure:
\begin{equation}\label{Aoblpredef12:21}
\begin{aligned}
\Aobl((u_h,c_h);(v_h,\mu_h)) & := \sum_{K\in\Th}\langle\gamma A\!:\!D^2u_h,\Delta v_h\rangle_K+\Bhalfobl((u_h,c_h);(v_h,\mu_h))\\
&~~~~~~~~~~-\sum_{K\in\Th}\langle\Delta u_h,\Delta v_h\rangle_K\quad\forall(u_h,c_h),(v_h,\mu_h)\in M_h.
\end{aligned}
\end{equation}
We \emph{claim} that the bilinear form $\Bhalfobl$ is coercive on $M_h\times M_h$, and that
\begin{equation}\label{cons11:20}
\Bhalfobl((w,c);(v_h,\mu)) = \sum_{K\in\Th}\langle\Delta w,\Delta v_h\rangle_K,\quad\forall(v_h,\mu)\in\dgo\times\mathbb R,
\end{equation}
when $w\in H^2_{\beta,0}(\X)\cap H^s(\X;\Th)$, $s>5/2$, and $c = \beta\cdot\nabla w|_{\partial\X}$.

It is then clear that~(\ref{cons11:20}) implies that
\begin{equation}\label{cons11:20:2}
\Aobl((w,c);(v_h,\mu)) =  \sum_{K\in\Th}\langle\gamma A\!:\!D^2u_h,\Delta v_h\rangle_K =: A_{\gamma,h}(w,v_h),\quad\forall(v_h,\mu)\in\dgo\times\mathbb R,
\end{equation}
for the aforementioned choice of $w$ and $c$. Note that
$A_{\gamma,h}$ is a numerical discretisation of $A_\gamma$, defined by~(\ref{agammadef}).

The bilinear form $\Bhalfobl:M_h\times M_h\to\mathbb R$ takes the following form
\begin{equation}\label{bhalfobldef}
\begin{aligned}
\Bhalfobl((u_h,c_h);(v_h,\mu_h))&:=\frac{1}{2}\Bstarobl((u_h,c_h);(v_h,\mu_h))+\frac{1}{2}\sum_{K\in\Th}\langle\Delta u_h,\Delta v_h\rangle_K\\
&~~~~~~+\Jhobl((u_h,c_h);(v_h,\mu_h)),
\end{aligned}
\end{equation}
where the bilinear forms $\Bstarobl,\Jhobl:M_h\times M_h\to\mathbb R$ satisfy
\begin{equation}\label{cons13:46}
\Bstarobl((w,c),(v_h,\mu)) = \sum_{K\in\Th}\langle\Delta w,\Delta v_h\rangle_K,\quad\mbox{and}\quad\Jhobl((w,c),(v_h,\mu)) = 0,
\end{equation}
for all $(v_h,\mu)\in\dgo\times\mathbb R$, when $w\in H^2_{\beta,0}(\X)\cap H^s(\X;\Th)$, $s>5/2$, and $c = \beta\cdot\nabla w|_{\partial\X}$.

Moreover, one can see that~(\ref{cons13:46}) implies~(\ref{cons11:20}), which in turn implies~(\ref{cons11:20:2}). We also remark that the bilinear form $\Jhobl$ is a jump penalty term that enforces regularity that is consistent with that of the true solution. In particular, if $w\in H^2_{\beta,0}(\X)\cap H^1_0(\X)$ (which is the space that the strong solution of~(\ref{1}) belongs to), and $c=\beta\cdot\nabla w|_{\partial\X}$, then we see that
\begin{equation}\label{internalregobl}
\jump{c} = \jump{w} = \jump{\nabla w\cdot n_F}=\jump{\nabla_\Ta w} = 0\quad\forall F\in\Ei,
\end{equation} and furthermore, since $\tau_F(\beta\cdot\nabla w) = c$ for all $F\in\Eb$, it follows that 
\begin{equation}\label{externalregobl}
\jump{\beta\cdot\nabla w-c}=\jump{\partial_{{\Ta}}(\beta\cdot\nabla w)} = 0\quad\forall F\in\Eb.
\end{equation}
$\Jhobl$ also enforces the oblique boundary condition, and leads to the bilinear form $\Bhalfobl$ being provably coercive (in a particular $H^2$-type norm on $M_h$). In particular we define $\Jhobl$ as follows:
\begin{equation}\label{Jobldef}
\begin{aligned}
&\Jhobl((u_h,\lambda),(v_h,\mu)) := \sum_{F\in\Ei}\mu_F\langle\jump{\nabla_\Ta u_h},\jump{\nabla_\Ta v_h}\rangle_F]\\
&~~~~~~+\sum_{F\in\Ei}[\mu_F\langle\jump{\nabla u_h\cdot n_F},\jump{\nabla v_h\cdot n_F}\rangle_F
+\eta_F\langle\jump{u_h},\jump{v_h}\rangle_F+\ell_F\langle\jump{\lambda},\jump{\mu}\rangle_F]\\
&~~~~~~+\sum_{F\in\Eb}\sigma_F\langle\beta\cdot\nabla u_h-\lambda,\beta\cdot\nabla v_h-\mu\rangle_F
,\\
\end{aligned}
\end{equation}
where the positive edge-dependent quantities $\mu_F$, $\eta_F$, $\ell_F$, and $\sigma_F$ will be specified later, and their particular choice will be made clear when we prove that $\Bhalfobl$ is coercive (see Lemma~\ref{Coercivity1}). Furthermore,~(\ref{internalregobl}) and~(\ref{externalregobl}) imply that
\begin{equation}\label{consident2}
\Jhobl((w,c),(v_h,\mu)) = 0,\quad\forall(v_h,\mu)\in\dgo\times\mathbb R,
\end{equation}
 when $w\in H^2_{\beta,0}(\X)\cap H^s(\X;\Th)$, $s>5/2$, and $c = \beta\cdot\nabla w|_{\partial\X}$.
 
The bilinear form $\Bstarobl$ plays a key role in identity~(\ref{cons11:20}), and its structure is motivated by Corollary~\ref{cor:applicabletotriangles} and identity~(\ref{star:13:16}), the statements of which we recall.

\emph{Statement of Corollary~\ref{cor:applicabletotriangles}:} Assume that $E\subset\R2$ is a bounded, Lipschitz, piecewise $C^{2}$ domain, and that $\beta\in C^1(\Gamma_n;\mathbb S^1)$ for each $C^{2}$ portion $\Gamma_n$ of $\partial E$, $n=1,\ldots,N$, $N\in\mathbb N$. Then, for any $u,v\in H^s(E)$, $s>5/2$, we have that
\begin{equation}\label{cor3.2.23statement}
\begin{aligned}
\int_ED^2u:D^2v&+\int_{\partial E}\left(\beta_1{\partial_{{\Ta}}}\beta_2-\beta_2{\partial_{{\Ta}}}\beta_1\right)(\betaperp\cdot\nabla u\,\betaperp\cdot\nabla v+\beta\cdot\nabla u\,\beta\cdot\nabla v)\\
&~~~~~~+\int_{\partial E}\left({\partial_{{\Ta}}}(\betaperp\cdot\nabla u)\beta\cdot\nabla v-{\partial_{{\Ta}}}(\beta\cdot\nabla u)\betaperp\cdot\nabla v\right)\\
&~~~~~~~~~~~~=\int_E\Delta u\,\Delta v.
\end{aligned}
\end{equation}
\emph{Identity~(\ref{star:13:16}):}
Let $\X\subset\mathbb R^2$ be a $C^2$ domain, and assume that $\beta\in C^1(\partial\X;\mathbb S^1)$. Then, on $\partial\X$, we have that
\begin{equation}\label{statementident}
\beta_1(\partial_{{\Ta}}\beta_2)-(\partial_{{\Ta}}\beta_1)\beta_2 = \partial_{{\Ta}}\Theta + \mc_{\partial\X}.
\end{equation}
\emph{Designing the bilinear form:} Let us assume that $v_h\in\dgo$, $w\in H^2_{\beta,0}(\X)\cap H^s(\X;\Th)$, $s>5/2$, and $c = \beta\cdot\nabla w|_{\partial\X}$. Furthermore, we will assume that $\X$ is a $C^2$ domain that is also piecewise $C^3$, and that $(\Th)_{h>0}$ is a regular of order $2$ family of triangulations on $\overline{\X}$, satisfying assumption~\ref{Meshconds}. 

Let us consider $K\in\Th$ that satisfies $|\partial K\cap\partial\X|\ne 0$ (this allows for elements with one curved edge that lies on $\partial\X$, but excludes elements that only intersect $\partial\X$ at a vertex). 
Note that $K\subset\R2$ is bounded with a Lipschitz continuous, piecewise $C^3$ boundary. $K$ also has three edges  $F^1_K$, $F^2_K$, $F^3_K$, (each of which are $C^3$ (and hence $C^2$) portions of $\partial K$) and three vertices $e^1_K$, $e^2_K$, $e^3_K$. Let $e^1_K$ and $e^2_K$ be the two vertices that lie on $\partial\X$, and let $F^1_K$ be the \emph{curved} side that lies on $\partial\X$ and also connects $e^1_K$ and $e^2_K$. Finally, let $F^2_K$ be the \emph{straight} edge of $K$ that connects $e^2_K$, and $e^3_K$. It then follows that $F_K^3$ is the remaining \emph{straight} edge that connects $e^3_K$, and $e^1_K$. 

Now define $\tilde{\beta}:\partial K\to\mathbb S^1$ by
\begin{equation*}
\left\{
\begin{aligned}
\tilde{\beta}|_{F_K^1} & = \beta,\\
\tilde{\beta}|_{F_K^2} & = \beta(e^2_K),\\
\tilde{\beta}|_{F_K^3} & = \beta(e^1_K),\\
\end{aligned}
\right.
\end{equation*}
and so $\tilde{\beta}\in C^1(F_K^j;\mathbb S^1)$, $j=1,2,3$, where $\partial K = \cup_{j=1}^3\overline{F_K^j}$.

Then, noting that $w,v_h\in H^s(K)$, $s>5/2$, applying~(\ref{cor3.2.23statement}) with $E:=K$, and $\beta:=\tilde{\beta}$, we obtain
\begin{equation*}
\begin{aligned}
&\int_K D^2_{11}wD^2_{22}v_h+D^2_{22}wD^2_{11}v_h-2D^2_{12}wD^2_{12}v_h = \\
&~~~~~~~\int_{\partial K}\left(\tilde{\beta}_1\partial_{{\Ta}}\tilde{\beta}_2-\tilde{\beta}_2\partial_{{\Ta}}\tilde{\beta}_1\right)({\tilde{\beta}}^\perp\cdot\nabla w\,{\tilde{\beta}}^\perp\cdot\nabla v_h+\tilde{\beta}\cdot\nabla w\,\tilde{\beta}\cdot\nabla v_h)\\
&~~~~~~~~~~~~~+\int_{\partial K}\left(\partial_{{\Ta}}({\tilde{\beta}}^\perp\cdot\nabla w)\tilde{\beta}\cdot\nabla v_h-\partial_{{\Ta}}(\tilde{\beta}\cdot\nabla w){\tilde{\beta}}^\perp\cdot\nabla v_h\right)\\
\end{aligned}
\end{equation*}
\begin{equation}\label{obl:star:2609}
\begin{aligned}
&~~~~~~~=\int_{F^1_K}\left(\beta_1\partial_{{\Ta}}\beta_2-\beta_2\partial_{{\Ta}}\beta_1\right)(\betaperp\cdot\nabla w\,\betaperp\cdot\nabla v_h+\beta\cdot\nabla w\,\beta\cdot\nabla v_h)\\
&~~~~~~~~~~~~~+\int_{F^1_K}\left(\partial_{{\Ta}}(\betaperp\cdot\nabla w)\beta\cdot\nabla v_h-\partial_{{\Ta}}(\beta\cdot\nabla w)\betaperp\cdot\nabla v_h\right)\\
&~~~~~~~~~~~~~+\int_{F^2_K\cup F^3_K}\left(\partial_{{\Ta}}({\tilde{\beta}}^\perp\cdot\nabla w)\tilde{\beta}\cdot\nabla v_h-\partial_{{\Ta}}(\tilde{\beta}\cdot\nabla w){\tilde{\beta}}^\perp\cdot\nabla v_h\right)
\end{aligned}
\end{equation}
Furthermore, upon noting that $\tilde{\beta}$, ${\tilde{\beta}}^\perp$, and the unit normal to $\partial K$, are all constant on $F^2_K$ and $F^3_K$, one can calculate the following:
\begin{equation}\label{obl:starstar:2609}
\begin{aligned}
&\int_{F^2_K\cup F^3_K}\left(\partial_{{\Ta}}({\tilde{\beta}}^\perp\cdot\nabla w)\tilde{\beta}\cdot\nabla v_h-\partial_{{\Ta}}(\tilde{\beta}\cdot\nabla w){\tilde{\beta}}^\perp\cdot\nabla v_h\right)\\
&~~~~~~=\int_{F^2_K\cup F^3_K}\Delta w\,(\nabla v_h\cdot n_{\partial K})-\nabla(\nabla w\cdot n_{\partial K})\cdot\nabla v_h.
\end{aligned}
\end{equation}
Since $F^1_K\subset\partial\X$, we may apply identity~(\ref{statementident}), obtaining
\begin{equation}\label{obl:starstarstar:2609}
\begin{aligned}
&\int_{F^1_K}\left(\beta_1\partial_{{\Ta}}\beta_2-\beta_2\partial_{{\Ta}}\beta_1\right)(\betaperp\cdot\nabla w\,\betaperp\cdot\nabla v_h+\beta\cdot\nabla w\,\beta\cdot\nabla v_h)\\
&~~~~~~=\int_{F^1_K}\left(\partial_{{\Ta}}\Theta+\mc_{F_K^1}\right)(\betaperp\cdot\nabla w\,\betaperp\cdot\nabla v_h+\beta\cdot\nabla w\,\beta\cdot\nabla v_h),
\end{aligned}
\end{equation}
where $\mc_{F_K^1} = \mc_{\partial\X}|_{F_K^1}$.

We now consider an element $K\in\Th$ that satisfies $|\partial K\cap\partial\X| = 0$. An application of integration by parts (noting that the unit outward normal to $\partial K$ is constant on each edge of $K$) yields
\begin{equation}\label{obl:alpha:2609}
\int_K(D^2w_h:D^2v_h)+\int_{\partial K}(\Delta w\,(\nabla v_h\cdot n_{\partial K})-\nabla(\nabla w\cdot n_{\partial K})\cdot\nabla v_h) = \int_K\Delta w\,\Delta v_h.
\end{equation}
One can also see that for any $K\in\Th$, and thus in particular, for those $K\in\Th$ that satisfy $|\partial K\cap\partial\X|\ne0$, 
\begin{equation}\label{obl:beta:2609}
\int_K(D^2w:D^2v_h+D^2_{11}wD^2_{22}v_h+D^2_{22}wD^2_{11}v_h-2D^2_{12}wD^2_{12}v_h) = \int_K\Delta w\,\Delta v_h.
\end{equation}
Applying identities~(\ref{obl:starstar:2609}) and~(\ref{obl:starstarstar:2609}) to~(\ref{obl:star:2609}), and summing~(\ref{obl:alpha:2609}) over all $K\in\Th$ such that $|\partial K\cap\partial\X|=0$, and~(\ref{obl:beta:2609}) over all $K\in\Th$ such that $|\partial K\cap\partial\X|\ne0$, we obtain
\begin{equation}\label{obl:epsilon:2609}
\begin{aligned}
&\sum_{K\in\Th}\int_KD^2w\!:\!D^2v_h +\sum_{F\in\Ei}\int_F\jump{\Delta w\,\nabla v_h\cdot n_F-\nabla(\nabla w\cdot n_F)\cdot\nabla v_h}\\
&~~+\sum_{F\in\Eb}\int_F(\partial_{{\Ta}}(\betaperp\cdot\nabla w)\beta\cdot\nabla v_h-\partial_{{\Ta}}(\beta\cdot\nabla w)\betaperp\cdot\nabla v_h)\\
&~~~~~~~~+\sum_{F\in\Eb}\int_F(\left(\partial_{{\Ta}}\Theta+\mc_F\right)\betaperp\cdot\nabla w\,\betaperp\cdot\nabla v_h+\left(\partial_{{\Ta}}\Theta+\mc_F\right)\beta\cdot\nabla w\,\beta\cdot\nabla v_h)\\
&~~~~~~~~~~~~~~=\sum_{K\in\Th}\int_K\Delta w\,\Delta v_h,
\end{aligned}
\end{equation}
where $n_F$ is now a \emph{fixed} choice of unit normal to $F$, and $\mc_F := \mc_{\partial\X}|_F$.
Utilising the tangential operator identities in~(\ref{identities}), we obtain
\begin{equation}\label{obl:gamma:2609}
\begin{aligned}
&\sum_{F\in\Ei}\int_F\jump{\Delta w\,\nabla v_h\cdot n_F-\nabla(\nabla w\cdot n_F)\cdot\nabla v_h}=\sum_{F\in\Ei}\int_F\jump{\Delta_\Ta w\,\nabla v_h\cdot n_F-\nabla_\Ta(\nabla w\cdot n_F)\cdot\nabla_\Ta v_h}.
\end{aligned}
\end{equation}
We then apply the identity (valid for any $f,g\in H^s(\X;\Th)$, $s>1/2$)
$$\sum_{F\in\Ei}\int_F\jump{fg} = \sum_{F\in\Ei}\int_F\jump{f}\avg{g} +\sum_{F\in\Ei}\int_F\avg{f}\jump{g},$$
along with~(\ref{internalregobl}), to~(\ref{obl:gamma:2609}), which gives us
\begin{equation}\label{12:55}
\begin{aligned}
&\sum_{F\in\Ei}\int_F\jump{\Delta w\,\nabla v_h\cdot n_F-\nabla(\nabla w\cdot n_F)\cdot\nabla v_h}\\
&~~~~~~=\sum_{F\in\Ei}\int_F(\jump{\Delta_\Ta w}\avg{\nabla v_h\cdot n_F}+\avg{\Delta_\Ta w}\jump{\nabla v_h\cdot n_F})\\
&~~~~~~~~~~~~-\sum_{F\in\Ei}\int_F(\jump{\nablaTa(\nabla w\cdot n_F)}\cdot\avg{\nablaTa v_h}+\avg{\nablaTa(\nabla w\cdot n_F)}\cdot\jump{\nablaTa v_h})\\
&~~~~~~=\sum_{F\in\Ei}\int_F(\avg{\Delta_\Ta w}\jump{\nabla v_h\cdot n_F}-\avg{\nablaTa(\nabla w\cdot n_F)}\cdot\jump{\nablaTa v_h}).
\end{aligned}
\end{equation}
Then, again by~(\ref{internalregobl}) we may consistently symmetrise the final right-hand side of~(\ref{12:55}), yielding
\begin{equation}\label{consandsymm13:02}
\begin{aligned}
&\sum_{F\in\Ei}\int_F\jump{\Delta w\,\nabla v_h\cdot n_F-\nabla(\nabla w\cdot n_F)\cdot\nabla v_h}\\
&~~~~~~=\sum_{F\in\Ei}\int_F(\avg{\Delta_\Ta w}\jump{\nabla v_h\cdot n_F}+\avg{\Delta_\Ta v_h}\jump{\nabla w\cdot n_F})\\
&~~~~~~~~~~~~-\sum_{F\in\Ei}\int_F(\avg{\nablaTa(\nabla w\cdot n_F)}\cdot\jump{\nablaTa v_h}+\avg{\nablaTa(\nabla v_h\cdot n_F)}\cdot\jump{\nablaTa w}).
\end{aligned}
\end{equation}
From~(\ref{externalregobl}), we obtain
\begin{equation}\label{obl:deltab:2609}
\sum_{F\in\Eb}\int_F\partial_{{\Ta}}(\beta\cdot\nabla w)\betaperp\cdot\nabla v_h = 0,
\end{equation}
for all $F\in\Eb$. Furthermore, on $F\in\Eb$, one has that 
\begin{equation}\label{innerident13:10}
\betaperp\cdot\nabla w\,\betaperp\cdot\nabla v_h+\beta\cdot\nabla w\,\beta\cdot\nabla v_h = \nabla w\cdot\nabla v_h.
\end{equation}
Applying~(\ref{consandsymm13:02}),~(\ref{obl:deltab:2609}), and~(\ref{innerident13:10}) to~(\ref{obl:epsilon:2609}) we obtain
\begin{equation}\label{obl:epsilon:2609:2}
\begin{aligned}
&\sum_{K\in\Th}\int_K(D^2w\!:\!D^2v_h)+\sum_{F\in\Ei}\int_F(\avg{\Delta_\Ta w}\jump{\nabla v_h\cdot n_F}+\avg{\Delta_\Ta v_h}\jump{\nabla w\cdot n_F})\\
&~~~~~~~~~~~-\sum_{F\in\Ei}\int_F(\avg{\nablaTa(\nabla w\cdot n_F)}\cdot\jump{\nablaTa v_h}+\avg{\nablaTa(\nabla v_h\cdot n_F)}\cdot\jump{\nablaTa w})\\
&~~~~~~~~~~~~~~~~~~~~~+\sum_{F\in\Eb}\int_F(\partial_{{\Ta}}(\betaperp\cdot\nabla w)\beta\cdot\nabla v_h+\left(\partial_{{\Ta}}\Theta+\mc_F\right)\nabla w\cdot\nabla v_h)\\
&~~~~~~~~~~~~~~~~~~~~~~~~~~~~~~~~~~=\sum_{K\in\Th}\int_K\Delta w\,\Delta v_h.
\end{aligned}
\end{equation}
So far, all of the applications of~(\ref{internalregobl}) and~(\ref{externalregobl}) have been made with consistency and symmetry in mind. We make a penultimate alteration, which is necessary for the coercivity of $\Bhalfobl$. In particular, notice that each term of each integrand on the left-hand side of~(\ref{obl:epsilon:2609:2}) either has a sign if we take $w=v_h$ (in particular, $D^2w\!:D^2v_h$ and $(\partial_{{\Ta}}\Theta+\mc_F)\nabla w\cdot\nabla v_h$), or consists of the product of two terms, one of which is present in the definition~(\ref{Jobldef}) of the jump stabilisation bilinear form, $\Jhobl$, except for the integrand $\partial_{{\Ta}}(\betaperp\cdot\nabla w)\beta\cdot\nabla v_h$.

To this end, let us denote by $e_F^+$ and $e_F^-$ the two vertices of an edge $F\in\Eb$, and notice that for any $\mu\in\R{}$,
\begin{equation}\label{obl:deltaa:2609}
\begin{aligned}
\sum_{F\in\Eb}\int_F\partial_{{\Ta}}(\betaperp\cdot\nabla w)\mu & = \sum_{F\in\Eb}(\betaperp\cdot\nabla w)\mu|^{e_F^+}_{e_F^-}\\
& = \mu\sum_{e\in\Vb}\jump{\betaperp\cdot\nabla w} = 0,
\end{aligned}
\end{equation}
where the jumps in~(\ref{obl:deltaa:2609}) are considered across boundary vertices $e\in\Vb$. Note that the final equality holds, due to the fact that $\betaperp\in C^1(\partial\X)$, and $\nabla w\in H^{1/2}(\partial\X)$, and thus neither function may jump across boundary vertices.

Applying~(\ref{obl:deltaa:2609}) to~(\ref{obl:epsilon:2609:2}), we obtain
\begin{equation}\label{obl:epsilon:2609:3}
\begin{aligned}
&\sum_{K\in\Th}\int_K(D^2w\!:\!D^2v_h)+\sum_{F\in\Ei}\int_F(\avg{\Delta_\Ta w}\jump{\nabla v_h\cdot n_F}+\avg{\Delta_\Ta v_h}\jump{\nabla w\cdot n_F})\\
&~~~~~~~~~~~~~~-\sum_{F\in\Ei}\int_F(\avg{\nablaTa(\nabla w\cdot n_F)}\cdot\jump{\nablaTa v_h}+\avg{\nablaTa(\nabla v_h\cdot n_F)}\cdot\jump{\nablaTa w})\\
&~~~~~~~~~~~~~~~~~~+\sum_{F\in\Eb}\int_F(\partial_{{\Ta}}(\betaperp\cdot\nabla w)(\beta\cdot\nabla v_h-\mu)+\left(\partial_{{\Ta}}\Theta+\mc_F\right)\nabla w\cdot\nabla v_h)\\
&~~~~~~~~~~~~~~~~~~=\sum_{K\in\Th}\int_K\Delta w\,\Delta v_h.
\end{aligned}
\end{equation}
Alas, the left-hand side of~(\ref{obl:epsilon:2609:3}) is not symmetric. However, by~(\ref{externalregobl}), it follows that
\begin{equation}\label{obl:deltab:2609:cd}
\sum_{F\in\Eb}\left\langle\partial_{{\Ta}}(\betaperp\cdot\nabla v_h),\beta\cdot\nabla w-c\right\rangle_F = 0,
\end{equation}
which, when applied to~(\ref{obl:epsilon:2609:3}), gives us
\begin{equation}\label{obl:epsilon:2609:4}
\begin{aligned}
&\sum_{K\in\Th}\int_K(D^2w\!:\!D^2v_h)+\sum_{F\in\Ei}\int_F(\avg{\Delta_\Ta w}\jump{\nabla v_h\cdot n_F}+\avg{\Delta_\Ta v_h}\jump{\nabla w\cdot n_F})\\
&~~~~~~~~~~~~~~-\sum_{F\in\Ei}\int_F(\avg{\nablaTa(\nabla w\cdot n_F)}\cdot\jump{\nablaTa v_h}+\avg{\nablaTa(\nabla v_h\cdot n_F)}\cdot\jump{\nablaTa w})\\
&~~~~~~~~~~~~~~~~~~+\sum_{F\in\Eb}\int_F(\partial_{{\Ta}}(\betaperp\cdot\nabla v_h)(\beta\cdot\nabla w-c)+\partial_{{\Ta}}(\betaperp\cdot\nabla w)(\beta\cdot\nabla v_h-\mu))\\
&~~~~~~~~~~~~~~~~~~~~+\sum_{F\in\Eb}\int_F(\left(\partial_{{\Ta}}\Theta+\mc_F\right)\nabla w\cdot\nabla v_h)\\
&~~~~~~~~~~~~~~~~~~=\sum_{K\in\Th}\int_K\Delta w\,\Delta v_h,
\end{aligned}
\end{equation}
consistently restoring symmetry.

We then define $\Bstarobl$ by the left-hand side of~(\ref{obl:epsilon:2609:4}). That is,
\begin{equation}\label{obl:B_h*def}
\begin{aligned}
&\Bstarobl((u_h,\lambda),(v_h,\mu)) := \sum_{K\in\Th}\langle D^2u_h,D^2v_h\rangle_K\\
&~~~~~~+\sum_{F\in\Ei}[\langle\divTa\gradTa\avg{u_h},\jump{\nabla v_h\cdot n_F}\rangle_F+\langle\divTa\gradTa\avg{v_h},\jump{\nabla u_h\cdot n_F}\rangle_F]\\
&~~~~~~-\sum_{F\in\Ei}[\langle\gradTa\avg{\nabla u_h\cdot n_F},\jump{\gradTa v_h}\rangle_F+\langle\gradTa\avg{\nabla v_h\cdot n_F},\jump{\gradTa u_h}\rangle_F]\\
&~~~~~~+\sum_{F\in\Eb}\left[\left\langle\left(\partial_{{\Ta}}\Theta+\mc_F\right)\nabla u_h,\nabla v_h\right\rangle_F\right]\\
&~~~~~~+\sum_{F\in\Eb}\left[\left\langle\partial_{{\Ta}}(\betaperp\cdot\nabla u_h),\beta\cdot\nabla v_h-\mu\right\rangle_F+\left\langle\partial_{{\Ta}}(\betaperp\cdot\nabla v_h),\beta\cdot\nabla u_h-\lambda\right\rangle_F\right],
\end{aligned}
\end{equation}
for all $(u_h,\lambda),(v_h,\mu)\in M_h$, and we recall that $n_F$ is a fixed choice of unit normal to $F$, and $\mc_F:=\partial_\Ta\cdot n_F=\mc_{\partial\X}|_F$ for $F\in\Eb$. It follows from~(\ref{obl:epsilon:2609:4}) that the bilinear form $\Bstarobl$ satisfies the first identity of~(\ref{cons13:46}). We are now ready to define the numerical method of this chapter.
\endgroup
\subsection{Finite element method}
Recall that we define $M_h:=\dgo\times\dgoo$. The definition of the finite element method first requires the definition of the jump stabilisation bilinear form
$J_h:M_h\times M_h\to\mathbb R$, given by
\begin{equation*}
\begin{aligned}
&J_h((u_h,\lambda),(v_h,\mu)) := \sum_{F\in\Ei}\mu_F\langle\jump{\nabla_\Ta u_h},\jump{\nabla_\Ta v_h}\rangle_F]\\
&~~~~~~+\sum_{F\in\Ei}[\mu_F\langle\jump{\nabla u_h\cdot n_F},\jump{\nabla v_h\cdot n_F}\rangle_F
+\eta_F\langle\jump{u_h},\jump{v_h}\rangle_F+\ell_F\langle\jump{\lambda},\jump{\mu}\rangle_F]\\
&~~~~~~+\sum_{F\in\Eb}\sigma_F\langle\beta\cdot\nabla u_h-\lambda,\beta\cdot\nabla v_h-\mu\rangle_F
,\\
\end{aligned}
\end{equation*}
where $\mu_F,\eta_F,\ell_F,\sigma_F$ are positive, face dependent parameters to be provided. Furthermore for $\theta\in[0,1]$, we define the bilinear form $B_{h,\theta}:M_h\times M_h\to\mathbb R$
\begin{equation*}
\begin{aligned}
&B_{h,\theta}((u_h,\lambda),(v_h,\mu)):=\theta B_{h,*}((u_h,\lambda),(v_h,\mu))+(1-\theta)\sum_{K\in\Th}\langle\Delta u_h,\Delta v_h\rangle_K\\
&~~~~~~~+J_h((u_h,\lambda),(v_h,\mu)),
\end{aligned}
\end{equation*} where $B_{h,*}$ is given by~(\ref{obl:B_h*def}). We now define
\begin{equation}\label{A:def:15:59}
A_h((u_h,\lambda),(v_h,\mu)):=\sum_{K\in\Th}[\langle\gamma Lu_h,\Delta v_h\rangle_K-\langle\Delta u_h,\Delta v_h\rangle_K]+B_{h,1/2}((u_h,\lambda),(v_h,\mu)).
\end{equation}
We can now state the finite element method: find $(u_h,{c}_h)\in M_h$ such that
\begin{equation}\label{FEM}
A_h((u_h,c_h),(v_h,\mu)) = \sum_{K\in\Th}\langle\gamma f,\Delta v_h\rangle_K\,\,\quad\forall (v_h,\mu)\in M_h.
\end{equation}
\subsection{Consistency of the method}
\begin{remark}[Extension of the bilinear forms]
The bilinear forms $B_{h,*}$ and $J_h$ are both defined on $M_h\times M_h$, but one must note that they are both well defined on $(H^s(\X;\Th)\cap H^2_{\beta,0}(\X)\times\dgoo)\times M_h$ for $s>5/2$, of which $(H^s(\X;\Th)\cap H^2_{\beta,0}(\X)\times\mathbb R)\times(\dgo\times\mathbb R)$ is a proper subset, that the functions in the following lemma belong to.
\end{remark}
\begingroup\color{black}
\begin{lemma}\label{obl:cons:lemma}
Let $\X\subset\mathbb R^2$ be a $C^2$ and piecewise $C^3$ domain, and let $\beta\in C^1(\partial\X;\mathbb S^1)$. Furthermore, assume that $\{\Th\}_{h}$ is a regular of order $2$ family of triangulations on $\overline{\X}$ satisfying Assumption~\ref{Meshconds}.
Let $(w,c)\in H^s(\X;\Th)\cap H^2_{\beta}(\X)\times\R{},$ $s>5/2$, where $\beta\cdot\nabla w|_{\partial\X}=c$.
Then, for every $(v_h,\mu)\in \dg\times\R{}$, we have the identities
\begin{equation}\label{obl:Cons:identity}
\Bstarobl((w,c),(v_h,\mu))=\sum_{K\in\Th}\langle\Delta w,\Delta v_h\rangle_K\quad\mbox{and}\quad \Jhobl((w,c),(v_h,\mu))=0.
\end{equation}
\end{lemma}
%
\emph{Proof:} Assume that the pair $(w,c)$ satisfies the hypotheses of the lemma. Then, the identities of~(\ref{obl:Cons:identity}) follow from~(\ref{obl:epsilon:2609:4}) and~(\ref{consident2}).$\quad\quad\square$
\endgroup

\subsection{Stability of the method}
We now aim to show that $\Bhtheobl$ is coercive in a particular norm on $M_h$.
Before we prove that $\Bhtheobl$ is coercive, we must define the norm in which the bilinear form is coercive. To this end, let us define the following family of functionals, $\|(\cdot,\cdot)\|_{h,\theta}:M_h\to[0,\infty)$ for $\theta\in(0,1]$:
\begin{equation}\label{oblique:dg:norm}
\begin{aligned}
\|(u_h,\lambda)\|_{h,\theta}^2&:=\sum_{K\in\Th}[\theta|u_h|_{H^2(K)}^2+(1-\theta)\|\Delta u_h\|_{L^2(K)}^2]\\
&~~~~~~~~~~~+c_*\Jhobl((u_h,\lambda),(u_h,\lambda))+\frac{\theta}{2}\sum_{F\in\Eb}\left\|\left(\partial_{{\Ta}}\Theta+\mc_F\right)^{1/2}\nabla u_h\right\|_{L^2(F)}^2,
\end{aligned}
\end{equation}
where $c_*$ is a positive constant to be determined.
\begin{lemma}\label{obl:normproofs}
Let $\X\subset\mathbb R^2$ be a $C^2$ and piecewise $C^3$ domain, and let $\beta\in C^1(\partial\X;\mathbb S^1)$. Assume that
$$\partial_{{\Ta}}\Theta+\mc_{\partial\X}>0\quad\mbox{on }\partial\X.$$
Furthermore, assume that $\{\Th\}_{h}$ is a regular of order $2$ family of triangulations on $\overline{\X}$ satisfying Assumption~\ref{Meshconds}. Then, for each $\theta\in(0,1]$, $\|\cdot\|_{h,\theta}:M_h\to[0,\infty)$ defines a norm on $M_h$.
\end{lemma}
\emph{Proof:} First we note that homogeneity and the triangle inequality are clear. Now let us assume that the pair $(v_h,\mu)\in M_h$ satisfies
$$\|(v_h,\mu)\|_{h,\theta} = 0,$$
for some $\theta\in(0,1]$. It then follows that $|v_h|_{H^2(\X;\Th)}=0$, and so $v_h$ is piecewise affine. Moreover
$$\jump{\mu}_F=\jump{v_h}_F=\jump{\nabla v_h}_F = 0\,\,\mbox{for}\,\,F\in\Ei,$$
and, as $\partial_{{\Ta}}\Theta+\mc_{\partial\X}>0$ on $\partial\X$, it follows that
$$\jump{\nabla v_h}_F = 0\,\,\mbox{for}\,\,F\in\Eb.$$
It then follows that $v_h$ is affine, i.e, $v_h=a^Tx+b$, with $a\in\mathbb R^d$, $b\in\mathbb R$, and that $\mu$ is constant. But then we see that
$$0=\nabla v_h|_F = a,$$
for $F\in\Eb$, and thus $a=0$, i.e., $v_h = b$.  Then, since $v_h\in\dgo$, $0=\int_\X v_h=|\X|b$, and so $b=0$, i.e., $v_h\equiv0$.

Finally, we see that $\Jhobl((v_h,\mu),(v_h,\mu)) = 0$, and it follows that
$$0=\beta\cdot\nabla v_h = \mu\,\,\mbox{on}\,\,\partial\X,$$
and so $\mu\equiv 0$.
Overall, we have obtained $(v_h,\mu)\equiv(0,0)$. $\quad\quad\square$
\begin{lemma}\label{Coercivity1}
Let $\X\subset\mathbb R^2$ be a $C^2$ and piecewise $C^3$ domain, and let $\beta\in C^1(\partial\X;\mathbb S^1)$. Assume that
\begin{equation}\label{16:04:23}
\partial_{{\Ta}}\Theta+\mc_{\partial\X}>0\quad\mbox{on }\partial\X.
\end{equation}
Furthermore, assume that $\{\Th\}_{h}$ is a regular of order $2$ family of triangulations on $\overline{\X}$ satisfying Assumption~\ref{Meshconds}.
Then, for
each constant $\kappa>1$, there exist positive constants $c_{\operatorname{stab}}$ and $c_*$, independent of $h$, and $\theta$, such that
\begin{equation}\label{obl:4.4}
\Bhtheobl((u_h,\lambda),(u_h,\lambda))\ge\kappa^{-1}\|(u_h,\lambda)\|_{h,\theta}^2\,\,\,\forall (u_h,\lambda)\in M_h,\,\forall\theta\in[0,1],
\end{equation}
whenever \begin{equation}\label{obl:4.5}
\mu_F\ge \frac{c_{\operatorname{stab}}}{\tilde{h}_F},\sigma_F\ge\frac{c_{\operatorname{stab}}}{\tilde{h}_F}\,\,\mbox{and}\,\,\eta_F,\ell_F>0.
\end{equation}
\end{lemma}
\emph{Proof:} We see that for $(u_h,\lambda)\in\dg\times\dgoo$,
\begin{equation*}
\begin{aligned}
&B_{h,\theta}((u_h,\lambda),(u_h,\lambda))  = \sum_{K\in\Th}[\theta\langle D^2u_h,D^2u_h\rangle_K+(1-\theta)\langle\Delta u_h,\Delta u_h\rangle_K]\\
&~~~~~~+2\theta\sum_{F\in\Ei}[\langle\divTa\nablaTa\avg{u_h},\jump{\nabla u_h\cdot n_F}\rangle_F-\langle\nablaTa\avg{\nabla u_h\cdot n_F},\jump{\nablaTa u_h}\rangle_F]\\
&~~~~~~+\theta\sum_{F\in\Eb}\left[\left\|(\partial_\Ta\Theta+\mc_F)^{1/2}\betaperp\cdot\nabla u_h\right\|_{L^2(F)}^2+\left\|(\partial_\Ta\Theta+\mc_F)^{1/2}\beta\cdot\nabla u_h\right\|_{L^2(F)}^2\right]\\
&~~~~~~+2\theta\sum_{F\in\Eb}\langle\partial_\Ta(\betaperp\cdot\nabla u_h),\beta\cdot\nabla u_h-\lambda\rangle_F\\
&~~~~~~+\sum_{F\in\Ei}[\mu_F\|\jump{\nablaTa u_h}\|_{L^2(F)}^2+\mu_F\|\jump{\nabla u_h\cdot n_F}\|_{L^2(F)}^2+\eta_F\|\jump{u_h}\|_{L^2(F)}^2+\ell_F\|\jump{\lambda}\|_{L^2(F)}^2]\\
&~~~~~~+\sum_{F\in\Eb}\sigma_F\|\beta\cdot\nabla u_h-\lambda\|_{L^2(F)}^2.
\end{aligned}
\end{equation*}
Now notice that for any $\alpha>0$,
\begin{equation}\label{I1:15:12}
\begin{aligned}
|I_1|&:=\left|\sum_{F\in\Ei}\langle\operatorname{div}_\Ta\nabla_\Ta\avg{u_h},\jump{\nabla u_h\cdot n_F}\rangle_F\right|\\
&\le\left(\sum_{F\in\Ei}\alpha\tilde{h}_F\|\operatorname{div}_\Ta\nabla_\Ta\avg{u_h}_F\|_{L^2(F)}^2\right)^{1/2}\left(\sum_{F\in\Ei}\frac{1}{\alpha\tilde{h}_F}\|\jump{\nabla u_h\cdot n_F}\|_{L^2(F)}^2\right)^{1/2},
\end{aligned}
\end{equation}
and (associating $F=\overline{K}\cap\overline{K'}$ for some $K,K'\in\Th$)
$$\|\operatorname{div}_\Ta\nabla_\Ta\avg{u_h}_F\|_{L^2(F)}^2\le\frac{1}{2}\|\operatorname{div}_\Ta\nabla_\Ta u_h|_{K}\|_{L^2(F)}^2+\frac{1}{2}\|\operatorname{div}_\Ta\nabla_\Ta u_h|_{K'}\|_{L^2(F)}^2.$$
\begingroup\color{black}Therefore, the trace inequality gives us
\begin{equation*}
\begin{aligned}
\frac{\alpha}{2}\sum_{F\in\Ei}\tilde{h}_F\|\operatorname{div}_\Ta\nabla_\Ta\avg{u_h}_F\|_{L^2(F)}^2&\le\frac{\alpha}{2}\sum_{F\in\Ei}\tilde{h}_F\sum_{K\in\Th:F\subset\partial K}\|D^2u_h\|_{L^2(\partial K)}^2\\
&\le\frac{\alpha C_{\operatorname{Tr}}}{2}\sum_{F\in\Ei}\tilde{h}_F\sum_{K\in\Th:F\subset\partial K}h_K^{-1}|u_h|_{H^2(K)}^2+h_K|u_h|_{H^3(K)}^2,\\
\end{aligned}
\end{equation*}
where $C_{\operatorname{Tr}}$ is the constant of the trace inequality, and is independent of $K$ and $h_K$. We now apply an inverse estimate above, noting that since a given function of the finite element space is a polynomial \emph{composed} with a nonaffine function $F_K$, by the chain rule, the inverse estimate takes the following form
$$|u_h|_{H^3(K)}^2\le C_Ih_K^{-2}(|u_h|_{H^2(K)}^2+|u_h|_{H^1(K)}^2),$$
where $C_I$ a constant independent of $K$ and $h_K$. This results in
\begin{equation*}
\begin{aligned}
\frac{\alpha}{2}\sum_{F\in\Ei}\tilde{h}_F\|\operatorname{div}_\Ta\nabla_\Ta\avg{u_h}\|_{L^2(F)}^2&\le\frac{\alpha C_{\operatorname{Tr}}C_I}{2}\sum_{F\in\Ei}\tilde{h}_F\sum_{K\in\Th:F\subset\partial K}h_K^{-1}(|u_h|_{H^2(K)}^2+|u_h|_{H^1(K)}^2),\\
&\le\frac{\alpha C_{\operatorname{Tr}}C_IC_{\mathcal{F}}}{2}\sum_{K\in\Th}|u_h|_{H^2(K)}^2+|u_h|_{H^1(K)}^2,
\end{aligned}
\end{equation*}
where the final inequality is due to the fact that the number of faces that make up a simplex $K\in\Th$ is bounded by $C_{\mathcal{F}}$. Applying the above estimate to~(\ref{I1:15:12}), we obtain
\begin{equation}\label{I1:15:12:2}
|I_1|\le \frac{\alpha C_{\operatorname{Tr}}C_IC_{\mathcal{F}}}{2}\sum_{K\in\Th}|u_h|_{H^2(K)}^2+|u_h|_{H^1(K)}^2+\sum_{F\in\Ei}\frac{1}{2\alpha\tilde{h}_F}\|\jump{\nabla u_h\cdot n_F}\|_{L^2(F)}^2.
\end{equation}
Similarly, for any $\alpha>0$
\begin{equation}
\begin{aligned}
|I_2|&:=\left|\sum_{F\in\Ei}\langle\nabla_\Ta\avg{\nabla u_h\cdot n_F},\jump{\nabla_\Ta u_h}\rangle_F\right|\\
&\le\frac{\alpha C_{\operatorname{Tr}}C_IC_{\mathcal{F}}}{2}\sum_{K\in\Th}|u_h|_{H^2(K)}^2+|u_h|_{H^1(K)}^2+\sum_{F\in\Ei}\frac{1}{2\alpha\tilde{h}_F}\|\jump{\nabla_\Ta u_h}\|_{L^2(F)}^2.
\end{aligned}
\end{equation}
Since $\betaperp\in C^1(\partial\X;\mathbb S^1)$, utilising the trace and inverse inequalities, we also see that for any $\alpha>0$, 
\begin{equation}\label{I3star}
\begin{aligned}
|I_3|&:=\left|\sum_{F\in\Eb}\langle\partial_\Ta(\betaperp\cdot\nabla u_h),\beta\cdot\nabla u_h-\lambda\rangle_F\right|\\
&\le \sum_{F\in\Eb}C(\beta)(\|\nabla u_h\|_{L^2(F)}+\|D^2u_h\|_{L^2(F)})\|\beta\cdot\nabla u_h-\lambda\|_{L^2(F)}\\
&\le\sum_{F\in\Eb}\frac{C(\beta)^2}{2\alpha\tilde{h}_F}\|\beta\cdot\nabla u_h-\lambda\|_{L^2(F)}^2+\frac{\alpha\tilde{h}_F}{2}\sum_{K\in\Th:F\subset\partial K}|u_h|_{H^2(\partial K)}^2+|u_h|_{H^1(\partial K)}^2\\
&\le \sum_{F\in\Eb}\left[\frac{C(\beta)^2}{2\alpha\tilde{h}_F}\|\beta\cdot\nabla u_h-\lambda\|_{L^2(F)}^2\right]+\frac{\alpha C_{\operatorname{Tr}}C_IC_{\mathcal{F}}}{2}\sum_{K\in\Th}|u_h|_{H^2(K)}^2+|u_h|_{H^1(K)}^2.
\end{aligned}
\end{equation}
Now, by the discrete Poincar\'e--Friedrichs' inequality of~\cite{MR1974504}, we obtain
\begin{equation*}
\begin{aligned}
\sum_{K\in\Th}|u_h|_{H^1(K)}^2\le C(\sigma)\left(\sum_{K\in\Th}|u_h|_{H^2(K)}^2+\sum_{F\in\Ei}\tilde{h}_F^{-1}\|\jump{\nabla u_h}\|_{L^2(F)}^2+\frac{1}{|\partial\X|}\sum_{F\in\Eb}\|\nabla u_h\|_{L^2(F)}^2\right)
\end{aligned}
\end{equation*}
\begin{equation}\label{poinbound}
\begin{aligned}
\le C(\sigma)\left(\sum_{K\in\Th}|u_h|_{H^2(K)}^2+\sum_{F\in\Ei}\tilde{h}_F^{-1}\|\jump{\nabla u_h}\|_{L^2(F)}^2+\frac{\Theta_*}{|\partial\X|}\sum_{F\in\Eb}\|(\partial_\Ta\Theta+\mc_F)^{1/2}\nabla u_h\|_{L^2(F)}^2\right)
\end{aligned}
\end{equation}
where $C(\sigma)$ depends only upon the shape-regularity constant of the family of meshes $(\Th)_h$, and the final inequality follows from the following observation:  for any $F\in\Eb$,
\begin{equation*}
\begin{aligned}
\|\nabla u\|_{L^2(F)}^2 & = \int_F\frac{1}{\partial_\Ta\Theta+\mc_F}(\partial_\Ta\Theta+\mc_F)^{1/2}|\nabla u|^2\\
&\le\left(\min_{F\in\Eb}\inf_F(\partial_\Ta\Theta+\mc_F)\right)^{-1}\|(\partial_\Ta\Theta+\mc_F)^{1/2}\nabla u\|_{L^2(F)}^2,
\end{aligned}
\end{equation*}
where $\Theta_*:=(\min_{F\in\Eb}\inf_F(\partial_\Ta\Theta+\mc_F))^{-1}$ is uniformly positive, due to~(\ref{16:04:23}).

Applying~(\ref{poinbound}) to~(\ref{I1:15:12:2})-(\ref{I3star}), and summing the resulting inequalities, we obtain
\begin{equation*}
\begin{aligned}
\sum_{i=1}^3|I_i|&\le\frac{3\alpha C_{\operatorname{Tr}}C_IC_{\mathcal{F}}(1+C(\sigma))}{2}\sum_{K\in\Th}\|D^2u_h\|_{L^2(K)}^2\\
&~~~+\sum_{F\in\Ei}\tilde{h}_F^{-1}\left(\frac{1}{2\alpha}+\frac{3\alpha C_{\operatorname{Tr}}C_IC_{\mathcal{F}}C(\sigma)}{2}\right)\left(\|\jump{\nabla _\Ta u_h}\|_{L^2(F)}^2+\|\jump{\nabla u_h\cdot n_F}\|_{L^2(F)}^2\right)\\
&~~~+\sum_{F\in\Eb}\frac{C(\beta)^2}{2\alpha\tilde{h}_F}\|\beta\cdot\nabla u_h-\lambda\|_{L^2(F)}^2+\frac{3\alpha C_{\operatorname{Tr}}C_IC_{\mathcal{F}}C(\sigma)\Theta_*}{2|\partial\X|}\|\nabla u\|_{L^2(F)}^2.
\end{aligned}
\end{equation*}
The above estimate implies that
$B_{h,\theta}((u_h,\lambda),(u_h,\lambda))\ge\sum_{i=1}^7A_i,$
where
\begin{equation*}
\begin{aligned}
&A_1:=\theta(1-3\alpha C_{\operatorname{Tr}}C_IC_{\mathcal{F}}(1+C(\sigma)))\sum_{K\in\Th}\|D^2u_h\|_{L^2(K)}^2,\,\,A_2:=(1-\theta)\sum_{K\in\Th}\|\Delta u_h\|_{L^2(K)}^2,\\
&A_3:=\sum_{F\in\Ei}\left(\mu_F-\frac{\theta}{\tilde{h}_F}(\alpha^{-1}+3\alpha C_{\operatorname{Tr}}C_IC_{\mathcal{F}}C(\sigma))\right)\|\jump{\nabla u_h\cdot n_F}\|_{L^2(F)}^2,\\
&A_4:=\sum_{F\in\Ei}\eta_F\|\jump{u_h}\|_{L^2(F)}^2+\ell_F\|\jump{\lambda}\|_{L^2(F)}^2,\\
&A_5:=\sum_{F\in\Ei}\left(\mu_F-\frac{\theta}{\tilde{h}_F}(\alpha^{-1}+3\alpha C_{\operatorname{Tr}}C_IC_{\mathcal{F}}C(\sigma))\right)\|\jump{\nabla_\Ta u_h}\|_{L^2(F)}^2,\\
&A_6:=\sum_{F\in\Eb}\left(\sigma_F-\frac{\theta C(\beta)^2}{\alpha\tilde{h}_F}\right)\|\beta\cdot\nabla u_h-\lambda\|^2_{L^2(F)},\\
&A_7:=\theta\left(1-\frac{3\alpha C_{\operatorname{Tr}}C_IC_{\mathcal{F}}C(\sigma)\Theta_*}{|\partial\X|}\right)\sum_{F\in\Eb}\|(\partial_\Ta\Theta+\mc_F)^{1/2}\nabla u_h\|_{L^2(F)}^2.
\end{aligned}
\end{equation*}\endgroup
Now let $\kappa>1$ be given. Then, since $\kappa^{-1}<1$, there exists $\alpha>0$ sufficiently small such that
$$\min\left\{1-3\alpha C_{\operatorname{Tr}}C_IC_{\mathcal{F}}(1+C(\sigma)),1-\frac{3\alpha C_{\operatorname{Tr}}C_IC_{\mathcal{F}}C(\sigma)\Theta_*}{|\partial\X|}\right\}>\kappa^{-1},$$
we then choose $c_{\operatorname{stab}}:=2\max\{\alpha^{-1}+3\alpha C_{\operatorname{Tr}}C_IC_{\mathcal{F}}C(\sigma),C(\beta)^2\alpha^{-1}\}$, $c_*=\kappa/2$ and note that by assumption, $\mu_F\ge c_{\operatorname{stab}}/\tilde{h}_F$ and $\sigma_F\ge c_{\operatorname{stab}}/\tilde{h}_F$. Therefore, for any $\theta\in[0,1]$, 
\begin{equation*}
\begin{aligned}
A_3&\ge\frac{1}{2}\sum_{F\in\Ei}\mu_F\|\jump{\nabla u_h\cdot n_F}_F\|_{L^2(F)}^2=\kappa^{-1}c_*\sum_{F\in\Ei}\|\jump{\nabla u_h\cdot n_F}_F\|_{L^2(F)}^2,\\
A_4&\ge\frac{1}{2}A_4=\kappa^{-1}c_*\sum_{F\in\Ei}\eta_F(\|\jump{u_h}\|_{L^2(F)}^2+\|\jump{\lambda}\|_{L^2(F)}^2),\\
A_5&\ge\frac{1}{2}\sum_{F\in\Ei}\mu_F\|\jump{\nabla_\Ta u_h}\|_{L^2(F)}^2=\kappa^{-1}c_*\sum_{F\in\Ei}\mu_F\|\jump{\nabla_\Ta u_h}_F\|_{L^2(F)}^2\\
A_6&\ge\frac{1}{2}\sum_{F\in\Eb}\sigma_F\|\beta\cdot\nabla u_h-\lambda\|_{L^2(F)}^2=\kappa^{-1}c_*\sum_{F\in\Eb}\sigma_F\|\beta\cdot\nabla u_h-\lambda\|^2_{L^2(F)}.
\end{aligned}
\end{equation*}
Thus, we obtain
\begin{equation*}
\begin{aligned}
\kappa B_{h,\theta}((u_h,\lambda),(u_h,\lambda))&\ge\sum_{K\in\Th}[\theta\|D^2u_h\|_{L^2(K)}^2+(1-\theta)\|\Delta u_h\|_{L^2(K)}^2]\\
&~~~~~+c_*J_h((u_h,\lambda),(u_h,\lambda))+\theta\sum_{F\in\Eb}\|(\partial_\Ta\Theta+\mc_F)^{1/2}\nabla u_h\|_{L^2(F)}^2.\quad\quad\square
\end{aligned}
\end{equation*}
We will now prove that $A_h$ is coercive in $\|\cdot\|_{h,1}$.
\begin{theorem}\label{Coercivity2}
Under the assumptions of Lemma~\ref{Coercivity1}, let $c_{\operatorname{stab}}$ and $c_*$, $\mu_F$, $\eta_F$, $\sigma_F$, and $\ell_F$ be chosen so that~(\ref{obl:4.4}) and~(\ref{obl:4.5}) hold with $\kappa < (1 - \varepsilon)^{-1/2}$. Let the operator $L$ be uniformly elliptic (and thus satisfy the Cordes condition~(\ref{2})). Then, the operator $A_h$ is coercive in $\|\cdot\|_{h,1}$. In particular, for any $(v_h,\mu)\in M_h$, we have
\begin{equation}\label{stability:est}
\|(v_h,\mu)\|_{h,1}^2\le\frac{2\kappa}{1-\kappa^2(1-\varepsilon)}A_h((v_h,\mu),(v_h,\mu)).
\end{equation}

Therefore, there exists a unique solution pair $(u_h,c_h)\in M_h$ of the numerical scheme~(\ref{FEM}). Moreover, the pair $(u_h,c_h)$ satisfies
\begin{equation}\label{unifbound}
\|(u_h,c_h)\|_{h,1}\le\frac{2\sqrt{2}\kappa\|\gamma\|_{L^\infty(\Omega)}}{1-\kappa^2(1-\varepsilon)}\|f\|_{L^2(\Omega)}.
\end{equation}
\end{theorem}
\emph{Proof:} Let $(v_h,\mu)\in M_h$, then we have that for any $K\in\Th$:
\begin{equation*}
\begin{aligned}
\langle\gamma Lv_h-\Delta v_h,\Delta v_h\rangle_K & \le \|(\gamma L-\Delta)v_h\|_{L^2(K)}\|\Delta v_h\|_{L^2(K)}\\
& \le\sqrt{1-\varepsilon}\|D^2v_h\|_{L^2(K)}\|\Delta v_h\|_{L^2(K)}.
\end{aligned}
\end{equation*}
Applying the Cauchy--Schwarz inequality with a parameter, and~(\ref{obl:4.4}) then gives us
\begin{equation*}
\begin{aligned}
&A_h((v_h,\mu);(v_h,\mu)) = B_{h,1/2}((v_h,\mu);(v_h,\mu))+\sum_{K\in\Th}\langle\gamma Lv_h-\Delta v_h,\Delta v_h\rangle_K+\|\Delta v_h\|_{L^2(K)}^2\\
&~\ge\kappa^{-1}\|(v_h,\mu)\|_{h,1/2}^2+\sum_{K\in\Th}[\|\Delta v_h\|_{L^2(K)}^2-\sqrt{1-\varepsilon}\|D^2v_h\|_{L^2(K)}\|\Delta v_h\|_{L^2(K)}]\\
&~\ge\kappa^{-1}\left(\sum_{K\in\Th}\frac{1}{2}\|D^2v_h\|_{L^2(K)}^2+\frac{1}{2}\|\Delta v_h\|_{L^2(K)}^2\right.\\
&~~~~~~~~~+\left.c_*J_h((v_h,\mu),(v_h,\mu))+\frac{1}{2}\sum_{F\in\Eb}\|(\partial_\Ta\Theta+\mc_F)^{1/2}\nabla u_h\|_{L^2(F)}^2\right)\\
&~~~~~~~~~~~~~~~~-\sum_{K\in\Th}\left[\frac{\kappa(1-\varepsilon)}{2}\|D^2v_h\|_{L^2(K)}^2+\frac{\kappa^{-1}}{2}\|\Delta v_h\|_{L^2(K)}^2\right]\\
& ~= \kappa^{-1}\sum_{K\in\Th}\frac{1-\kappa^2(1-\varepsilon)}{2}\|D^2v_h\|_{L^2(K)}^2\\
&~~~~~~~~~~~~~~~~~~~~~~~~+\kappa^{-1}\left(c_*J_h((v_h,\mu),(v_h,\mu))+\frac{1}{2}\sum_{F\in\Eb}\|(\partial_\Ta\Theta+\mc_F)^{1/2}\nabla u_h\|_{L^2(F)}^2\right)\\
&~\ge\frac{1-\kappa^2(1-\varepsilon)}{2\kappa}\left[\sum_{K\in\Th}\|D^2v_h\|_{L^2(K)}^2+c_*J_h((v_h,\mu),(v_h,\mu))+\frac{1}{2}\sum_{F\in\Eb}\|(\partial_\Ta\Theta+\mc_F)^{1/2}\nabla u_h\|_{L^2(F)}^2\right]\\
&~=\frac{1-\kappa^2(1-\varepsilon)}{2\kappa}\|(v_h,\mu)\|_{h,1}^2,
\end{aligned}
\end{equation*}
where the penultimate inequality follows from the fact that $1>(1-\kappa^2(1-\varepsilon))/2>0$.Thus, we obtain
\begin{equation}\label{cbound}
\|(v_h,\mu)\|_{h,1}^2\le\frac{2\kappa}{1-\kappa^2(1-\varepsilon)}A_h((v_h,\mu),(v_h,\mu)).
\end{equation}
By Lemma~\ref{obl:normproofs}, $\|\cdot\|_{h,1}$ is a norm on $M_h$, and so it follows that there exists a unique pair $(u_h,c_h)\in M_h$ such that
$$A_h((u_h,c_h),(v_h,\mu)) = \sum_{K\in\Th}\langle\gamma f,\Delta v_h\rangle_K\,\,\quad\forall ({v}_h,\mu)\in M_h.$$

Finally, taking $(v_h,\mu)=(u_h,c_h)$ in~(\ref{cbound}) gives us:
\begin{equation*}
\begin{aligned}
\|(u_h,c_h)\|_{h,1}^2&\le\frac{2\kappa}{1-\kappa^2(1-\varepsilon)}A_h((u_h,c_h),(u_h,c_h))\\
& = \frac{2\kappa}{1-\kappa^2(1-\varepsilon)}\sum_{K\in\Th}\langle\gamma f,\Delta u_h\rangle_K\\
&\le \frac{2\kappa\|\gamma\|_{L^\infty(\Omega)}}{1-\kappa^2(1-\varepsilon)}\sum_{K\in\Th}\|f\|_{L^2(K)}\|\Delta u_h\|_{L^2(K)}\\
&\le \frac{2\sqrt{2}\kappa\|\gamma\|_{L^\infty(\Omega)}}{1-\kappa^2(1-\varepsilon)}\|f\|_{L^2(K)}\|(u_h,c_h)\|_{h,1},
\end{aligned}
\end{equation*}
note that the factor of $\sqrt{2}$ comes from the fact that $\|\Delta u_h\|_{L^2(K)}\le\sqrt{2}\|D^2u_h\|_{L^2(K)}$ for $K\in\Th$. Dividing through by $\|(u_h,c_h)\|_{h,1}$, we obtain~(\ref{unifbound}).$\quad\quad\square$
\end{section}
\begin{section}{Error analysis}\label{section:5}
Herein we will denote $a\lesssim b$ for $a,b\in\mathbb R$, if there exists a constant $C>0$, such that
$$a \le Cb,$$
independent of $\mathbf{h}:=\{h_K:K\in\Th\}$, and $u$, but otherwise possibly dependent on the polynomial degree, $p$, the shape-regularity
constants of $\Th$, $C_\mathcal{T}$, $\sigma$, $c$, $\mathbf{s}$, etc.
\begin{theorem}\label{error:est:thm}
%
\begingroup\color{black}Let $\X\subset\mathbb R^2$ be a  $C^2$ and piecewise $C^{m+1}$ domain, $m\in\mathbb N$, $m\ge 2$, and let $\beta\in C^1(\partial\X;\mathbb S^1)$. Assume that
$$\partial_{{\Ta}}\Theta+\mc_{\partial\X}>0\quad\mbox{on }\partial\X.$$
Furthermore, assume that $\{\Th\}_{h}$ is a regular of order $m$ family of triangulations on $\overline{\X}$ satisfying Assumption~\ref{Meshconds}.
Let $(u,c)\in H^2_{\beta,0}(\Omega)\times\mathbb R$ be the unique strong solution of~(\ref{1}). Assume that $u\in H^\mathbf{s}(\Omega;\Th)$ with $s_K>5/2$ for all $K\in\Th$. 

Let $c_{\operatorname{stab}}$, $c_*$, $\mu_F$ and $\sigma_F$ be
chosen as in Theorem~\ref{Coercivity2}, and choose $\eta_F\lesssim 1/\tilde{h}_F^3$, $\mu_F,\sigma_F\lesssim 1/\tilde{h}_F$, $F\in\Eib$. Furthermore, for $F\in\Eb$, let $\tilde{h}_F^{1-2t^*_F}\lesssim\ell_F$, where $t^*_F:=\max_{K\in\Th:|\partial K\cap\overline{F}|\ne0}t_K$, and $t_K:=\min\{p+1,s_K,m+1\}$.
Then, there exists
a constant $C > 0$, independent of $h$, and $u$, but depending on $\max_K s_K$, such that
\begin{equation}\label{5.1}
\begin{aligned}
\|(u-u_h,c-c_h)\|_{h,1}&\le C\left(\sum_{K\in\Th}h_K^{2t_K-4}\|u\|^2_{H^{s_K}(K)}\right)^{1/2}+\left(\sum_{e\in\Vb:\jump{c_h}_e\ne 0}h_{K_e}^{2t_{K_e}-4}\|u\|_{H^{s_{K_e}}(K_e)}^2\right)^{1/2},
\end{aligned}
\end{equation}
where, for a given $e\in\Vb$ such that $\jump{c_h}_e\ne 0$, $K_e\in\Th$ has $e$ as a vertex.
Note that for the special case of quasi-uniform meshes, denoting $s:=\min_Ks_K$, the a priori error bound~(\ref{5.1}) simplifies to
\begin{equation*}
\|(u-u_h,c-c_h)\|_{h,1}\le Ch^{\min(p+1,s,m)-2}\|u\|_{H^s(\Omega)}.
\end{equation*}
\endgroup
\end{theorem}
\emph{Proof:} Let us take $z_h\in\dg$, and denote by $\psi_h:=z_h-u_h$, $\xi_h:=z_h-u$, and $\mu_h:=c-c_h$. Then, we see that
\begin{equation}\label{triineq23:57}
\|(u-u_h,c-c_h)\|_{h,1} = \|(\xi_h+\psi_h,\mu_h)\|_{h,1}
\le \|(\xi_h,0)\|_{h,1}+\|(\psi_h,\mu_h)\|_{h,1}.
\end{equation}
The proof we present relies on the existence of a $z_h\in \dg$ and a constant $C$, independent of $u$, $h_K$, but dependent on $\max_Ks_K$, such that for each $K\in\Th$, each nonnegative integer $q\le\min\{s_K,m\}$, and each multi-index $\alpha$ with $|\alpha|<s_K-1/2$, we have
\begin{equation}\label{opt:interp}
\begin{aligned}
\|u-z_h\|_{H^q(K)}&\lesssim h_K^{t_K-q}\|u\|_{H^{s_K}(K)},\\
\|D^\alpha(u-z_h)\|_{L^2(\partial K)}&\lesssim Ch^{t_K-|\alpha|-1/2}\|u\|_{H^{s_K}(K)}.
\end{aligned}
\end{equation}
The error estimates given by the first inequality in~(\ref{opt:interp}) is given in~\cite{babuvska1987hp} in the context of meshes consisting of simplices that do not have curved faces. These results, however, still hold when elements of the mesh are curved. First one must note that the second inequality in~(\ref{opt:interp}) follows from the trace inequality, followed by an application of the first inequality in~(\ref{opt:interp}). Furthermore, in~\cite{MR1014883}, the first bound in~(\ref{opt:interp}) is derived (see Corollary 4.1 in~\cite{MR1014883}) for integer values of $s_K$. However, for non-integer values of $s_K$, the estimate can be proven via scaling.

Due to our assumptions upon the parameters $\mu_F,\eta_F,\ell_F$ and $\sigma_F$, by applying the estimates in~(\ref{opt:interp}), we obtain
$$\|(\xi_h,0)\|_{h,1}\lesssim \left(\sum_{K\in\Th}h_K^{2t_K-4}\|u\|^2_{H^{s_K}(K)}\right)^{1/2},$$
thus, by~(\ref{triineq23:57}), it is sufficient to obtain the following estimate: 
\begin{equation}\label{andwe'redone01:13}
\|(\psi_h,\mu_h)\|_{h,1}\lesssim\left(\sum_{K\in\Th}h_K^{2t_K-4}\|u\|^2_{H^{s_K}(K)}\right)^{1/2}+\left(\sum_{e\in\Vb:\jump{c_h}_e\ne 0}h_{K_e}^{2t_{K_e}-4}\|u\|_{H^{s_{K_e}}(K_e)}^2\right)^{1/2}.
\end{equation}
Now, applying the coercivity result from Theorem~\ref{Coercivity2}, we obtain
\begin{equation}\label{coercpsi16:36}
\begin{aligned}
\|(\psi_h,\mu_h)\|_{h,1}^2
&\lesssim A_h((\psi_h,\mu_h),(\psi_h,\mu_h))\\
& = A_h((z_h-u_h,c-c_h),(\psi_h,\mu_h))\\
& =  A_h((z_h,c),(\psi_h,\mu_h)) - A_h((u_h,c_h),(\psi_h,\mu_h)).
\end{aligned}
\end{equation}
We then utilise the consistency identity~(\ref{obl:Cons:identity}), noting that $c$ is constant, and the fact that the pair $(u_h,c_h)\in M_h$ satisfies~(\ref{FEM}), yielding
\begin{equation*}
\begin{aligned}
A_h((u_h,c_h),(\psi_h,\mu_h))&=\sum_{K\in\Th}\langle\gamma f,\Delta\psi_h\rangle_K\\
& = A_h((u,c),(\psi_h,c))\\
& = A_h((u,c),(\psi_h,c-c_h)+(0,c_h))\\
& = A_h((u,c),(\psi_h,\mu_h))+A_h((u,c),(0,c_h)).
\end{aligned}
\end{equation*}
We apply the above identity to~(\ref{coercpsi16:36}), which results in
\begin{equation*}
\begin{aligned}
\|(\psi_h,\mu_h)\|_{h,1}^2 & \lesssim A_h((z_h,c),(\psi_h,\mu_h))-A_h((u,c),(\psi_h,\mu_h))-A_h((u,c),(0,c_h))\\
& = A_h((\xi_h,c),(\psi_h,\mu_h))-A_h((u,c),(0,c_h)).
\end{aligned}
\end{equation*}
From this, we obtain
$\|(\psi_h,\mu_h)\|_{h,1}\lesssim \sum_{i=1}^6A_i,$
where
\begin{equation*}
\begin{aligned}
A_1&:=\sum_{K\in\Th}\langle D^2\xi_h,D^2\psi_h\rangle_K,\,
A_2:=\sum_{K\in\Th}\langle(\gamma L-\Delta)\xi_h,\Delta\psi_h\rangle_K,\\
A_3&:=\sum_{K\in\Th}\frac{1}{2}\langle\Delta\xi_h,\Delta\psi_h\rangle_K,\,
A_4:=\frac{1}{2}B_{h,*}((\xi_h,0),(\psi_h,\mu_h)),\\
A_5&:=J_h((\xi_h,0),(\psi_h,\mu_h)),\,A_6:=-B_{h,1/2}((u,c),(0,c_h)).
\end{aligned}
\end{equation*}
We see that
\begin{equation}\label{boundlist:0}
\begin{aligned}
&|A_1|,|A_2|,|A_3|\lesssim\left(\sum_{K\in\Th}\|D^2\xi_h\|_{L^2(K)}^2\right)^{1/2}\|(\psi_h,\mu_h)\|_{h,1},\mbox{ and }\\
&|A_5|\le J_h((\xi_h,0),(\xi_h,0))^{1/2}\|(\psi_h,\mu_h)\|_{h,1}.
\end{aligned}
\end{equation}
Applying the first estimate in~(\ref{opt:interp}) to the estimates in~(\ref{boundlist:0}), we obtain
$$|A_1|,|A_2|,|A_3|\lesssim\left(\sum_{K\in\Th}h_K^{2t_K-4}\|u\|^2_{H^{s_K}(K)}\right)^{1/2}\|(\psi_h,\mu_h)\|_{h,1}.$$
We also see that
$|A_5|\lesssim(e_1+e_2+e_3)^{1/2}\|\psi_h\|_{h,1},$
where, based on the assumption that $\eta_F\lesssim1/\tilde{h}_F^3$, and $\sigma_F\lesssim1/\tilde{h}_F$,\begin{equation*}
\begin{aligned}
e_1 &:=\sum_{F\in\Ei}\mu_F[\|\jump{\nabla\xi_h\cdot n_F}\|_{L^2(F)}^2+\|\jump{\nabla_\Ta\xi_h}\|_{L^2(F)}^2]
\lesssim\sum_{F\in\Ei}\frac{1}{\tilde{h}_F}\|\nabla\xi_h\|_{L^2(F)}^2,\\
e_2&:=\sum_{F\in\Ei}\eta_F\|\jump{\xi_h}\|_{L^2(F)}^2\lesssim\sum_{F\in\Ei}\frac{1}{\tilde{h}_F^3}\|\xi_h\|_{L^2(F)}^2,\\
e_3&:=\sum_{F\in\Eb}\sigma_F\|\beta\!\cdot\!\nabla\xi_h\|_{L^2(F)}^2\lesssim\sum_{F\in\Eb}\frac{1}{\tilde{h}_F}\|\nabla\xi_h\|_{L^2(F)}^2,
\end{aligned}
\end{equation*}
which are all bounded above by $C\sum_{K\in\Th}h_K^{2t_K-4}\|u\|^2_{H^{s_K}(K)}$, due to~(\ref{opt:interp}).
Now we must obtain a bound for $A_4=(1/2)B_{h,*}((\xi_h,0),(\psi_h,\mu_h))$. One can see that $B_{h,*}((\xi_h,0),(\psi_h,\mu_h)) =: \sum_{i=1}^6I_i$, for which
\begin{equation*}
\begin{aligned}
|I_1|&:=\left|\sum_{K\in\Th}\langle D^2\xi_h,D^2\psi_h\rangle_K\right|\lesssim\left(\sum_{K\in\Th}\|D^2\xi_h\|_{L^2(K)}^2\right)^{1/2}\left(\sum_{K\in\Th}\|D^2\psi_h\|_{L^2(K)}^2\right)^{1/2},\\
|I_2|&:=\left|\sum_{F\in\Ei}[\langle\divTa\nablaTa\avg{\xi_h},\jump{\nabla\psi_h\cdot n_F}\rangle_F+\langle\divTa\nablaTa\avg{\psi_h},\jump{\nabla\xi_h\cdot n_F}\rangle_F]\right|\\
&\lesssim\left(\sum_{F\in\Ei}\tilde{h}_F\|D^2\xi_h\|_{L^2(F)}^2\right)^{1/2}\left(\sum_{F\in\Ei}\frac{1}{\tilde{h}_F}\|\jump{\nabla\psi_h\cdot n_F}\|_{L^2(F)}^2\right)^{1/2}\\
&~~~~~~~+\left(\sum_{F\in\Ei}\frac{1}{\tilde{h}_F}\|\nabla\xi_h\|_{L^2(F)}^2\right)^{1/2}\left(\sum_{F\in\Ei}\tilde{h}_F\|\divTa\nablaTa\avg{\psi_h}\|_{L^2(F)}^2\right)^{1/2}\\
&\lesssim\left(\left(\sum_{F\in\Ei}\tilde{h}_F\|D^2\xi_h\|_{L^2(F)}^2\right)^{1/2}+\left(\sum_{F\in\Ei}\frac{1}{\tilde{h}_F}\|\nabla\xi_h\|_{L^2(F)}^2\right)^{1/2}\right)\|(\psi_h,\mu_h)\|_{h,1},\\
|I_3|&:=\left|-\sum_{F\in\Ei}[\langle\nablaTa\avg{\nabla\xi_h\cdot n_F},\jump{\nablaTa\psi_h}\rangle_F+\langle\nablaTa\avg{\nabla\psi_h\cdot n_F},\jump{\nablaTa\xi_h}\rangle_F]\right|\\
&\lesssim\left(\left(\sum_{F\in\Ei}\tilde{h}_F\|D^2\xi_h\|_{L^2(F)}^2\right)^{1/2}+\left(\sum_{F\in\Ei}\frac{1}{\tilde{h}_F}\|\nabla\xi_h\|_{L^2(F)}^2\right)^{1/2}\right)\|(\psi_h,\mu_h)\|_{h,1},\\
\end{aligned}
\end{equation*}
\begin{equation*}
\begin{aligned}
|I_4|&:=\left|\sum_{F\in\Eb}\langle(\partial_\Ta\Theta+\mc_F)\nabla\xi_h,\nabla\psi_h\rangle_F\right|=\left|\sum_{F\in\Eb}\langle\nabla\xi_h,\langle(\partial_\Ta\Theta+\mc_F)\nabla\psi_h\rangle_F\right|\\
&\lesssim\left(\sum_{F\in\Eb}\|\nabla\xi_h\|_{L^2(F)}^2\right)^{1/2}\|(\psi_h,\mu_h)\|_{h,1},\\
|I_5|&:=\left|\sum_{F\in\Eb}\langle\partial_\Ta(\betaperp\!\cdot\!\nabla\xi_h),\beta\!\cdot\!\nabla\psi_h-\mu_h\rangle_F\right|\\
&\lesssim\left(\sum_{F\in\Eb}\tilde{h}_F\|D^2\xi_h\|_{L^2(F)}^2\right)^{1/2}\left(\sum_{F\in\Eb}\frac{1}{\tilde{h}_F}\|\beta\!\cdot\!\nabla\psi_h-\mu_h\|_{L^2(F)}^2\right)^{1/2}\\
&~~~~~~+\left(\sum_{F\in\Eb}\|\nabla\xi_h\|_{L^2(F)}^2\right)^{1/2}\left(\sum_{F\in\Eb}\|\beta\!\cdot\!\nabla\psi_h-\mu_h\|_{L^2(F)}^2\right)^{1/2}\\
&\lesssim\left(\left(\sum_{F\in\Eb}\tilde{h}_F\|D^2\xi_h\|_{L^2(F)}^2\right)^{1/2}+\left(\sum_{F\in\Eb}\|\nabla\xi_h\|_{L^2(F)}^2\right)^{1/2}\right)\|(\psi_h,\mu_h)\|_{h,1},\\
|I_6|&:=\left|\sum_{F\in\Eb}\langle\partial_\Ta(\betaperp\!\cdot\!\nabla\psi_h),\beta\!\cdot\!\nabla\xi_h\rangle_F\right|\\
&\lesssim\left(\sum_{F\in\Eb}\frac{1}{\tilde{h}_F}\|\nabla\xi_h\|_{L^2(F)}^2\right)^{1/2}\left(\sum_{F\in\Eb}\tilde{h}_F\|\partial_\Ta(\betaperp\!\cdot\!\nabla\psi_h)\|_{L^2(F)}^2\right)^{1/2}\\
&\lesssim\left(\sum_{F\in\Eb}\frac{1}{\tilde{h}_F}\|\nabla\xi_h\|_{L^2(F)}^2\right)^{1/2}\|(\psi_h,\mu_h)\|_{h,1},
\end{aligned}
\end{equation*}
note that obtaining the final inequality in the estimate for $I_6$ is analogous to~(\ref{I3star}).

Applying both estimates from~(\ref{opt:interp}) to the estimates for $I_1,\ldots,I_6$, we obtain
$$|A_4|\le\sum_{i=1}^6|I_i|\lesssim\left(\sum_{K\in\Th}h_K^{2t_K-4}\|u\|^2_{H^{s_K}(K)}\right)^{1/2}\|(\psi_h,\mu_h)\|_{h,1}.$$
It then follows that
\begin{equation}\label{01:10}
\|(\psi_h,\mu_h)\|_{h,1}^2\le\sum_{i=1}^6|A_i|\lesssim\left(\sum_{K\in\Th}h_K^{2t_K-4}\|u\|_{H^{s_K}(K)}^2\right)^\frac{1}{2}\|(\psi_h,\mu_h)\|_{h,1}+A_6.
\end{equation}
Based upon our assumptions upon $\ell_F$, and Sobolev embeddings, it follows that 
\begin{equation}\label{A_6bound}
\begin{aligned}
A_6&\lesssim\left(\sum_{e\in\Vb:\jump{c_h}_e\ne 0}(h_e\ell_F)^{-1}|\nabla u(e)|^2\right)^{1/2}\|(\psi_h,\mu_h)\|_{h,1}\\
&\lesssim\left(\sum_{e\in\Vb:\jump{c_h}_e\ne 0}h_{K_e}^{2t_{K_e}-4}\|u\|_{H^{s_{K_e}}({K_e})}^2\right)^{1/2}\|(\psi_h,\mu_h)\|_{h,1},
\end{aligned}
\end{equation}
where, for a given $e\in\Vb$ such that $\jump{c_h}_e\ne 0$,  $K_e\in\Th$ has $e$ as a vertex, and a face $F\in\Ei$ satisfies $F\subset\partial K$, and $h_e=\tilde{h}_F$ (recall the definition~(\ref{hedef}) of $h_e$).

The first inequality of~(\ref{A_6bound}) holds due to the following argument. We see that
\begin{equation*}
\begin{aligned}
A_6 & = -B_{h,1/2}((u,c),(0,c_h))\\
& = \sum_{F\in\Eb}\left[\frac{1}{2}\langle\partial_\Ta(\betaperp\cdot\nabla u),c_h\rangle_F+\sigma_F\langle\beta\cdot\nabla u-c,c_h\rangle_F\right]+\sum_{F\in\Ei}\ell_F\langle\jump{c},\jump{c_h}\rangle_F\\
& = \frac{1}{2}\sum_{F\in\Eb}\langle\partial_\Ta(\betaperp\cdot\nabla u),c_h\rangle_F,
\end{aligned}
\end{equation*}
where the final equality holds due to the fact that $\beta\cdot\nabla u-c|_F=0$ for all $F\in\Eb$, and as $c$ is constant, it cannot jump across internal edges. Upon integrating by parts on each $F\in\Eb$, we see that
\begin{equation*}
\begin{aligned}
A_6 & = \frac{1}{2}\sum_{F\in\Eb}\langle\partial_\Ta(\betaperp\cdot\nabla u),c_h\rangle_F\\
& = \frac{1}{2}\sum_{e\in\Vb}\jump{(\betaperp\cdot\nabla u)c_h}_e\\
& =  \frac{1}{2}\sum_{e\in\Vb}\jump{\betaperp\cdot\nabla u}_e\avg{c_h}_e+\avg{\betaperp\cdot\nabla u}_e\jump{c_h}_e\\
& =  \frac{1}{2}\sum_{e\in\Vb}\avg{\betaperp\cdot\nabla u}_e\jump{c_h}_e\\
& = \frac{1}{2}\sum_{e\in\Vb:\jump{c_h}_e\ne0}\avg{\betaperp\cdot\nabla u}_e\jump{c_h}_e.
\end{aligned}
\end{equation*}
The penultimate equality holds, due to the fact that $\betaperp\in C^1(\partial\X;\mathbb S^1)$, and so cannot jump across vertices, furthermore, $\nabla u\in H^{1/2}(\partial\X)$, and thus, since $\partial\X$ is a one-dimensional hypersurface, neither can $\nabla u$. 

For a given $e\in\Vb$, we have two cases to consider: either $e\cap\Ei =\emptyset$ or $e\cap\Ei\ne\emptyset$. In the first case, $e=\overline{F}\cap\overline{F'}$, for some $F,F'\in\Eb$, where $F$ and $F'$ are faces of \emph{one} $K\in\Th$. Since $c_h$ is piecewise constant, it must be constant on $K$. Thus, $c_h|_F=c_h|_{F'}$, and so, $\jump{c_h}_e=0$.
This tells us that if $e\in\Vb$ such that $\jump{c_h}_e\ne0$, then $e\cap\Ei\ne\emptyset$.
 
When $e\cap\Ei\ne\emptyset$, we see that, $e=\bigcap_{m=1}^{N_e}\overline{F_m}$, for some finite collection $F_1,\ldots,F_{N_e}\in{\Ei}$ such that $\overline{F_m}\cap\partial\X\ne\emptyset$ for each $m=1,\ldots,N_e$ (that is, $e$ may be expressed as the intersection of several faces that have a nonempty intersection with the boundary of $\X$). Furthermore, ordering the collection of faces (without relabelling them) with respect to an anticlockwise orientation, each $F_i$ in this collection satisfies 
$$F_m=\overline{K_m}\cap\overline{K_{m+1}},\quad m=1,\ldots,N_e,$$
for some $K_m,K_{m+1}\in\Th$ that have a nonempty intersection with $\partial\X$. We thus see that
\begin{equation*}
\begin{aligned}
\jump{c_h}_e & = c_h|_{K_{N_e+1}}-c_h|_{K_1}\\
& = \sum_{m=1}^{N_e}c_h|_{K_{m+1}}-c_h|_{K_{m}}\\
&\le\sum_{m=1}^{N_e}\left|c_h|_{K_{m+1}}-c_h|_{K_{m}}\right|\\
&=\sum_{m=1}^{N_e}|\jump{c_h}_{F_m}|,
\end{aligned}
\end{equation*}
which yields,
\begin{equation*}
\begin{aligned}
|\jump{c_h}_e|^2& \le\left(\sum_{m=1}^{N_e}|\jump{c_h}_{F_m}|\right)^2\\
& =\sum_{m,n=1}^{N_e}|\jump{c_h}_{F_m}||\jump{c_h}_{F_n}|\\
&\le\frac{1}{2}\sum_{m,n=1}^{N_e}|\jump{c_h}_{F_m}|^2+|\jump{c_h}_{F_n}|^2\\
& = N_e\sum_{m=1}^{N_e}|\jump{c_h}_{F_m}|^2\\
& = N_e\sum_{F\in\Ei\,:\,\overline{F}\cap\,e\ne\emptyset}|\jump{c_h}_{F}|^2.
\end{aligned}
\end{equation*}
Thus, applying the Cauchy--Schwarz inequality for vectors in $\mathbb R^n$, we obtain
\begin{equation*}
\begin{aligned}
A_6 & =  \frac{1}{2}\sum_{e\in\Vb:\jump{c_h}_e\ne 0}\avg{\betaperp\cdot\nabla u}_e\jump{c_h}_e\\
&\le\frac{1}{2}\left(\sum_{e\in\Vb:\jump{c_h}_e\ne 0}\alpha_e^{-1}|\avg{\betaperp\cdot\nabla u}_e|^2\right)^{1/2}\left(\sum_{e\in\Vb:\jump{c_h}_e\ne 0}\alpha_e|\jump{c_h}_e|^2\right)^{1/2}\\
&\le\frac{1}{2}\left(\sum_{e\in\Vb:\jump{c_h}_e\ne 0}\alpha_e^{-1}|\avg{\betaperp\cdot\nabla u}_e|^2\right)^{1/2}\left(\sum_{e\in\Vb:\jump{c_h}_e\ne 0}\alpha_eN_e\left(\sum_{F\in\Ei\,:\,\overline{F}\cap\,e\ne\emptyset}|\jump{c_h}_F|^2\right)\right)^{1/2}.
\end{aligned}
\end{equation*}
for any collection of positive real numbers $\{\alpha_e\}_{e\in\Vb:\jump{c_h}_e\ne 0}$. Let us denote
$$r_h:=\frac{1}{2}\left(\sum_{e\in\Vb:\jump{c_h}_e\ne 0}\alpha_e^{-1}|\avg{\betaperp\cdot\nabla u}_e|^2\right)^{1/2}.$$ Note that the value $\max_{e\in\Vb:\jump{c_h}_e\ne 0}N_e$ is bounded independently of the mesh size, due to the mesh regularity assumptions, and so
$$A_6\lesssim r_h\left(\sum_{e\in\Vb:\jump{c_h}_e\ne 0}\alpha_e\left(\sum_{F\in\Ei\,:\,\overline{F}\cap\,e\ne\emptyset}|\jump{c_h}_F|^2\right)\right)^{1/2}.$$ We now choose $\alpha_e$ such that
$$\alpha_e = \min_{F\in\Ei\,:\,\overline{F}\cap\,e\ne\emptyset}\ell_F|F|.$$
We also see that if $F\in\Ei$, then $F$ can intersect at most one vertex belonging to $\{e\in\Vb:\jump{c_h}_e\ne 0\}$; thus
$$\bigcup_{e\in\Vb:\jump{c_h}_e\ne 0}\{F\in\Ei\,:\,\overline{F}\cap\,e\ne\emptyset\}\subset\{F\in\Ei:F\cap\partial\X\ne\emptyset\}.$$
It then follows that
\begin{equation}\label{a6bnd}
A_6\lesssim r_h\left(\sum_{F\in\Ei:F\cap\partial\X\ne\emptyset}\ell_F|F||\jump{c_h}_F|^2\right)^{1/2}.
\end{equation}
Since $\jump{c_h}_F$ is constant on $F$, one also has that 
$$|F||\jump{c_h}_F|^2=|\jump{c_h}_F|^2\int_F1=\int_F|\jump{c_h}_F|^2= \|\jump{c_h}_F\|_{L^2(F)}^2,$$ hence,
\begin{equation*}
\begin{aligned}
\sum_{F\in\Ei:F\cap\partial\X\ne\emptyset}\ell_F|F||\jump{c_h}_F|^2& = \left(\sum_{F\in\Ei:\overline{F}\cap\partial\X\ne\emptyset}\ell_F\|\jump{c_h}\|_{L^2(F)}^2\right)\\
&\le\left(\sum_{F\in\Ei}\ell_F\|\jump{c_h}\|_{L^2(F)}^2\right)\\
& = \left(\sum_{F\in\Ei}\ell_F\|\jump{c-c_h}\|_{L^2(F)}^2\right)\\
&\lesssim\|(\psi_h,\mu_h)\|_{h,1}^2.
\end{aligned}
\end{equation*}
Applying this to~(\ref{a6bnd}), we obtain
\begin{equation}\label{5.6local}
A_6\lesssim r_h\|(\psi_h,\mu_h)\|_{h,1}.
\end{equation}
Overall, we have obtained
\begin{equation*}
\|(\psi_h,\mu_h)\|_{h,1}^2\lesssim\left(r_h+\left(\sum_{K\in\Th}h_K^{2t_K-4}\|u\|^2_{H^{s_K}(K)}\right)^{1/2}\right)\|(\psi_h,\mu_h)\|_{h,1}.
\end{equation*}
Thanks to the mesh condition~(\ref{meshcond2}), we have that $$\alpha_e^{-1}\lesssim |F|^{-1}\ell_F^{-1}$$
for any $F\in\Ei$ such that $\overline{F}\cap\,e\ne\emptyset$. Furthermore, $|F|^{-1}\lesssim\tilde{h}_e^{-1}$. Thus,
\begin{equation}\label{rhbnd}
\begin{aligned}
r_h&=\frac{1}{2}\left(\sum_{e\in\Vb:\jump{c_h}_e\ne 0}\alpha_e^{-1}|\avg{\betaperp\cdot\nabla u}_e|^2\right)^{1/2}\\
&\lesssim\left(\sum_{e\in\Vb:\jump{c_h}_e\ne 0}\alpha_e^{-1}|\nabla u(e)|^2\right)^{1/2}\\
&\lesssim\left(\sum_{e\in\Vb:\jump{c_h}_e\ne 0}(h_e\ell_F)^{-1}|\nabla u(e)|^2\right)^{1/2},
\end{aligned}
\end{equation}
applying~(\ref{rhbnd}) to~(\ref{5.6local}), we obtain
$$A_6\lesssim\left(\sum_{e\in\Vb:\jump{c_h}_e\ne 0}(h_e\ell_F)^{-1}|\nabla u(e)|^2\right)^{1/2}\|(\psi_h,\mu_h)\|_{h,1},$$
which is the first inequality of~(\ref{A_6bound}). We obtain the final estimate of~(\ref{A_6bound}), by first noting that for $e\in\Vb$ such that $\jump{c_h}_e\ne 0$, $e$ is a vertex of some $K_e\in\Th$, and $u\in H^{s_{K_e}}(K_e)$, with $s_{K_e}>5/2$, and so, by Sobolev embeddings and scaling, we obtain
$$|\nabla u(e)|^2\lesssim h_{K_e}^{-1}\|u\|_{H^{s_{K_e}}(K_e)}^2.$$
Furthermore, we may choose $K_e$ such that there is an $F\in\Ei$ that satisfies $F\subset\partial K_e$, and $h_e^{-1}=\tilde{h}_F^{-1}\lesssim h_{K_e}^{-1}$, and thus by our assumptions on $\ell_F$, we have $\ell_F^{-1}\lesssim\tilde{h}_F^{2t^*_F-1}\lesssim h_{K_e}^{2t_K-1}$. Thus, we obtain
\begin{align*}
\sum_{e\in\Vb:\jump{c_h}_e\ne 0}(h_e\ell_F)^{-1}|\nabla u(e)|^2&\lesssim\sum_{e\in\Vb:\jump{c_h}_e\ne 0}h_e^{-1}\ell_F^{-1}h_{K_e}^{-2}\|u\|_{H^{s_{K_e}}(K_e)}^2\\
&\lesssim\sum_{e\in\Vb:\jump{c_h}_e\ne 0}h_{K_e}^{2t_{K_e}-4}\|u\|_{H^{s_{K_e}}(K_e)}^2.
\end{align*}

Combining our estimates for $|A_1|,\ldots,|A_5|$, and $A_6$, yields
\begin{equation*}
\begin{aligned}
\|(\psi_h,\mu_h)\|_{h,1}^2\lesssim\sum_{i=1}^6A_i
&\lesssim\left(\left(\sum_{K\in\Th}h_K^{2t_K-4}\|u\|^2_{H^{s_K}(K)}\right)^{1/2}\right.\\
&\left.~~+\left(\sum_{e\in\Vb:\jump{c_h}_e\ne 0}h_{K_e}^{2t_{K_e}-4}\|u\|_{H^{s_{K_e}}(K_e)}^2\right)^{1/2}\right)\|(\psi_h,\mu_h)\|_{h,1}.
\end{aligned}
\end{equation*}
Finally, upon noting that
$\|(\xi_h,0)\|_{h,1}\lesssim\left(\sum_{K\in\Th}h_K^{2t_K-4}\|u\|^2_{H^{s_K}(K)}\right)^{1/2},$
we obtain
\begin{equation*}
\begin{aligned}
&\|(u-u_h,c-c_h)\|_{h,1}\le\|(\xi_h,0)\|_{h,1}+\|(\psi_h,\mu_h)\|_{h,1}\\
&\lesssim\left(\sum_{K\in\Th}h_K^{2t_K-4}\|u\|^2_{H^{s_K}(K)}\right)^{1/2}+\sum_{e\in\Vb:\jump{c_h}_e\ne 0}h_{K_e}^{2t_{K_e}-4}\|u\|_{H^{s_{K_e}}(K_e)}^2,
\end{aligned}
\end{equation*}
as desired. $\quad\quad\square$
\subsection{An error estimate in the case of conformal regularity}
The hypotheses of Theorem~\ref{error:est:thm} includes the sufficient condition that the strong solution, $u$, is piecewise-sufficiently regular, so that one may substitute $(u,c)$ into the left hand argument of the operator, $A_h$. In the following lemma, we provide an error estimate for strong solutions $u\in H^2_{\beta,0}$, i.e., the expected conformal regularity of strong solutions implied by Theorem~\ref{linear:existence+uniqueness}. As in estimate~(\ref{5.1}), one can see the error contribution arising from the inconsistency of $c_h$ belonging to $\dgoo$ as opposed to $\mathbb R$. Similarly, this contribution is zero if $c_h$ does not jump across boundary vertices. This shows that our method provides an approximation that is at least as accurate in the $\|(\cdot,\cdot)\|_{h,1}$-norm, as a $H^2$-conforming finite element method. 
\begin{lemma}\label{min:reg:lemma}
Let $\X\subset\mathbb R^2$ be a $C^2$ and piecewise $C^3$ domain, and let $\beta\in C^1(\partial\X;\mathbb S^1)$. Assume that
$$\partial_{{\Ta}}\Theta+\mc_{\partial\X}>0\quad\mbox{on }\partial\X.$$
Furthermore, assume that $\{\Th\}_{h}$ is a regular of order $2$ family of triangulations on $\overline{\X}$ satisfying Assumption~\ref{Meshconds}.
Let $(u,c)\in H^2_{\beta,0}(\Omega)\times\mathbb R$ be the unique strong solution of~(\ref{1}). Let $c_{\operatorname{stab}}$, $c_*$, and $\mu_F$ be
chosen as in Theorem~\ref{Coercivity2}, and choose $\eta_F\lesssim 1/\tilde{h}_F^3$, $\sigma_F\lesssim 1/\tilde{h}_F$, and $\tilde{h}_F^{1+r+p^*}\lesssim\ell_F$, where $r>0$, and $p^*:=2\operatorname{sgn}(p-2)$. Then, we have the following error estimate
\begin{equation}\label{min:error:est}
\begin{aligned}
&\|(u-u_h,c-c_h)\|_{h,1}\lesssim\inf_{z_h\in V}\left\{\|u-z_h\|_{H^2(\X)}+\left[\sum_{F\in\Eb}\frac{1}{\tilde{h}_F}\|\beta\!\cdot\!\nabla(u-z_h)\|_{L^2(F)}^2\right]^\frac{1}{2}\right.\\
&~~~~~~~~\left.+\left[\sum_{F\in\Eb}\frac{1}{\tilde{h}_F}\|\partial_\Ta(\beta\!\cdot\!\nabla(u-z_h))\|_{L^2(F)}^2\right]^\frac{1}{2}\right\}+\left(\sum_{e\in\Vb:\jump{c_h}_e\ne 0}\frac{\ell_F^{-1}}{h_{K_e}^{1+p^*}}\|u\|_{H^{2}(K_e)}^2\right)^\frac{1}{2}.
\end{aligned}
\end{equation}
where $V:=\dgo\cap H^2(\X)$.
\end{lemma}
\emph{Proof:} First we assume that $z_h\in H^2(\X)\cap\dgo$. Then, we see that
$$\|(u-u_h,c-c_h)\|_{h,1}\le\|(\xi_h,0)\|_{h,1}+\|(\psi_h,\mu_h)\|,$$
where $\xi_h,\psi_h$ and $\mu_h$ are given as in the proof of Theorem~\ref{error:est:thm}. Since we only assume that $z_h$ is in $H^2(\X)\cap\dgo$, only the consistency properties of the bilinear form $A_h$ that depend on the piecewise regularity and $H^2$-regularity of $z_h$ hold. In particular~(\ref{obl:deltab:2609}) does not hold. This results in:
\begin{equation*}
\begin{aligned}
\|(\psi_h,\mu_h)\|_{h,1}^2&\lesssim A_h((z_h,c),(\psi_h,\mu_h))-A_h((u_h,c_h),(\psi_h,\mu_h))\\
& = \sum_{K\in\Th}\langle\gamma(Lz_h-f),\Delta\psi_h\rangle_K\\
&~~~~+\frac{1}{2}\sum_{F\in\Eb}\langle\partial_\Ta(\betaperp\cdot\nabla z_h),c_h\rangle_F-\langle\partial_\Ta(\betaperp\cdot\nabla\psi_h),\beta\cdot\nabla z_h-c\rangle_F\\
&~~~~+\sum_{F\in\Eb}\frac{1}{2}\langle\partial_\Ta(\beta\cdot\nabla z_h),\betaperp\cdot\nabla\psi_h\rangle_F+\sigma_F\langle\beta\cdot\nabla z_h-c,\beta\cdot\nabla\psi_h-\mu_h\rangle_F\\
& = \sum_{K\in\Th}\langle\gamma L(z_h-u),\Delta\psi_h\rangle_K+B_{h,1/2}((z_h,c),(0,c_h))\\
&~~~~+\frac{1}{2}\sum_{F\in\Eb}[\langle\partial_\Ta(\beta\cdot\nabla z_h-c),\betaperp\cdot\nabla\psi_h\rangle_F-\langle\partial_\Ta(\betaperp\cdot\nabla\psi_h),\beta\cdot\nabla z_h-c\rangle_F]\\
&~~~~~~~~+\sum_{F\in\Eb}\sigma_F\langle\beta\cdot\nabla z_h-c,\beta\cdot\nabla\psi_h-\mu_h\rangle_F\\
&\lesssim\left(|z_h-u|_{H^2(\X)}+\left(\sum_{e\in\Vb:\jump{c_h}_e\ne 0}\frac{\ell_F^{-1}}{h_{K_e}}\|z_h\|_{H^{3}(K_e)}^2\right)^{1/2}\right)\|(\psi_h,\mu_h)\|_{h,1}\\
&~~~~+\frac{1}{2}\sum_{F\in\Eb}[\langle\partial_\Ta(\beta\cdot\nabla z_h-c),\betaperp\cdot\nabla\psi_h\rangle_F-\langle\partial_\Ta(\betaperp\cdot\nabla\psi_h),\beta\cdot\nabla z_h-c\rangle_F]\\
&~~~~~~~~+\sum_{F\in\Eb}\sigma_F\langle\beta\cdot\nabla z_h-c,\beta\cdot\nabla\psi_h-\mu_h\rangle_F,
\end{aligned}
\end{equation*}
where the bound for the term $B_{h,1/2}((z_h,c),(0,c_h))$ is obtained analogously to the first estimate of~(\ref{A_6bound}) in the proof of Theorem~\ref{error:est:thm} (except we utilise the fact that $z_h$ is piecewise $H^3$-regular), and we recall that $K_e\in\Th$ has $e$ as a vertex, and $F\in\Ei$ is an edge of $K_e$ that satisfies $h_e=\tilde{h}_F$. Furthermore, we have the following bounds:
\begin{equation*}
\begin{aligned}
&\sum_{F\in\Eb}\sigma_F\langle\beta\cdot\nabla z_h-c,\beta\cdot\nabla\psi_h-\mu_h\rangle_F\lesssim\left(\sum_{F\in\Eb}\sigma_F\|\beta\cdot\nabla z_h-c\|_{L^2(F)}^2\right)^{1/2}\|(\psi_h,\mu_h)\|_{h,1},\\
&\sum_{F\in\Eb}\langle\partial_\Ta(\betaperp\cdot\nabla\psi_h),\beta\cdot\nabla z_h-c\rangle_F\lesssim\\
&~~~~~~~~~~~~~~~~~~~~~\left(\sum_{F\in\Eb}\sigma_F\|\beta\cdot\nabla z_h-c\|_{L^2(F)}^2\right)^{1/2}\left(\sum_{F\in\Eb}\|\nabla\psi_h\|_{L^2(F)}^2+\tilde{h}_F\|D^2\psi_h\|_{L^2(F)}^2\right)^{1/2}\\
&~~~~~~~~~~~~~~~~~~~~~~~~~~~~~~~~~~~~~~~~~~~~~\lesssim\left(\sum_{F\in\Eb}\sigma_F\|\beta\cdot\nabla z_h-c\|_{L^2(F)}^2\right)^{1/2}\|(\psi_h,\mu_h)\|_{h,1}.
\end{aligned}
\end{equation*}
Note that the first bound follows directly from an application of the Cauchy--Schwarz inequality, and the definition of the $\|(\cdot,\cdot)\|_{h,1}$-norm, and the second inequality is obtained analogously to~(\ref{I3star}) in the proof of Lemma~\ref{Coercivity1}. 
We also see that
\begin{equation}\label{trickyterm}
\sum_{F\in\Eb}\langle\partial_\Ta(\beta\cdot\nabla z_h-c),\betaperp\cdot\nabla\psi_h\rangle_F\lesssim\left(\sum_{F\in\Eb}\frac{1}{\tilde{h}_F}\|\partial_\Ta(\beta\cdot\nabla z_h-c)\|_{L^2(F)}^2\right)^{1/2}\|(\psi_h,\mu_h)\|_{h,1}.
\end{equation}
Since $\sigma_F\lesssim1/\tilde{h}_F$, overall, we have obtained
\begin{equation}\label{star:tues}
\begin{aligned}
&\|(\psi_h,\mu_h)\|_{h,1}^2\lesssim\left(|z_h-u|_{H^2(\X)}+\left(\sum_{F\in\Eb}\frac{1}{\tilde{h}_F}\|\beta\cdot\nabla z_h-c\|_{L^2(F)}^2\right)^{1/2}\right.\\
&\left.+\left(\sum_{F\in\Eb}\frac{1}{\tilde{h}_F}\|\partial_\Ta(\beta\cdot\nabla z_h-c)\|_{L^2(F)}^2\right)^{1/2}+\left(\sum_{e\in\Vb:\jump{c_h}_e\ne 0}\frac{\ell_F^{-1}}{h_{K_e}}\|z_h\|_{H^{3}(K_e)}^2\right)^{1/2}\right)\|(\psi_h,\mu_h)\|_{h,1}.
\end{aligned}
\end{equation}
Finally, by Lemma 2.3 of~\cite{MR1014883}, if $p=2$, we have that $\|z_h\|_{H^3(K)}\lesssim\|z_h\|_{H^2(K)}$ for all $K\in\Th$, and if $p\ge3$, $\|z_h\|_{H^3(K)}\lesssim h_K^{-1}\|z_h\|_{H^2(K)}$. I.e., for any $p\ge2$, $\|z_h\|_{H^3(K)}^2\lesssim h_K^{-p^*}\|z_h\|_{H^2(K)}^2$, where $p^*:=2\operatorname{sgn}(p-2)$. Furthermore, by our assumptions on $\ell_F$, it follows that $\ell_F^{-1}/h_K^{p^*+1}\lesssim1$, thus
\begin{equation}\label{pmuh2}
\begin{aligned}
\frac{\ell_F^{-1}}{h_K}\|z_h\|_{H^3(K)}^2\lesssim\frac{\ell_F^{-1}}{h_K^{1+p^*}}\|z_h\|_{H^2(K)}^2&\lesssim\frac{\ell_F^{-1}}{h_K^{1+p^*}}(\|u-z_h\|_{H^2(K)}^2+\|u\|_{H^2(K)}^2)\\
&\lesssim\|u-z_h\|_{H^2(K)}^2+\frac{\ell_F^{-1}}{h_K^{1+p^*}}\|u\|_{H^2(K)}^2.
\end{aligned}
\end{equation}
Applying~(\ref{pmuh2}) to~(\ref{star:tues}), dividing through by $\|(\psi_h,\mu_h)\|_{h,1}$, and recalling that $c=\beta\cdot\nabla u|_{\partial\X}$, we obtain
\begin{align*}
\|(\psi_h,\mu_h)\|_{h,1}&\lesssim\|z_h-u\|_{H^2(\X)}+\left(\sum_{F\in\Eb}\frac{1}{\tilde{h}_F}\|\beta\cdot\nabla(z_h-u)\|_{L^2(F)}^2\right)^{1/2}\!\\
&+\left(\sum_{F\in\Eb}\frac{1}{\tilde{h}_F}\|\partial_\Ta(\beta\cdot\nabla (z_h-u))\|_{L^2(F)}^2\right)^{1/2}+\left(\sum_{e\in\Vb:\jump{c_h}_e\ne 0}\frac{\ell_F^{-1}}{h_{K_e}^{1+p*}}\|u\|_{H^{2}(K_e)}^2\right)^{1/2}.
\end{align*}
We then see that
\begin{equation*}
\begin{aligned}
\|(\xi_h,0)\|_{h,1}^2\lesssim|u-z_h|_{H^2(\X)}^2+\sum_{F\in\Eb}\left[\frac{1}{\tilde{h}_F}\|\beta\cdot\nabla(z_h-u)\|_{L^2(F)}^2+\|(\partial_\Ta\Theta+\mc_F)^{1/2}\nabla(u-z_h)\|_{L^2(F)}^2\right].
\end{aligned}
\end{equation*}
Furthermore, the trace operator is continuous from $H^1(\X)\to L^2(\partial\X)$, and so
$$\sum_{F\in\Eb}\|(\partial_\Ta\Theta+\mc_F)^{1/2}\nabla(u-z_h)\|_{L^2(F)}^2\lesssim\|\nabla(u-z_h)\|_{L^2(\partial\X)}^2\lesssim\|u-z_h\|_{H^2(\X)}^2.$$
Thus, we obtain
\begin{equation*}
\begin{aligned}
&\|(u-u_h,c-c_h)\|_{h,1}\lesssim\|z_h-u\|_{H^2(\X)}+\left(\sum_{F\in\Eb}\frac{1}{\tilde{h}_F}\|\beta\cdot\nabla(z_h-u)\|_{L^2(F)}^2\right)^\frac{1}{2}\!\\
&~~~~~~~+\left(\sum_{F\in\Eb}\frac{1}{\tilde{h}_F}\|\partial_\Ta(\beta\cdot\nabla (z_h-u))\|_{L^2(F)}^2\right)^{1/2}+\left(\sum_{e\in\Vb:\jump{c_h}_e\ne 0}\frac{\ell_F^{-1}}{h_{K_e}^{1+p^*}}\|u\|_{H^{2}(K_e)}^2\right)^{1/2}.
\end{aligned}
\end{equation*}
Note that our choice of $z_h\in H^2(\X)\cap\dgo$ was arbitrary, thus we may take an infimum over $V$ above, yielding estimate~(\ref{min:error:est}).$\quad\quad\square$
\begin{remark}[Stabilisation parameter choice]
Note that our assumption $\ell_F^{-1}\lesssim\tilde{h}_F^{1+p^*+r}$ implies that 
$\ell_F^{-1}/h_K^{1+p^*}\lesssim h_K^r$, controlling the contribution of the final term on the right-hand side of~(\ref{min:error:est}).
\end{remark}
\begin{remark}[Conforming finite element methods]
Notice that the error bound~(\ref{min:error:est}) incorporates the error arising from the approximation of the oblique derivative (i.e., the second and third term in the infimum). Firstly, one must note that if $\beta\cdot\nabla z_h|_{F}\mbox{ is constant}$, for each $F\in\Eb$, then the third term in the infimum vanishes. This occurs if $z_h\in V_c=\dgo\cap H^2_{\beta,0}(\X)$. However, it is not immediately clear that the space $V_c$ is non-empty for arbitrary $\beta\in C^1(\partial\X;\mathbb S^1)$ satisfying the hypotheses of Lemma~\ref{FullMT}. 

Furthermore, if $\beta\cdot\nabla z_h|_{e}=c_z$ for some constant $c_z$, and all $e\in\Vb$, then one may subtract an arbitrary constant on each face in~(\ref{trickyterm}) (i.e., $\betaperp\cdot\nabla\psi_h$ may be replaced with $\betaperp\cdot\nabla\psi_h-c_F$, for arbitrary constants $c_F$), and then applying the Cauchy--Schwarz inequality with a parameter, and Poincar\'e's inequality, improves the order of the estimate from $\mathcal{O}(\tilde{h}_F^{-1})$ to $\mathcal{O}(\tilde{h}_F)$.

These considerations also imply that the scheme introduced in this paper is at least as accurate as any conforming method seeking a numerical solution $u_h\in V$, where $\dgo\cap H^2_{\beta,0}(\X)\subset V\subset \dgo\cap H^2(\X)$.
\end{remark}
\end{section}
\begin{section}{Numerical results}\label{Experiments:sec:6}
\subsection{Implementation}
\emph{Software and code:} The experiments in this section have been implemented in the most recent version of the Firedrake software~\cite{Rathgeber2016,Luporini2015} (as of 3rd July 2018), which interfaces directly with PETSc~\cite{petsc-user-ref,petsc-efficient} running through a Python interface~\cite{Dalcin2011,Chaco95}. A working Firedrake script, Curved-oblique-DGFEM.py, used to generate the experiments of this Chapter is available in the Github repository:\\ https://github.com/ekawecki/FiredrakeNDV.

\emph{Linear systems and condition numbers:}
The bilinear form $\Aobl$ defined by~(\ref{A:def:15:59}) can also be considered to be similar to those present in finite element methods for fourth-order elliptic boundary-value problems (see~\cite{suli2007hp,MR2142191} for example), in the sense that the evaluation of $\Aobl((u_h,\lambda_h);(v_h,\mu_h))$ for $(u_h,\lambda_h),(v_h,\mu_h)\in M_h$ involves the integration of products of second order partial derivatives. This typically leads to the matrix $\mathbf{A}_h$, describing the linear system given by~(\ref{FEM}), to have a Euclidean norm condition number of order $h^{-4}$. This can pose difficulties when applying iterative methods to solve the linear system, and thus to ensure that we solve the linear system with sufficiently high accuracy as the mesh size $h$ decreases, we apply the Iterative refinement algorithm, i.e., Algorithm 1.1 of~\cite{carson2017new}. We implement the Iterative refinement algorithm by using the choices depicted in the code snippet below, in the Firedrake ``solve" function.
\begin{lstlisting}[language=Python]
        # implementing nullspace, as solution should have zero sum
        V_basis = VectorSpaceBasis(constant=True)
        nullspace = MixedVectorSpaceBasis(S, [V_basis, S[1]])
        
        # begin timing of linear system solve
        t = time()
        
        # solving linear system
        solve(A_gamma == L,Uh,nullspace = nullspace,
              solver_parameters = {"mat_type": "aij",
              "snes_type": "newtonls",
              "ksp_type": "preonly",
              "pc_type": "lu",
              "snes_monitor": False,
              "snes_rtol": 1e-16,
              "snes_atol": 1e-25})
        # end timing of linear system solve
        tt.append(time()-t)
\end{lstlisting}
One can also see that when executing the script in Firedrake, we record the runtimes by way of the sixth and last line above, so that we only record the time that it takes to solve the linear system. 

Furthermore, in the solver choices we include ``nullspace = nullspace", where ``nullspace" is defined on line 3, and ``mat\_type = aij". The first choice imposes that the numerical solution $u_h$ (from the pair $(u_h,c_h)\in M_h$ that satisifies~(\ref{FEM})) has a zero-sum, and the latter informs the solver that the solution consists of two parts, i.e., $u_h$ and $c_h$, and that the system may be treated in block formation.

\emph{Two-dimensional curved boundary approximation:} When implementing curved finite elements, we use a piecewise quadratic polynomial mapping to obtain a higher order approximation of the domain boundary. This is implemented in Firedrake by first using Gmsh~\cite{geuzaine2008gmsh} (version 3.0.1) to generate an affine triangulation $\X_h$ that approximates $\X$. We then define the \emph{continuous} Lagrange finite element space $\mathbb V:=\{v\in C(\overline{\X_h};\mathbb R^2):v\in\mathbb P^2(K;\mathbb R^2)\,\forall K\in\X_h\}$.
Then, we take $\psi_i:\omega_i\to\mathbb R^2$, $\omega_i\subset\mathbb R$, $i=1,\ldots,n$, to be the collection of charts that locally describe $\partial\X$, and denote $\{x_j\}_{j=1}^N$ to be the degrees of freedom of $\mathbb V$. We partition the collection of degrees of freedom by defining $J_{\operatorname{ext}} = \{j\in\{1,\ldots,N\}:x_j\in\partial\X_h\}$, and $J_{\operatorname{int}} = \{1,\ldots,N\}\setminus J_{\operatorname{ext}}$, and so $\{x_j\}_{j=1}^N = \{x_j\}_{j\in J_{\operatorname{int}}}\cup\{x_j\}_{j\in J_{\operatorname{ext}}}$. We then
 define the the function $T\in\mathbb V$ by
\begin{equation}\label{T:def}
\left\{
\begin{aligned}
T(x_j) &= x_j,\quad j\in J_{\operatorname{int}},\\
T(x_j) &=\psi_i(x_j),\quad j\in J_{\operatorname{ext}},\quad i\in\{1,\ldots,n\}\mbox{ such that }x_j\in\omega_i.
\end{aligned}
\right.
\end{equation}
Finally, we define our computational finite element space $\dgcomp:= \{v\in L^2(\X):v\circ T^{-1}\in\mathbb P^p(\hat K)\}$. 
This procedure is implemented in Firedrake, in the code snippet below, utilising the Firedrake ``Mesh" function. In this case $\X$ is the unit disk, and so there is only one chart, $\psi:=x/|x|$. Furthermore, when we refine the mesh in our experiments, the meshes at each refinement level are not related to one another (the one exception being Experiment~\ref{exp2}). That is, there is no hierarchical mesh structure, i.e., at each refinement level, we ``remesh". A collection of the meshes used for the computations of this thesis can be found in the folder ``Meshes" in the Github repository: https://github.com/ekawecki/FiredrakeNDV.
\begin{lstlisting}[language=Python]
        #Affine mesh of the unit disk, generated in Gmsh
        mesh = Mesh("quasiunifrefdisk.msh")
        # Implementing quadratic domain approximation
        V = FunctionSpace(mesh, "CG", 2)
        # Defining a function that identifies the curved portion of the boundary
        bdry_indicator = Function(V)
        bc = DirichletBC(V, Constant(1.0), 1)
        bc.apply(bdry_indicator)
        # Defining the continuous, piecewise quadratic vector-valued finite element space
        VV = VectorFunctionSpace(mesh, "CG", 2)
        T = Function(VV)
        T.interpolate(SpatialCoordinate(mesh))
        # Defining the function T given by (6.1)
        T.interpolate(conditional(abs(1-bdry_indicator) < 1e-5, T/sqrt(inner(T,T)), T))
        # Defining the curved mesh
        mesh = Mesh(T)
        # Defining the space V_{h,p}^{comp}
        FES = FunctionSpace(mesh,"DG",deg)
\end{lstlisting}
\begin{remark}[Computational parameters]\label{parameters:sec:6}
In the following experiments, we employ the following parameter choices: $c_{\operatorname{stab}} = 2.5$, $\mu_F\!=\! 2c_{\operatorname{stab}}(p-1)^2/\tilde{h}_F$, $\eta_F\!=\!15(p-1)^4/16\tilde{h}_F^3$, $\sigma_F = 2c_{\operatorname{stab}}p^2/\tilde{h}_F^2$, and $\ell_F = c_{\operatorname{stab}}\tilde{h}_F^{-3}$. The order of the computational parameters with respect to $\tilde{h}_F$ were guided by the hypotheses of Theorem~\ref{error:est:thm}. The orders with respect to $p$ for $\eta_F$ and $\mu_F$ were guided by the experiments of~\cite{MR3077903}. Finally, the value of $c_{\operatorname{stab}}$, and the orders with respect to $p$ of $\sigma_F$ and $\ell_F$ (in the case of $\ell_F$, the parameter is in fact independent of $p$) were obtained experimentally. Furthermore, we have that $\mc_F=1$ on each $F\in\Eb$, due to the fact that the experiments are on the unit disc $\X:=\{x\in\mathbb R^2:|x|<1\}$.
\end{remark}
\newpage
\subsection{Experiments}
In this section, we test the robustness of the scheme~(\ref{FEM}), with the computational domain $\X$ taken to be the unit disk, approximated in the same manner as in present in~\cite{Kawecki:77:article:On-curved}, Section 3.4. We consider various elliptic operators, $L$, that satisfy the Cordes condition~(\ref{2}). In each case, we see that the convergence rates are of the expected order in the various broken Sobolev norms considered, and in particular in the $\|\!\cdot\!\|_{h,1}$--norm, for which we have proven the error bound~(\ref{5.1}).
\subsubsection{Experiment 1}\label{exp1}
In this experiment, we consider the following problem
\begin{equation}\label{exp:1:obl}
\left\{
\begin{aligned}
\Delta u &= f,\quad\mbox{in}\quad\X,\\
\beta\cdot\nabla u & \mbox{ is constant on }\partial\X,
\end{aligned}
\right.
\end{equation}
where $\X=\{x\in\mathbb R^2:|x|<1\}$, and $\beta\equiv n_{\partial\X}$. In this case $f$ is chosen so that the solution of~(\ref{exp:1:obl}) is given by $u(x) = \frac{1}{6}|x|^6-\frac{1}{2}|x|^2+\frac{5}{24}.$ Notice that in this case, the compatibility constant $c=0$, $\gamma=1$, and $\Theta\equiv\partial_\Ta\Theta\equiv0$

In this experiment, we successively increase the degree, $p$, of the finite element space $\dgcompo$ from $2$ to $4$, and for each fixed degree we refine the mesh quasi--uniformly, we observe that the experimental orders of convergence in the $\|\cdot\|_{h,1}$-norm are optimal, that is $\|(e_h^u,e_h^c)\|_{h,1}=\mathcal{O}(h^{p-1})$, where $(e_h^u,e_h^c):=(u-u_h,c-c_h)$. We also observe that $\|e_h^c\|_{L^2(\partial\X)}=\mathcal{O}(h^{p})$.
We plot the error values in the $\|\cdot\|_{h,1}$-norm, and plot the error arising in the approximation of the compatibility constant in Figure~\ref{plot:s6:exp1}, and report the exact values in Tables~\ref{exp1errs:s6} and~\ref{exp1Cerrs:s6}, with the corresponding experimental orders of convergence given in brackets.
\begin{figure}[h]
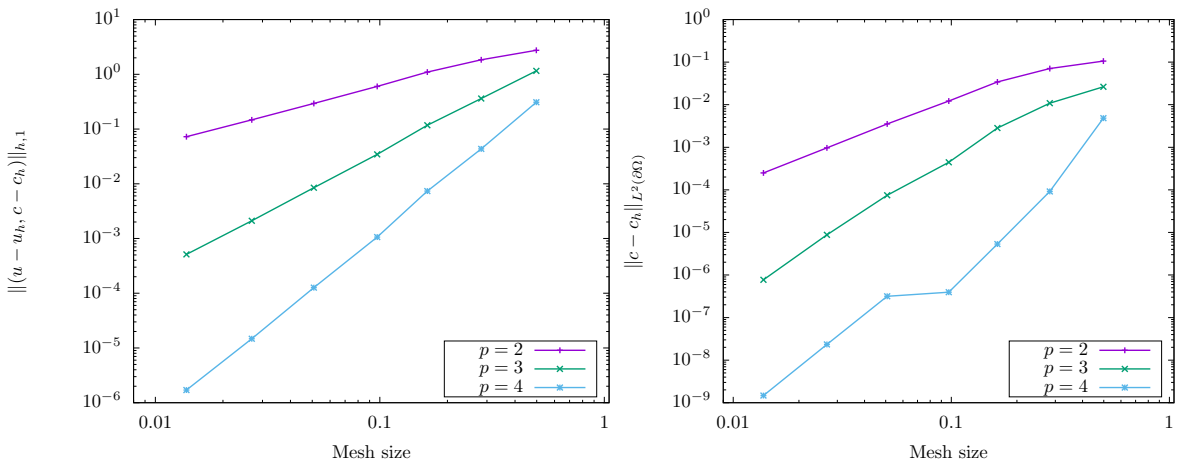

\begin{center}
  \begin{tabular}{lrr}
  \hspace{-1.3cm}
    \scalebox{0.7}{%
    \renewenvironment{table}[1][]{\ignorespaces}{\unskip}%
    \input{sec_6_exp_1.tex}
    \unskip
    }
    &
    \hspace{-2cm}
    \scalebox{0.7}{%
    \renewenvironment{table}[1][]{\ignorespaces}{\unskip}%
    \input{sec_6_exp_1_C.tex}
    \unskip
    }
  \end{tabular}
 \end{center}
 \caption{Convergence rates for the numerical scheme applied to problem~(\ref{exp:1:obl}). We provide the error values $\|(u-u_h,c-c_h)\|_{h,1}$ (left), and $\|c-c_h\|_{L^2(\partial\X)}$ (right).
We observe that the convergence rates in the $\|\cdot\|_{h,1}$ norm are optimal with respect to the choice of polynomial degree, $p$. That is, $\|(u-u_h,c-c_h)\|_{h,1} = \mathcal{O}(h^{p-1})$. Furthermore, we observe that $\|c-c_h\|_{L^2(\partial\X)}=\mathcal{O}(h^p)$.
}\label{plot:s6:exp1}
\end{figure}
\begin{table}[H]
\begin{center}
\begin{tabular}{|c|cc|cc|cc|}\hline
\cline{1-7} 
\multicolumn{1}{|c|}{Mesh size} & \multicolumn{2}{|c|}{$p=2$} &  \multicolumn{2}{|c|}{$p=3$} & \multicolumn{2}{|c|}{$p=4$}\\\hline
\cline{1-7}
$ 0.4981 $&$ 2.75 $&$$&$ 1.16 $&$$&$ 3.09\times 10^{-1} $&$$\\
$ 0.2828 $&$ 1.84 $&$( 0.70 )$&$ 3.61\times 10^{-1} $&$( 2.06 )$&$ 4.33\times 10^{-2} $&$( 3.47 )$\\
$ 0.1627 $&$ 1.10 $&$( 0.94 )$&$ 1.17\times 10^{-1} $&$( 2.03 )$&$ 7.35\times 10^{-3} $&$( 3.21 )$\\
$ 0.0973 $&$ 6.00\times 10^{-1} $&$( 1.17 )$&$ 3.45\times 10^{-2} $&$( 2.39 )$&$ 1.06\times 10^{-3} $&$( 3.76 )$\\
$ 0.0508 $&$ 2.93\times 10^{-1} $&$( 1.10 )$&$ 8.46\times 10^{-3} $&$( 2.16 )$&$ 1.27\times 10^{-4} $&$( 3.27 )$\\
$ 0.0269 $&$ 1.47\times 10^{-1} $&$( 1.08 )$&$ 2.11\times 10^{-3} $&$( 2.18 )$&$ 1.47\times 10^{-5} $&$( 3.38 )$\\
$ 0.0138 $&$ 7.24\times 10^{-2} $&$( 1.06 )$&$ 5.12\times 10^{-4} $&$( 2.11 )$&$ 1.70\times 10^{-6} $&$( 3.22 )$\\\hline
\end{tabular}
\caption{Error values in the $\|\cdot\|_{h,1}$-norm and EOCs for Experiment~\ref{exp1}.}\label{exp1errs:s6}
\end{center}
\end{table}
\begin{table}[H]
\begin{center}
\begin{tabular}{|c|cc|cc|cc|}\hline
\cline{1-7} 
\multicolumn{1}{|c|}{Mesh size} & \multicolumn{2}{|c|}{$p=2$} &  \multicolumn{2}{|c|}{$p=3$} & \multicolumn{2}{|c|}{$p=4$}\\\hline
\cline{1-7}
$ 0.4981 $&$ 1.06\times 10^{-1} $&$$&$ 2.63\times 10^{-2} $&$$&$ 4.82\times 10^{-3} $&$$\\
$ 0.2828 $&$ 7.06\times 10^{-2} $&$( 0.72 )$&$ 1.09\times 10^{-2} $&$( 1.56 )$&$ 9.19\times 10^{-5} $&$( 6.99 )$\\
$ 0.1627 $&$ 3.42\times 10^{-2} $&$( 1.31 )$&$ 2.83\times 10^{-3} $&$( 2.44 )$&$ 5.32\times 10^{-6} $&$( 5.16 )$\\
$ 0.0973 $&$ 1.22\times 10^{-2} $&$( 2.01 )$&$ 4.44\times 10^{-4} $&$( 3.60 )$&$ 3.94\times 10^{-7} $&$( 5.06 )$\\
$ 0.0508 $&$ 3.53\times 10^{-3} $&$( 1.91 )$&$ 7.48\times 10^{-5} $&$( 2.74 )$&$ 3.16\times 10^{-7} $&$( 0.34 )$\\
$ 0.0269 $&$ 9.70\times 10^{-4} $&$( 2.03 )$&$ 8.79\times 10^{-6} $&$( 3.36 )$&$ 2.34\times 10^{-8} $&$( 4.08 )$\\
$ 0.0138 $&$ 2.50\times 10^{-4} $&$( 2.02 )$&$ 7.72\times 10^{-7} $&$( 3.63 )$&$ 1.48\times 10^{-9} $&$( 4.13 )$\\\hline
\end{tabular}
\caption{$\|c-c_h\|_{L^2(\partial\X)}$ error values and EOCs for Experiment~\ref{exp1}.}\label{exp1Cerrs:s6}
\end{center}
\end{table}
\subsubsection{Experiment 2}\label{exp3}
In this experiment, we consider the 
following problem
\begin{equation}\label{exp:3:obl}
\left\{
\begin{aligned}
\sum_{i,j=1}^2(1+\delta_{ij})\frac{x_i}{|x_i|}\frac{x_j}{|x_j|} D^2_{ij}u&= f,\quad\mbox{in}\quad\X,\\
\beta\!\cdot\!\nabla u &\mbox{ is constant on }\partial\X,
\end{aligned}
\right.
\end{equation}
where $\X=\{x\in\mathbb R^2:|x|<1\}$. We take $\beta$ to be the anti-clockwise rotation of the normal by the angle $\varphi(x_1,x_2):=\pi/4+\arctan(\textstyle{\frac{x_2}{x_1}})$, for $(x_1,x_2)\in\partial\X$. 
Furthermore, the function $f$ on the right-hand side of~(\ref{exp:3:obl}) is chosen so that the solution $u$ is given by $u(x_1,x_2) = \frac{1}{4}\cos(\pi(x_1^2+x_2^2))-\frac{1}{\pi}\int_\X\frac{1}{4}\cos(\pi(x_1^2+x_2^2)).$ Notice that in this case the compatibility constant $c=0$, $\gamma=2/5$, $\varepsilon=3/5$, $\Theta = \pi/4+\varphi(x_1,x_2)$, and $\partial_{{\Ta}}\Theta \equiv1$.
%
This experiment serves to demonstrate the robustness of this method with respect to the choice of oblique vector, $\beta$, and choice of discontinuous coefficents, $A\in L^\infty(\X;\mathbb R^{2\times 2})$. In particular, $\beta$ performs a full rotation around the normal vector, and the coefficients $A_{12},A_{21}$, are discontinuous across the lines $\{x\in\X:x_1=0\}$ and  $\{x\in\X:x_2=0\}$.

In this experiment, we successively increase the degree, $p$, of the finite element space $\dgcompo$ from $2$ to $4$, and for each fixed degree we refine the mesh quasi--uniformly. In Tables~\ref{exp3errs:s6} and~\ref{exp3errsH1:s6}, we report the error values in the $\|\cdot\|_{h,1}$-norm and the $|\cdot|_{H^1(\X;\Th)}$-seminorm, respectively, with the corresponding experimental orders of convergence given in brackets. 
We observe the optimal convergence rates $\|(u-u_h,c-c_h)\|_{h,1}=\mathcal{O}(h^{p-1})$, and $|u-u_h|_{H^1(\X;\Th)}=\mathcal{O}(h^{p})$.
\begin{table}[H]
\begin{center}
\begin{tabular}{|c|cc|cc|cc|}\hline
\cline{1-7} 
\multicolumn{1}{|c|}{Mesh size} & \multicolumn{2}{|c|}{$p=2$} &  \multicolumn{2}{|c|}{$p=3$} & \multicolumn{2}{|c|}{$p=4$}\\\hline
\cline{1-7}
$ 0.4981 $&$ 8.64 $&$$&$ 5.05 $&$$&$ 1.13 $&$$\\
$ 0.2828 $&$ 6.86 $&$( 0.41 )$&$ 1.03 $&$( 2.80 )$&$ 2.50\times 10^{-1} $&$( 2.66 )$\\
$ 0.1627 $&$ 3.66 $&$( 1.14 )$&$ 3.10\times 10^{-1} $&$( 2.18 )$&$ 6.72\times 10^{-2} $&$( 2.38 )$\\
$ 0.0973 $&$ 1.86 $&$( 1.32 )$&$ 1.16\times 10^{-1} $&$( 1.91 )$&$ 1.41\times 10^{-2} $&$( 3.04 )$\\
$ 0.0508 $&$ 8.80\times 10^{-1} $&$( 1.15 )$&$ 2.96\times 10^{-2} $&$( 2.10 )$&$ 1.81\times 10^{-3} $&$( 3.15 )$\\
$ 0.0269 $&$ 4.42\times 10^{-1} $&$( 1.08 )$&$ 8.31\times 10^{-3} $&$( 1.99 )$&$ 2.41\times 10^{-4} $&$( 3.17 )$\\
$ 0.0138 $&$ 2.21\times 10^{-1} $&$( 1.04 )$&$ 2.15\times 10^{-3} $&$( 2.02 )$&$ 3.02\times 10^{-5} $&$( 3.10 )$\\\hline
\end{tabular}
\caption{Error values in the $\|\cdot\|_{h,1}$-norm and EOCs for Experiment~\ref{exp3}.}\label{exp3errs:s6}
\end{center}
\end{table}
\begin{table}[H]
\begin{center}
\begin{tabular}{|c|cc|cc|cc|}\hline
\cline{1-7} 
\multicolumn{1}{|c|}{Mesh size} & \multicolumn{2}{|c|}{$p=2$} &  \multicolumn{2}{|c|}{$p=3$} & \multicolumn{2}{|c|}{$p=4$}\\\hline
\cline{1-7}
$ 0.4981 $&$ 1.01 $&$$&$ 4.26\times 10^{-1} $&$$&$ 6.14\times 10^{-2} $&$$\\
$ 0.2828 $&$ 3.32\times 10^{-1} $&$( 1.96 )$&$ 4.33\times 10^{-2} $&$( 4.04 )$&$ 8.78\times 10^{-3} $&$( 3.43 )$\\
$ 0.1627 $&$ 1.13\times 10^{-1} $&$( 1.95 )$&$ 7.71\times 10^{-3} $&$( 3.13 )$&$ 1.42\times 10^{-3} $&$( 3.30 )$\\
$ 0.0973 $&$ 3.34\times 10^{-2} $&$( 2.37 )$&$ 1.55\times 10^{-3} $&$( 3.12 )$&$ 2.15\times 10^{-4} $&$( 3.67 )$\\
$ 0.0508 $&$ 8.45\times 10^{-3} $&$( 2.11 )$&$ 2.11\times 10^{-4} $&$( 3.07 )$&$ 1.10\times 10^{-5} $&$( 4.57 )$\\
$ 0.0269 $&$ 2.14\times 10^{-3} $&$( 2.16 )$&$ 3.57\times 10^{-5} $&$( 2.79 )$&$ 7.45\times 10^{-7} $&$( 4.23 )$\\
$ 0.0138 $&$ 5.32\times 10^{-4} $&$( 2.08 )$&$ 6.10\times 10^{-6} $&$( 2.64 )$&$ 4.96\times 10^{-8} $&$( 4.04 )$\\\hline
\end{tabular}
\caption{Error values in the $|\cdot|_{H^1(\X;\Th)}$-seminorm and EOCs for Experiment~\ref{exp3}.}\label{exp3errsH1:s6}
\end{center}
\end{table}
%
\subsubsection{Experiment 3}\label{exp2}
In this experiment, we consider problem~(\ref{exp:3:obl}),
 where $\X=\{x\in\mathbb R^2:|x|<1\}$, and $\beta$ is a $\pi/4$ anticlockwise rotation of the normal, $n_{\partial\X}$. That is
$\beta = \frac{1}{\sqrt{2}}([n_{\partial\X}]^1-[n_{\partial\X}]^2,[n_{\partial\X}]^1+[n_{\partial\X}]^2)^T.$  In this case, $f$ is chosen so that the solution of~(\ref{exp:3:obl}) is given by $u(x) = |x|^{1.5}-0.75|x|^2-\frac{1}{\pi}\int_\X(|x|^{1.5}-0.75|x|^2).$
Notice that in this case the compatibility constant $c=0$, $\gamma=2/5$, $\varepsilon=3/5$, $\Theta\equiv\pi/4$, and $\partial_{{\Ta}}\Theta \equiv0$.

In this experiment, the true solution $u\in H^2(\X)$, and, in particular, $u\in H^{5/2-\delta}(\X)$ for arbitrary $\delta>0$. However, the $H^s$-broken Sobolev regularity of $u$ fails for $s\ge5/2$, and we must appeal to the conformal regularity estimate of Lemma~\ref{min:reg:lemma}. 
In this experiment we successively increase the degree, $p$, of the finite element space $\dgcompo$ from $2$ to $4$. 

Furthermore, we compute the numerical solution both on sequence of meshes refined towards the origin (where the solution lacks regularity, an example of such a mesh is given in Figure~\ref{adaptmesh:test}), and on a sequence of quasi-uniformly refined meshes (that in particular does not prioritise refinement towards the origin). We plot the error arising in both cases (adapted mesh refinement and non adapted mesh refinement) in the broken $H^2$-seminorm, against the number of DoFs in Figure~\ref{adaptmesh:test}. 
For $p=2,3,4$, we see a reduction in error from the adapted mesh sequence. In particular, for $p=3$ and $p=4$ we see a reduction in the order of error from $\mathcal{O}(\operatorname{ndofs}^{-1/4})$ to $\mathcal{O}(\operatorname{ndofs}^{-1/2})$.
\begin{figure}[h]
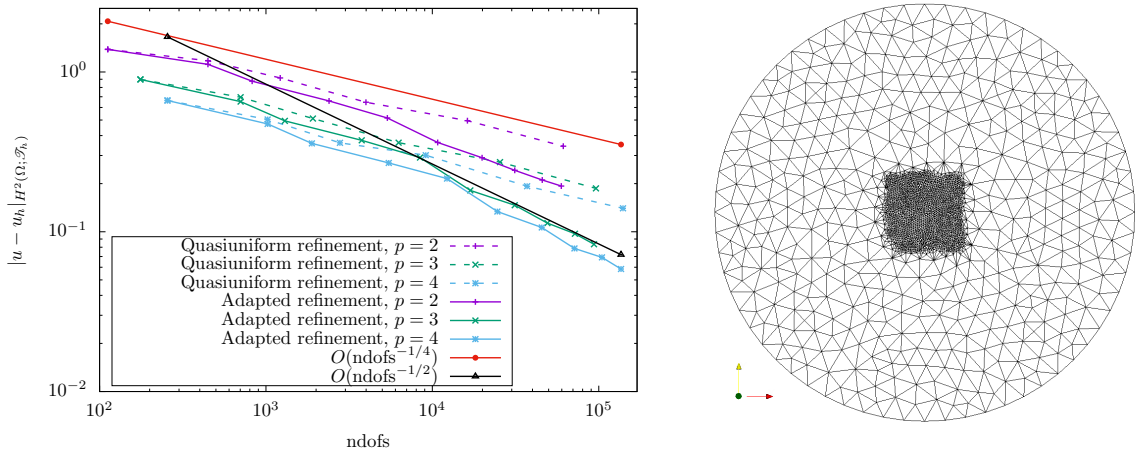

\begin{center}
  \begin{tabular}{lrr}
  \hspace{-5.9cm}
    \scalebox{0.7}{%
    \renewenvironment{table}[1][]{\ignorespaces}{\unskip}%
    \input{paper_exp_nonsmooth.tex}
    \unskip
    }
    &
    \hspace{-0.5cm}
    \scalebox{0.065}{%
    \renewenvironment{table}[1][]{\ignorespaces}{\unskip}%
   \input{diskmesh.tex}
    \unskip
    }
  \end{tabular}
 \end{center}
 \caption{Convergence rates for the numerical scheme applied to problem~(\ref{exp:3:obl}), with a true solution of conformal regularity. On the left, we provide the error values in the $|\cdot|_{H^2(\X)}$ seminorm, where the numerical scheme is implemented on a quasiuniformly refined mesh, and an adapted mesh, with refinement towards the origin. On the right we provide an example of this adapted mesh, at refinement level $7$, consisting of $4532$ elements.
}
 \end{figure}\label{adaptmesh:test}
\end{section}
\begin{section}{Conclusion}\label{Conclusion}
We have extended the framework introduced in~\cite{MR3077903}, and~\cite{Kawecki:77:article:On-curved} allowing for domains with curved boundaries, as well as oblique boundary conditions. In doing so, we have introduced a new DGFEM for elliptic equations in nondivergence form, that satisfy the Cordes condition. 

The computational domain we considered was the unit disc; in order to verify the error estimates
present in Section 3 we used a mesh consisting of curved triangles with edges were defined by
polynomial mappings. It would be an interesting avenue for future research to consider oblique boundary-value problems in dimensions three and higher; this would require one to prove the Miranda--Talenti estimates~(\ref{MT}) in higher dimensions, which is currently an open problem.

The finite element approximation of solutions to elliptic problems in nondivergence form with oblique boundary conditions is a challenging problem, and as such appears to be underrepresented in the available literature. This paper provides and analyses a new method, which appears to be the first discontinuous Galerkin finite element method for oblique boundary-value problems; we were successful in proving both a stability estimate~(\ref{stability:est}), guaranteeing existence and uniqueness of the numerical solution, and an apriori error estimate~(\ref{5.1}) that is optimal with respect to the polynomial degree.


\end{section}
\bibliographystyle{plain}
\bibliography{toskaweckiBIB}
\end{document}